\documentclass[preprint]{imsart}

\RequirePackage[OT1]{fontenc}
\RequirePackage{amsthm,amsmath,amssymb}
\RequirePackage[numbers]{natbib}
\RequirePackage[colorlinks,citecolor=blue,urlcolor=blue]{hyperref}
\usepackage{figsize,subfigure,graphicx,ifthen,calc,psfrag}
\startlocaldefs
\numberwithin{equation}{section}
\theoremstyle{plain}

\endlocaldefs

\newtheorem{lemma}{Lemma}[section]
\newtheorem{proposition}{Proposition}[section]
\newtheorem{theorem}{Theorem}[section]
\newtheorem{remark}{Remark}[section]

\def\R{\mathbb{R}}
\def\Z{\mathbb{Z}}
\def\N{\mathbb{N}}
\def\C{\mathbb{C}}

\begin{document}

\begin{frontmatter}
\title{Sieve Bootstrap for Functional Time Series}
\runtitle{Functional Sieve Bootstrap}
%\thankstext{T1}{Footnote to the title with the ``thankstext'' command.}

\begin{aug}
\author{\fnms{Efstathios} \snm{Paparoditis}\thanksref{t1}\ead[label=e1]{stathisp@ucy.ac.cy}}
%\author{\fnms{Second} \snm{Author}\thanksref{t3,m1,m2}\ead[label=e2]{second@somewhere.com}}
%\and
%\author{\fnms{Third} \snm{Author}\thanksref{t1,m2}
%\ead[label=e3]{third@somewhere.com}
%\ead[label=u1,url]{http://www.foo.com}}

%\thankstext{t1}{Some comment}
\thankstext{t1}{Supported in part   by a University of Cyprus Research Grant.}
%\thankstext{t3}{Second supporter of the project}
\runauthor{Efstathios Paparoditis}

\affiliation{University of Cyprus
%\thanksmark{m1} 
%and Another University\thanksmark{m2}
}

\address{UNIVERISTY OF CYPRUS\\
DEPARTMENT OF MATHEMATICS AND STATISTICS\\
P.O.Box 20537\\
CY-1678 NICOSIA\\
CYPRUS\\
\printead{e1}\\
\phantom{E-mail:\ }%\printead*{e2}
}
%\address{Address of the Third author\\
%Usually a few lines long\\
%Usually a few lines long\\
%\printead{e3}\\
%\printead{u1}}
\end{aug}

\begin{abstract}
A   bootstrap procedure for functional time series is proposed which  exploits   a general  vector autoregressive representation  of 
 the  time series of Fourier coefficients  appearing in  the Karhunen-Lo\`eve expansion of the  functional process. A double sieve-type bootstrap method is developed 
 which  avoids the estimation of process operators and generates functional pseudo-time series that appropriately mimic the dependence structure of the functional  time series at hand. The method  uses   a finite  set of functional principal components to capture  
  the essential driving parts of the  infinite dimensional process   and   a  finite order vector autoregressive process 
  to  imitate  the temporal dependence structure of the corresponding vector  time series of Fourier coefficients.  By allowing the 
  number of functional principal components 
   as well as the autoregressive order used  to increase   to infinity  (at some appropriate rate) as  the sample size increases, 
   consistency of the functional sieve bootstrap can be  established. We  demonstrate this by proving a basic bootstrap  central limit theorem  for  
    functional finite Fourier transforms and  by  establishing bootstrap validity  in the context of  a  fully functional testing problem.
   A novel procedure to select the number of functional principal components is introduced while simulations illustrate the  good  finite sample performance of the new bootstrap method proposed.
\end{abstract}

\begin{keyword}[class=MSC]
\kwd[Primary ]{62M10, 62M15}
%\kwd{60K35}
\kwd[; secondary ]{62G09}
\end{keyword}

\begin{keyword}
\kwd{Bootstrap, Fourier transform, Principal components, Karhunen-Lo\`eve expansion, Spectral density operator}
%\kwd{\LaTeXe}
\end{keyword}

\end{frontmatter}

\section{Introduction}
Statistical inference for  time series stemming  
from stationary  functional processes  has attracted considerable interest  during the last decades 
and  progress has been made in several directions. Estimation and testing procedures  
have been developed 
for  a  wide range of inference problems  and for  large  classes of stationary functional processes; see  Bosq (2000),   H\"ormann and Kokoszka (2012) and Horv\'ath and Kokoszka (2012). 
However,  the asymptotic results derived,   typically   depend  in a complicated way on difficult  to estimate, infinite dimensional  characteristics of the underlying 
 functional  process.  This  restricts considerably the  implementability  
 of  asymptotic approximations   when used  in practice to judge the uncertainty of estimation procedures or to calculate critical values of tests.
In  such situations, bootstrap methods  can provide  useful  alternatives.
%, i.e.,   for  statistical inference problems based on Hilbert space valued random variables,  the bootstrap becomes an  important %tool.

 Bootstrap procedures  for Hilbert space-valued  time series proposed so far  in the literature,  
are mainly attempts to adapt, to the infinite dimensional functional set-up,  of bootstrap methods    that have been  developed 
for the finite dimensional (i.e., mostly univariate)  time series case; cf. Lahiri (2003).   Politis and Romano (1994) considered applications of the  
stationary bootstrap to functional, Hilbert-valued time series and showed  its validity  for 
 the sample mean for functional processes satisfying certain mixing  and boundeness conditions.  Dehling et al.  (2015)  considered applications of the non-overlapping block bootstrap to U-statistics for so called near epoch dependent functional processes and Sharipov et al. (2016) to change point analysis.  Franke  and Nyarigue (2016) and  Zhou and Politis (2016)  developed  some theory for different  residual-based bootstrap procedures  applied to a  first order functional autoregressive process.
 % resp. to a nonparametric functional autoregression.  
Notice that the transmission  of  other bootstrap methods  for real-valued time series to the functional set-up,  like for instance  of the autoregressive-sieve  bootstrap,  Kreiss (1988) and Kreiss et al. (2011), 
 seems  to be difficult  mainly   due to problems  associated with the estimation  (of  an 
 with sample size increasing number) of infinite dimensional autoregressive operators.

%Apart from the aforementioned contributions,
%  where  validity of the bootstrap methods used, have been established, 
% there 
Applications of bootstrap procedures  to certain inference problems  in functional time series analysis  have been also considered in the literature.
% without, however,  developing  the accompanying theory.   
For  instance, for  the construction of prediction intervals,  Fern\'andez De Castro et al. (2005) used  an approach based on  
resampling pairs of functional observations  by means of  kernel-driven resampling probabilities. The same authors also apply a parametric,  residual-based bootstrap approach  using  an estimated first order  functional autoregression  with i.i.d. resampling of appropriately defined functional residuals. For the same prediction problem,
 Hyndman and Shang (2009)  applied different bootstrap approaches  including 
 bootstrapping  the functional curves by randomly disturbing the forecasted scores  using residuals obtained from univariate autoregressive  fits.  Aneiros-Perez et al. (2011) considered the nonparametric functional autoregressive models, while Mingotti et al. (2015) the case of the integrated functional autoregressive model.   Apart from the lack of theoretical justification,   the aforementioned bootstrap applications do not  provide  a general bootstrap  methodology  for functional time series as they are  designed for and their   applicability  
is  restricted  to the  particular inference problem considered; 
% Thus and despite the fact that  for  statistical inference based on Hilbert space valued random variables,  the bootstrap becomes an  %important tool, bootstrapping functional time series seems to be  a  less developed area;
% and the theory of bootstrapping functional time series appears to be  rather incomplete; 
see  also McMurry and Politis (2011) and Shang (2016) for an overview.

In this paper  a general and easy to implement  bootstrap procedure  
for functional time series   is proposed  which generates  bootstrap replicates $ X_1^\ast$, $X_2^\ast$, $ \ldots, X^\ast_n$ 
of a functional time series $X_1, X_2, \ldots, X_n$ and  is  applicable   to a 
large class of stationary functional processes. 
The procedure avoids the explicit estimation of process operators and  exploits some  basic  properties of the  stochastic process of 
Fourier coefficients (scores)  appearing in the well-known Karhunen-Lo\`eve expansion of the functional random variables.
 It is in particular  shown, that under quite general  assumptions,  the   stochastic process of  Fourier coefficients obeys a so-called  vector autoregressive representation and this  representation plays a key role in developing a   bootstrap procedure for   the functional time series at hand.
% which a basic building block .
More specifically,   to capture the essential driving  functional parts of the underlying infinite dimensional  process,  
 the first $m$ functional principal components  are used and 
the corresponding  $m$-dimensional time series of  Fourier coefficients  is   bootstrapped   using a $p$th order vector autoregression 
  fitted to the vector  time series of sample  Fourier coefficients.  In this way, a  
 $m$-dimensional pseudo-time series of Fourier coefficients is  
 generated  which  imitates   the temporal dependence structure  of  
  the  vector time series of  sample Fourier coefficients. 
Using the (truncated) Karhunen-Lo\`eve expansion, these   pseudo-Fourier coefficients  are  then   transformed   to  functional
  bootstrap replicates of the main driving, principal components,  of the observed functional time series.   
Adding  to these replicates    an  appropriately resampled  functional noise,  
  leads  finally to the bootstrapped  functional pseudo-time series $ X_1^\ast, X_2^\ast, \ldots, X_n^\ast$.
   
In a certain sense,  our bootstrap procedure works by using   a  finite rank (i.e., $m$-dimensional) approximation of the infinite dimensional   
structure  of the underlying functional process and a $p$th order vector autoregressive approximation of its  infinite order   temporal dependence structure.
% of the time series of Fourier coefficients. 
To  achieve consistency and to capture appropriately the entire infinite dimensional  structure 
of the functional process,   the number  $m$ of  functional principal components  
used as well as the order $p$ of the  vector autoregression  applied,  are allowed 
%to bootstrap the  $m$-dimensional time series 
%of Fourier coefficients,  both, 
to increase to infinity (at some appropriate rate) as the sample size $n$ increases to infinity. This double sieve property justifies the use of the term ``sieve bootstrap"   for the bootstrap  procedure proposed.  
%allows for  the dimension   $m$  and  for the autoregressive  order $p$ to increase to infinity as the sample size increases to infinity.

We show that under  quite general   conditions, this    bootstrap procedure succeeds in   imitating   correctly   the entire infinite dimensional autocovariance
 structure of the underlying functional process.   Notice that apart from the problem that  instead of the unknown true scores, the  time series of estimated 
  scores  is  used,  the asymptotic analysis of our bootstrap procedure faces additional  challenges    which are  caused  by  the fact that vector autoregressions of  increasing order and of increasing dimension are considered and that 
  the lower bound of the corresponding spectral density matrix approaches zero as the dimension of the vector time series of scores  used, increases
  to infinity.  
%,  where the underlying infinite dimensional 
%spectral density operator  is approximated by an increasing  sequence of finite rank  (spectral density) operators and the infinite order temporal %dependence   structure of the Fourier coefficients by an increasing sequence of  finite orders. 
%This is why  the term functional sieve bootstrap is coined 
%for our bootstrap procedure. 
%This motivates the name functional sieve bootstrap used for our bootstrap procedure. 
 We  demonstrate how the new  bootstrap procedure proposed can be successfully applied to different inference problems in functional time series analysis. In particular, we apply the proposed sieve bootstrap procedure to the problem of estimating the distribution of the functional Fourier transform which is  fundamental in a multitude of applications and  has   attracted  interest in the functional time series literature; 
% and it  includes the 
%sample mean as a special case; 
see  Cerovecki and H\"ormann (2015) for some recent developments. 
In this context, a basic bootstrap central limit theorem  is  established 
 which shows    validity of the functional sieve  bootstrap  for  this  important  class of statistics.
 %,  which includes the sample mean as a special case.  
  Furthermore, we consider applications of the functional sieve bootstrap  in the context of fully 
 functional testing and to the two sample mean problem 
 and  show how  this  bootstrap procedure  can be applied to consistently estimate  the complicated distribution of the test statistic of interest under the null. 
 %Validity 
%for  non-parametric estimators of the spectral density operator is also established. 

Using the time series of Fourier coefficients in the context of functional time series  analysis has been considered  by many authors in a variety of applications. Among others we mention  Hyndman and Shang (2009) who,  for functional autoregressive models and for the sake of prediction, used  univariate autoregressions  fitted to the scalar  time series of scores. In the same context and more related  to the approach proposed in this paper,       
a multivariate approach of prediction 
%for  functional autoregressive processes,  
has been proposed by Aue et al.  (2014) which works by fitting  a vector autoregressive model  to  the multivariate  time series of scores. 

The paper is organized as follows. Section 2 derives some basic properties and discuss the autoregressive  representations of the vector process of Fourier coefficients appearing in the Karhunen-Lo\`eve expansion of the functional process.  Apart from being 
useful for bootstrap purposes, these properties are of interest on their own.
The functional sieve bootstrap procedure proposed is described  in Section 3 where     some properties of the  bootstrap functional pseudo-time series are  also discussed. Asymptotic validity of the new bootstrap  procedure applied to  finite Fourier transforms and to fully functional testing 
%and nonparametric estimators of the spectral density operator 
is  established in Section 4. Section 5 
proposes some novel practical, data driven rules  to choose  the bootstrap parameters  and 
presents some numerical  simulations which  investigate  the finite sample performance of the functional sieve bootstrap. Comparisons   with   three  variants of  block bootstrap methods are also given. 
Technical proofs and auxiliary lemmas are deferred to Section 6.

\section{The Process of   Fourier Coefficients} \label{se.FouCoe}
\subsection{The functional set-up}
We consider a (functional) stochastic process ${\bf X}= \{X_t, t \in \Z\}$ where for each $t$ (interpreted as time), 
$ X_t $ is a random element of the separable Hilbert space ${\mathcal H}:= L^2([0,1],\R)$ with parametrization $ \tau \rightarrow X_t(\tau) \in  \R$ for $ \tau \in [0,1]$.  As usual we denote by $ \langle\cdot,\cdot\rangle$ the inner product in $ {\mathcal H}$ and by $ \| \cdot \|$ the  induced norm  defined for 
 $ x, y \in {\mathcal H}$ as
$\langle x,y\rangle=\Big(\int_{[0,1]}x(t)y(t)dt\Big)^{1/2}$ and $ \|x\|= \big(\langle x,x\rangle\big)^{1/2}$ respectively.  Furthermore, for  matrices $ A $  and $B$ we denote by $ \|A\|_F$ the Frobenius norm, we write $ A\geq B$ or $ B \leq A$ if $ A-B$ is non-negative hermitian while for an operator $ T$,  $ \|T\| $ denotes its operator norm and   $ \|T\|_{HS}$  its Hilbert-Schmidt norm, if $ T$ is  a Hilbert-Schmidt operator.

For the  underlying functional process ${\bf X}$ it is assumed that its dependence structure  satisfies the following assumption.

{\bf Assumption 1} \ $ {\bf X}$ is a purely non deterministic,   $L^4$-${\mathcal M}$ approximable  process.
%which  satisfies $ E\|X_0\|^4 <\infty$. 

The general notion of $ L^p-\mathcal{M}$ approximability refers to stochastic  process  ${\bf X} =\{X_t, t \in \Z\}$ with $X_t$ taking values in ${\mathcal H}$,  $ E\|X_t\|^p <\infty$,  and where the random element  
  $ X_t$ admits the representation $ X_t=f(\varepsilon_t,\varepsilon_{t-1}, \ldots)$. Here  the $\varepsilon_t$'s are i.i.d. random elements in  $ {\mathcal H}$ and $f$ some measurable function $ f :  {\mathcal H}^\infty \rightarrow {\mathcal H}$.   If for  $ 
\{\widetilde{\varepsilon}_t,t\in\Z\}$ an independent copy of $ \{ \varepsilon_t, t\in\Z\}$ and $X_t^{({\mathcal M})}=f(\varepsilon_t,\varepsilon_{t-1}, \ldots, \varepsilon_{t-{\mathcal M}+1}, \widetilde{\varepsilon}_{t-{\mathcal M}},\widetilde{\varepsilon}_{t-{\mathcal M}-1}, \ldots)$,  the condition  
\[  \sum_{k=1}^{\infty} \big(E \|X_k-X_k^{(k)} \|^p\big)^{1/p} < \infty,\]
is satisfied, then ${\bf X}$ is called $ L^p-\mathcal{M}$ approximable.   $L^p-\mathcal{M}$ approximability is a notion of weak dependence which applies to many commonly used functional time series models, like linear functional processes, functional ARCH processes, etc.; see  H\"ormann and Kokoszka (2010) for more details. 

Let  $\mu := EX_0 \in {\mathcal H}$  be the mean of  ${\bf X}$ which  by stationary is independent of $ t$ and for which 
we assume $ \mu=0$ for simplicity. We denote by 
 $C_h$ the autocovariance operator 
$ C_h : {\mathcal H}\rightarrow{\mathcal H}$ at lag $h \in \Z$  defined by $ C_h(\cdot)=E\langle X_t-\mu,\cdot \rangle (X_{t+h}-\mu) $. Associated with the autocovariance operator is the autocovariance function  $c_h: [0,1]\times[0,1] \rightarrow \R$ with $ c_h(\tau,\nu) =E(X_t(\tau)-\mu(\tau))(X_{t+h}(\nu)-\mu(\nu))$, $ \tau,\nu \in [0,1]$, that is, $ C_h$ is an integral operator with kernel function  $c_h$.

% $  \sum_{h} (1+|h|)^r \|C_h\|_S <\infty$ for some $ r \geq 0$.
% %%%% $\la x,y\ra$

Assumption 1 implies  that $ \sum_{h} \|C_h\|_{HS} <\infty$ and that   for every $ \omega \in \R$ the spectral density operator
% $ {\mathcal F}_\omega : L^2([0,1]) \rightarrow L^2([0,1],\C)$ 
\[   {\mathcal F}_\omega(x) = (2\pi)^{-1}\sum_{h\in \Z} C_h(x) e^{-i h \omega},\ \  x \in {\mathcal H}\]
is well defined, continuous in $\omega$,  selfadjoint  and trace class,  H\"ormann et al. (2015); see also Panaretos and Tavakoli  (2013) for  similar properties under different weak dependence  conditions.
% which require   summability of  functional cumulants.
% where 
%also  examples of stationary processes are discussed possessing  a spectral density operator. 
In what follows we will strengthen somehow the assumption on the norm summability of the autocovariance operator  to the following requirement.

{\bf Assumption 2}  \  $  \sum_{h} (1+|h|)^r \|C_h\|_{HS} <\infty$ for some $ r \geq 0$.
% %%%% $\la x,y\ra$

 Furthermore, we  will assume that   the spectral density operator $ {\mathcal F}_\omega$ satisfies the following condition. 

{\bf Assumption 3} \  
For all $\omega \in [0,\pi]$,   
   the operator  $ {\mathcal F}_{\omega} $ is 
of full rank, i.e., $ \mbox{ker}({\mathcal F}_\omega)=0$. 
%The eigenvalues $ \nu_j(\omega)$, $ j=1,2, \ldots $  of $ {\mathcal F}_\omega$ satisfy  $ \nu_j(\omega) \sim j^{-s} $ for some $s > %1 ??????$.
%The spectral operator $ {\mathcal F}_{\lambda} $ satisfies
%  $ \inf_{\lambda\in[0,\pi]}\langle {\mathcal F}_\lambda(x),x \rangle \geq \delta \\ KANN GAR NICHT GEHEN WENN TRACE CLASS, VILLEICHT %NUR  >  0 $ for all $ x \in (L^2[0,1],\C)\setminus\{0\}$.

 For real-valued univariate processes,  $ \mbox{ker}({\mathcal F}_\omega)=0$ is equivalent to the condition 
   that the spectral density is everywhere in $[0,\pi]$ strictly positive   while   for multivariate 
   process 
   to the non-singularity of the spectral density matrix for every frequency $\omega \in [0,\pi]$.  
  Notice  that all eigenvalues $\nu_j(\omega)$, $ j=1,2, \ldots $ of ${\mathcal F}_\omega$ are positive and   that   $ \sum_{j=1}^\infty \nu_j(\omega) < \infty$ by the trace class property of $ {\mathcal F}_\omega$. 
  %The second requirement of Assumption 3 controls the decay  of $ \nu_j(\omega)$  to zero.  
  %which is needed  in the following since we will allow for the number of projections used 
   %to increase to infinity; see  Section XXX
   %Assumption 2  is equivalent to the requirement  that  
  % the spectral operator $ {\mathcal F}_{\lambda} $ is 
%of full rank, i.e., $ \mbox{ker}(F_\lambda)=0$ for all $ \lambda \in [0,\pi]$.

\subsection{Vector autoregressive representation}
Since  
$ C_0=\int_{-\pi}^\pi {\mathcal F}_{\omega} d\omega$, 
 the positivity of $ {\mathcal F}_\omega$  implies 
that the covariance operator $ C_0$ has full rank, that is,  its eigenvalues  $\lambda_j$ satisfy $\lambda_j>0$ for  all $ j \geq 1$.
% and the  corresponding eigenfunctions $ \{v_j,j \geq 1\}$  form an orthonormal basis.
By the symmetry and compacteness of $C_0$, the random element $ X_t$ admits 
 the well known Karhunen-Lo\`eve representation
 % that is  
 %orthonormalized eigenfunctions $ v_j, j\geq 1$ corresponding  to  
\begin{equation} \label{eq.xKL}
X_t=\sum_{j=1}^{\infty} \langle X_t,v_j\rangle  v_j, \ \  t \in \Z,
\end{equation}
where $v_j$, $j=1,2, \ldots $, are  the  orthonormalized eigenfunctions that correspond 
to the eigenvalues $ \lambda_j$, $ j =1,2, \ldots $,  of  $C_0$.
%, where the later are   assumed to be in descending order. 
For $ t \in \Z$, let $ \xi_{j,t} :=\langle X_t,v_j\rangle $,    $ j \geq 1$,
%, denote by    $  \widetilde{M}=\{s: \lambda_s  >0\}$ 
 and consider any subset 
 of indices  $ M=\{j_1, j_2, \ldots, j_m\}  \subset  \N$
 %\{s: \lambda_s  \ \mbox{eigenvalue of $C_0$} \}$  
 with  $ j_1 < j_2<\ldots < j_m$, $m < \infty$. Later on, we will concentrate  on  the specific  set $ M=\{1,2, \ldots, m\}$  which will be %corresponds to   
 the set of the  
 $m$ largest eigenvalues of the covariance operator $ C_0$. 
 
 Consider  now the   
 $m$-dimensional process $ {\bf \xi}^{(M)}= \{ {\bf \xi}^{(M)}_{t}=(\xi_{j_s,t}^{(M)},  s=1,2, \ldots,m )^{\top}, t \in \Z\}$. Observe that 
   $ {\bf \xi}^{(M)}$  is strictly stationary, purely non deterministic and  has mean zero,  i.e., 
  $ E({\bf \xi}^{(M)}_{t})=(\langle EX_t,v_{j_s}\rangle , \ s \in M)=0$. Furthermore, 
 its    autocovariance matrix function $ \Gamma_{\xi ^{(M)}}(h)=
 E({\bf \xi}^{(M)}_{t} {\bf \xi}^{(M)^{T}}_{t+h} )$, $ h \in \Z$,  is given by  $ \Gamma_{\xi^{(M)}}(h)= \big(\langle C_h(v_{j_s}),v_{l_r}\rangle \big)_{s,r=1,2, \ldots, m} $ and satisfies by Assumption 2, 
 \begin{align} \label{eq.gammasum}
 \sum_{h =-\infty}^{\infty}(1+|h|)^r\| \Gamma_{{\xi}^{(M)}}(h)\|_F  &  =  \sum_{h=-\infty}^{\infty}(1+|h|)^r \Big(\sum_{s,r=1}^m  \langle C_h(v_{j_s}),v_{l_r}\rangle^2\Big)^{1/2}  \nonumber \\
 & \leq \sum_{h=-\infty}^{\infty}(1+|h|)^r \|C_h\|_{HS} < \infty.
 \end{align} 
Note that    the bound on the right hand side above is independent of  the set $M$ and  that although by construction  it holds  true that $ Cov(\xi^{(M)}_{r_1,t}, \xi^{(M)}_{r_2,t})=0$ for $ r _1\neq r_2$,  the random variables  $ \xi_{r_1,t} $ and $  \xi_{r_2,s}$  may be correlated for  $t\neq s$. The summability property (\ref{eq.gammasum})   implies that  the 
$m$-dimensional vector  process   $ {\bf \xi}^{(M)}$ possesses a continuous   spectral density matrix  $  f_{{\bf \xi}^{(M)}}(\cdot)$    which is given by 
 \[  f_{{\bf \xi}^{(M)}}(\omega)= (2\pi)^{-1} \sum_{h=-\infty}^{\infty} \Gamma_{\xi^{(M)}}(h) e^{-i\omega h },\ \  \omega \in \R.\]
Moreover,   $  f_{{\bf \xi}^{(M)}} $ satisfies the following boundeness conditions.
 
\begin{lemma} \label{le.eigen} Under Assumption 1 and   3 and Assumption 2 with $r=0$,  the spectral density 
 $ f_{{\bf \xi}^{(M)}}$  satisfies 
 %the following boundeness conditions, 
 \begin{equation}  \label{eq.fbounds} 
 \delta_M  I_m \leq   f_{{\bf \xi}^{(M)}}(\omega) \leq c\, I_{m}, \ \ \mbox{for all $\omega \in [0,\pi]$},
 \end{equation}
  where $ \delta_M$ and $ c$ are real numbers ($\delta_M$ depends on the set $M$),  such that  
  $0< \delta_M \leq  c < \infty$ and $ I_m$ is the $ m\times m$ unity matrix.   
\end{lemma}
\vspace*{0.1cm}

The continuity and the  boundeness properties of the spectral density matrix $ f_{\xi^{(M)}}(\cdot)$ stated in Lemma   \ref{le.eigen},  
imply  
%by   Theorem XXX of  . That is       
that the process    $ {\bf \xi}^{(M)}$ obeys  a so called vector autoregressive representation;  
Cheng and Pourahmadi (1983), see also Wiener and Masani  (1958). That is,    there exist
 an infinite  sequence of $m \times m$-matrices  $ \{A_j^{(M)} , j \in \N\}$ and a  full rank 
 $m$-dimensional white noise process $ \{e^{(M)}_t, t\in \Z\}$, 
such that  $ \xi_t^{(M)} $ can be expressed as 
%$  \xi_t^{(M)}$ can be expressed using his  history as  
\begin{equation} \label{eq.ARxi}
{\bf \xi}^{(M)}_t = \sum_{j=1}^{\infty} A_j^{(M)} \xi^{(M)}_{t-j} + e_t^{(M)}, \ t \in \Z,
\end{equation}
where the coefficients matrices satisfy $ \sum_{j\in \N} (1+j)\|A_j^{(M)}\|_F < \infty$ 
%(where the upper bound is independent of  $M$) 
and $\{ e_t^{(M)}, t \in \Z\}$ is a zero mean white noise  innovation process, that is    $ E(e_t^{(M)})=0$ and $ E(e_t^{(M)} e^{(M)\top}_s) = \delta_{t,s} \Sigma_e^{(M)}$,
with $ \delta_{t,s}=1$ if $ t=s$,  $ \delta_{t,s}=0$ otherwise and $ \Sigma_e^{(M)}$ a full rank $m \times m$ covariance matrix.  
We stress here the fact that       (\ref{eq.ARxi}) 
 does not describe a model for the process of  Fourier coefficients  $ \xi_t^{(M)}$ and should not be confused with the  so-called linear,  infinite order vector autoregressive (VAR($\infty$)) process driven by  independent, identically distributed (i.i.d.) innovations.
 %, which is the case for invertible linear processes. 
  In fact, representation (\ref{eq.ARxi})  is   the  autoregressive analogue of the well-known (moving average) Wold representation
  of  $ \xi_t^{(M)}$ with respect to the same white noise  innovation process $ \{e^{(M)}_t,t\in\Z\}$. This autoregressive  representation  is valid for any stationary and purely non deterministic process the spectral density matrix of which is continuous and  satisfies the boundness conditions (\ref{eq.fbounds});  see also Cheng and Pourahmadi (1983) 
   and Pourahmadi (2001)  for details.   In contrast to the Wold representation, 
 the autoregressive  representation  (\ref{eq.ARxi}) seems to be  more appealing for statistical purposes, since it  express  the vector time series of  Fourier coefficients $ \xi^{(M)}_t $ as a 
 function of  its  
   (in principle) observable past values $ \xi_{t-j}^{(M)}$, $ j=1,2, \ldots $.

In what follows we  assume that   the eigenvalues are in descending order, i.e.,  that 
$ \lambda_1>\lambda_2>  \cdots  >\lambda_m >0$ 
%for every $m \in \N$ 
and  we consider  
the set $ M=\{1,2, \ldots, m\}$  of  the   $m$ largest eigenvalues of $C_0$. The corresponding normalized  
eigenfunctions  (principal components) are denoted by  $v_j$, $ j=1,2, \ldots, m$  and are  (up to a sign) uniquely identified. 
Furthermore,  by Parseval's identity, 
% Recall that 
  %$ E\|X_t-\mu\|^2=\sum_{j=1}^\infty\lambda_j$, that is,  
   the quantity $ \sum_{j=1}^{m}\lambda_j$ describes   
 the variance of $ X_t$  captured  by the first $m$ functional principal components.
 % which  correspond
  %to the first $m$ largest eigenvalues.
To simplify notation we   surpass in the following  the upper index   $(M)$ and  write simple $ \xi_t$ for $ \xi^{(M)}_t $  respectively $ f_\xi$ for $f_{\xi{(M)}} $,  keeping in mind that 
the $ j$th component $ \xi_{j,t}=\langle X_t, v_j\rangle$ of $ \xi_t=(\xi_{1.t}, \xi_{2,t}, \ldots, \xi_{m,t})^{\top} $ is obtained using the orthonormalized eigenfunction $v_j$ which corresponds to the 
$j$th largest eigenvalue $\lambda_j$ of $C_0$,  $j=1,2, \ldots,m$.    Furthermore, we write $A_j(m) $, $e_t(m)$,  $\delta_m$ and $\Sigma_e(m) $ for $ A_j^{(M)}$, $ e_t^{(M)}$, $ \delta_M$    and $\Sigma_e^{(M)} $,  respectively.

\section{The Functional Sieve Bootstrap Procedure}
\subsection{The bootstrap procedure}
The basic idea of  our  procedure is to generate pseudo-replicates $ X_1^\ast, X_2^\ast, \ldots, X_n^\ast$ of the functional time series at hand by first bootstrapping the $m$-dimensional time series of Fourier coefficients $ \xi_t=(\xi_{1,t}, \xi_{2,t}, \ldots, \xi_{m,t})^{\top}$, $t=1,2, \ldots, n$,
corresponding  to the first $m$ principal components. 
%in  the  Karhunen-Lo\`eve expansion of $X_t$.  
This $m$-dimensional time series of Fourier coefficients   is bootstrapped using  the  autoregressive representation of  $\xi_t$  discussed in Section 2.2.  The generated $m$-dimensional pseudo-time series of Fourier coefficients is  then  transformed  to functional principal pseudo-components  by means of  the  
 truncated Karhunen-Lo\`eve expansion $ \sum_{j=1}^m \xi_{j,t}v_j$. Adding to this   an  appropriately resampled  functional noise leads  to the  
%Before giving a precise description of  this bootstrap procedure, we   fix some notation. 
%Let   $\widehat{v}_j$, $ j=1,2 \ldots, m$,  be  the (sign-adjusted) sample analogue of $ v_j$ which,   
 %given a functional time series $ X_1, X_2, \ldots, X_n$  at hand,  are  obtained  as the eigenfunctions of the sample 
 %covariance operator  $ \widehat{C}_0=n^{-1}\sum_{t=1}^{n} (X_t-\overline{X}_n) \otimes (X_t-\overline{X}_n)$ corresponding to %the $m$ largest estimated eigenvalues $ \widehat{\lambda}_j$, $j=1,2, \ldots, m$ of   $\widehat{C}_0 $. Here,   
%$ \overline{X}_n=n^{-1}\sum_{t=1}^{n}X_t$ denotes  the sample mean of the functional time series.
%More precisely, the
%bootstrap procedure proposed to generate the 
functional pseudo-time series  $ X_1^\ast, X_2^\ast, \ldots, X_n^\ast$.  However, since  the $\xi_t$'s are not observed, we work with the time series of estimates scores. This   idea is precisely described in the following functional sieve bootstrap algorithm.  

\begin{itemize}
\item[] {\it Step 1}: \ Select a number $ m=m(n)$ of functional  principal components  and an autoregressive order $p=p(n)$, both finite and depending on $n$. 
\item[] {\it Step 2}: \ Let 
%$ \widehat{\bf \xi}_t=(\widehat{\xi}_{1,t}, \widehat{\xi}_{2,t}, \ldots, \widehat{\xi}_{m,t})^{'} $, $ t=1,2, \ldots, n$ be the time series of %estimated 
%Fourier coefficients, 
\[   \widehat{\bf \xi}_t=(\widehat{\xi}_{j,t}= \langle X_t, \widehat{v}_j\rangle , j=1,2, \ldots, m)^{\top}, \ \  t=1,2, \ldots, n, \]
be the $m$-dimensional time series of estimated Fourier coefficients, where $ \widehat{v}_j$, $j=1,2, \ldots,m$ are the estimated  
eigenfunctions corresponding to the estimated  eigenvalues $\widehat{\lambda}_1 > 
\widehat{\lambda}_2> \cdots > \widehat{\lambda}_m$ of  the sample covariance operator $ \widehat{C}_0=n^{-1}\sum_{t=1}^{n} (X_t-\overline{X}_n) \otimes (X_t-\overline{X}_n)$,    
$ \overline{X}_n=n^{-1}\sum_{t=1}^{n}X_t$.
% is   the sample mean of the functional time series.
\item[] {\it Step 3}: \  Let $ \widehat{X}_{t,m}= \sum_{j=1}^{m}\widehat{\xi}_{j,t} \widehat{v}_j $ and define  the functional 
residuals $ \widehat{U}_{t,m}= X_t - \widehat{X}_{t,m}$, $ t=1,2, \ldots, n$.
\item[] {\it Step 4}: \  Fit a $ p$th order  vector autoregressive process to  the $m$-dimensional time series $ \widehat{\xi}_t$, $ t=1,2, \ldots, n$, denote by $ \widehat{A}_{j,p}(m), j=1,2, \ldots, p$,  the  
%Yule-Walker 
estimates of the autoregressive matrices 
%$ A_{j,p}(m)$, $j=1,2, \ldots, p$, 
and by  
$ \widehat{e}_{t,p}$
%, $ t=p+1, p+2, \ldots, n$  be 
the  residuals, 
\[ \widehat{e}_{t,p} = \widehat{\xi}_t - \sum_{j=1}^p \widehat{A}_{j,p} (m)\widehat{\xi}_{t-j}, \  t=p+1, p+2, \ldots, n.\]
Different estimators $ \widehat{A}_{j,p}(m)$, $j=1,2, \ldots, p$ can be used,  but we focus in the following on  Yule-Walker estimators; cf. Brockwell and Davis (1991). 
\item[] {\it Step 5}: \  Generate a $m$-dimensional pseudo time series of scores $ \xi^\ast_t=(\xi_{1,t}^\ast, \xi_{2,t}^\ast, \ldots, \xi_{m,t}^\ast )$, $ t=1,2,\ldots, n$, using 
\[ \xi_t^\ast= \sum_{j=1}^{p} \widehat{A}_{j,p}(m) \xi^\ast_{t-j} + e^\ast_t, \]
where  $ e^\ast_t$, $ t=1,2, \ldots, n$ are i.i.d. random vectors  having as distribution the empirical distribution of  the centered residual vectors  $ \widetilde{e}_{t,p}=\widehat{e}_{t,p} -\overline{\widehat{e}}_n$, $ t=p+1,p+2, \ldots, n$  and $\overline{\widehat{e}}_n=(n-p)^{-1}\sum_{t=p+1}^{n}\widehat{e}_{t,p} $.  
\item[] {\it Step 6}: \ 
Generate a pseudo-functional time series $ X_1^\ast, X^\ast_2, \ldots, X^\ast_n$,  where  
\begin{equation}  \label{eq.bootX}
 X^\ast_t = \sum_{j=1}^m \xi^\ast_{j,t}\widehat{v}_j + U^\ast_t, \ \ \ \ \ t=1,2, \ldots, n,
 \end{equation}  
and  $ U^\ast_1, U^\ast_2, \ldots, U^\ast_n$ are  i.i.d. random functions obtained by choosing with replacement from the set of centered functional residuals 
$ \widehat{U}_{t,m} -\overline{\widehat{U}}_n$, $ t=1,2, \ldots, n$ and   $ \overline{\widehat{U}}_n=n^{-1}\sum_{t=1}^{n}\widehat{U}_{t,m}$.
 \end{itemize}

Some comments regarding  the above algorithm are in order. 
Notice first   that $ X_1^\ast, X_2^\ast, \ldots, X_n^\ast$ are functional pseudo-random variables and  that the  
 autoregressive representation of  the vector time series of Fourier coefficients  is solely used   as a  tool  to  bootstrap the  $m$ main 
 functional principal components of the  functional  time series at hand.
  In fact, it is this  autoregressive representation 
 which allows the   generation of the  pseudo-time series of  Fourier coefficients $ \xi_1^\ast, \xi^\ast_2, \ldots, \xi_n^\ast$  in Step 4 and Step 5    in a way that  
  imitates the dependence structure of the  sample Fourier coefficients $ \xi_1, \xi_2, \ldots, \xi_n$.
 These pseudo-Fourier coefficients are transformed to bootstrapped main principal components  by means of the truncated and estimated Karhunen-Lo\`eve expansion  which 
%    and which are  used   in the  (m-truncated)  Karhunen-Lo\`eve  representation. This
 together  with the  additive  functional noise  $U^\ast_t$, 
 lead to   the new functional 
pseudo-observations $X^\ast_1, X^\ast_2, \ldots, X_n^\ast$.  

The estimated eigenfunctions  $\widehat{v}_j$ used in Step 2 may point in an opposite direction than the eigenfunctions $ v_j$. In asymptotic derivations  this  is commonly  taking care off    by considering the sign corrected estimator  $ \widehat{s}_j \widehat{v}_j$, where the (unobserved) random  variable $ \widehat{s}_j$ is given by $ \widehat{s}_j = \mbox{sign}(\langle \widehat{v}_j,v_j\rangle )$. 
%This modification  makes  $ \widehat{s}_j \widehat{v}_j $ and $v_j $ pointing in the same direction.  
However, since in our setting adding this sign correction will not affect the asymptotic results derived, we assume  for simplicity throughout this paper, that $ \widehat{s}_j=1$, for $j =1,2, \ldots, m$.  
%Notice that other multivariate bootstrap approaches,  like for instance a variant of the 
% block bootstrap, can be also used to bootstrap the vector time series of sample Fourier coefficients and to generated the  $ \xi_t^+$'s. Second, 
\begin{remark} \label{re.cent} {\rm  To simplify notation we have assumed that  the mean of $ {\bf X}$ is zero. 
If $ EX_t = \mu \neq 0$ then  the  sieve bootstrap algorithm  can   be appropriately modified  by  defining  the pseudo-random element $X_t^\ast$ in 
{\it Step 6}   as $ X_t^\ast = \overline{X}_n + \sum_{j=1}^{m} \xi_{j,t}^\ast \widehat{v}_j + U_t^\ast$,  $ t=1,2,  \ldots, n$.
Notice that since  under Assumption 1,  $ \| \overline{X}_n -\mu\| = O_P(n^{-1/2})$, see H\"ormann and Kokoszka (2012),   the asymptotic results derived in this paper are not affected, i.e.,  $EX_t=0$ is not a stringent assumption.}
% this case is  the mean can by ccurately removed in a preproccecing step    
\end{remark}

\begin{remark} {\rm Modifications  of the above basic bootstrap algorithm are possible   which concern 
 the  resampling schemes used  to  generate the vector of pseudo-innovations $ e^\ast_t$ and/or   the bootstrap functional noise  $U_t^\ast$.  
 To elaborate, and as we will see in the sequel, for general stationary processes satisfying Assumption 1, the  applied  i.i.d. resampling used to generate  the pseudo-innovations $ e^\ast_t$ in {\it Step 5}, suffices in order  to  capture the entire, infinite dimensional second order structure of the underlying functional process $ {\bf X}$. However, a  modification of this i.i.d.  resampling scheme 
 may be needed    if   
 %statistics  of the functional time series are considered, the distribution of which depend
 higher order dependence  characteristics of the underlying functional process beyond those of order two,  should also  be correctly mimicked by the functional pseudo-time series $ X_1^\ast, X_2^\ast, \ldots, X_n^\ast$. In such a 
 case,  the  i.i.d. resampling  used to generate the $ e^\ast_t$'s 
 % and/or  the $ U_t^\ast$'s, 
 in {\it Step  5} 
 %and Step 6 respectively, 
 can be replaced by other resampling  
 schemes  (i.e.,   block bootstrap schemes) 
that are able to  capture    higher order dependence characteristics of the white noise process  $\{ e_t, t \in \Z\}$ appearing in (\ref{eq.ARxi}).}
 % of the "residual" series $ \widehat{e}_{t,p}$ respectively   $ \widehat{U}_{t,m}$. 
%  However,  in the following we restrict ourselves   to the simple  case of i.i.d. resampling of both bootstrap noise series $ e^+_t$ and $ U^\ast_t$, %since this resampling  scheme suffices 
% for proving  validity of the bootstrap procedure proposed when it is applied to statistics the distribution of which depends  only 
% on the autocovariance structure of the underlying functional process. The sample mean as well as  nonparametric estimators of the 
% spectral density operator are examples of such statistics that will be discussed in the next section. 
\end{remark}

\subsection{Some properties of the bootstrap functional process}
As usual, all considerations made regarding the bootstrap procedure are made conditionally  on the observed functional time series $ X_1, X_2, \ldots, X_n$. The generation mechanism of the pseudo-time series $ X_1^\ast, X_2^\ast, \ldots, X_n^\ast$, enables us to consider 
 the bootstrap functional process ${\bf X}^\ast=  \{X_t^\ast, t \in \Z\}$,  where  for $ t\in \Z$, $ X_t^\ast=\sum_{j=1}^m{\bf 1}_j^{\top} \xi_{t}^\ast \widehat{v}_j + U^\ast_t$,  with  
$ \{ \xi^\ast_t=(\xi^\ast_{1,t}, \ldots, \xi^\ast_{m,t})^{\top}, t \in \Z\}$  generated 
as 
%by   the vector autoregressive process  
$ \xi^\ast_t= \sum_{j=1}^p\widehat{A}_{j,p}(m) \xi^\ast_{t-j} + e^\ast_t$ and   the  
$ U^\ast_t $'s  are  i.i.d. functional random variable taking values in the set $\{ \widehat{U}_{t,m} -\overline{\widehat{U}}_n, t=1,2, \ldots, n\}$ with probability $ 1/n$. In the above notation  $ {\bf 1}_j$ is the $m$-dimensional vector 
${\bf 1}_j=(0,\ldots, 0, 1, 0, \ldots, 0)^{\top}$,  where the unity  appears  in the $ j$th position.

It is easy to see that  
$ {\bf X}^\ast$ is a strictly stationary functional process with mean function $ E^\ast X^\ast_t =0 $ and autocovariance operator $ C^\ast_h : {\mathcal H} \rightarrow {\mathcal H}$ given,  for $ h \in \Z$,  by
\begin{align*}
C_h^\ast (\cdot) 
 & = 
\sum_{j_1=1}^{m}\sum_{j_2=1}^m {\bf 1}^{'}_{j_1}{\bf \Gamma}^\ast_h {\bf 1}_{j_2}\langle \widehat{v}_{j_1},\cdot\rangle
 \widehat{v}_{j_2} + I(h=0)E^\ast\langle U^\ast_t, \cdot\rangle U^\ast_{t},
\end{align*}
where $ \Gamma_h^\ast =E^\ast (\xi_{t}^\ast \xi_{t+h}^{\ast^{T}} )$  is the $m\times m$ autocovariance matrix at lag $h$ of $\{  \xi^\ast_t, t \in \Z\}$.
% of the  vector autoregressive  process $ \xi_t^+=\sum_{j=1}^{p}\widehat{A}_{j,p}\xi^+_{t-j} + e^+_t$. 
%Notice that for every $h \in \Z$, 
$ C^\ast_h$ is a Hilbert-Schmidt operator since it is, for $h\neq 0$,   a   finite rank  operator 
%$ \sum_{j_2=1}^{m}\sum_{j_1=1}^m {\bf 1}^{'}_{j_1}{\bf \Gamma}^\ast_h {\bf 1}_{j_2}\langle \widehat{v}_{j_1},\cdot\rangle
% \widehat{v}_{j_2}$ 
while  for $h=0$ it is   the  sum of  a  finite rank operator and of the (Hilbert-Schmidt) empirical covariance operator of the functional pseudo-innovations 
 $  C^\ast_U=E^\ast\langle U^\ast_t, \cdot\rangle U^\ast_{t}=n^{-1}\sum_{t=1}^{n} \langle \widehat{U}_{t,m} -\overline{\widehat{U}}_n,\cdot  \rangle ( \widehat{U}_{t,m} -\overline{\widehat{U}}_n)$.
 
%Furthermore, since  for every $m\in \N$, 
%   \[  \sum_{h\in \Z} \|C^\ast_h\|_{HS} \leq   \sum_{h \in \Z} \|\Gamma^\ast_h\| _F +  I(h=0) \|C^\ast_U\|_{HS}=O_P(1),\]
% the bootstrap process   ${\bf X}^\ast$ possesses  for $ \omega \in \R$  the spectral density operator $ {\mathcal F}^\ast_{ \omega,m} : L^2([0,1],\C) \rightarrow L^2([0,1],\C)$ given  by 
%   \begin{equation} \label{eq.bootsdo}
%   {\mathcal F}^\ast_{ \omega,m}(x) = (2\pi)^{-1} \sum_{h \in \Z}C^\ast_h(x) e^{-ih \omega},\ \ \ \  x \in {\mathcal H}.
%   \end{equation}
 
   If the (estimated) 
 vector autoregressive process used   to generate the time series  of pseudo-scores $ \xi_t^\ast$  is stable, then the 
  dependence structure  of  the bootstrap process $ {\bf X}^\ast$  can be  precisely described. This is stated 
 in the following proposition. Notice that the required stability condition of   the  estimated autoregressive polynomial is fulfilled,  if  for instance,   $ \widehat{A}_{j,p}$, $j=1, 2, \ldots., p$, are the Yule-Walker estimators; cf. Brockwell and Davis (1991),  Ch. 11.4.

\begin{proposition} \label{pr.Lp_m}  If  $ p, m \in \N$  is such that 
%for  $m\in\N$ and $ p \in \N$ 
the estimator $ \widehat{A}_{j,p}$, $j=1,2, \ldots, p$, used in Step 4 of the functional  sieve bootstrap algorithm  is well defined and satisfies  
$ det(\widehat{A}_{p,m}(z)) \neq 0$ for  all $ |z| \leq 1$,  where $\widehat{A}_{p,m}(z) = I _m- \sum_{j=1}^p \widehat{A}_{j,p}(m) z^j $, $ z \in \C$, then,  conditionally on $ X_1, X_2, \ldots, X_n$,   the  bootstrap process  $ {\bf X}^\ast$  is   $ L^2-{\mathcal M}$ approximable.
\end{proposition}

 The $ L^2-{\mathcal M}$ approximability of $ {\bf X}^\ast$ implies that  $ \sum_{h} \|C^\ast_h\|_{HS} < \infty$, see  H\"ormann et al. (2015), which can be also easily verified  since   
   \[  \sum_{h\in \Z} \|C^\ast_h\|_{HS} \leq   \sum_{h \in \Z} \|\Gamma^\ast_h\| _F +  I(h=0) \|C^\ast_U\|_{HS}=O_P(1).\]
Furthermore, and because of the $ L^2-{\mathcal M}$ approximability property,  the bootstrap process   ${\bf X}^\ast$ possesses  for every $ \omega \in \R$  a spectral density operator $ {\mathcal F}^\ast_{ \omega,m} $
% : L^2([0,1],\C) \rightarrow L^2([0,1],\C)$ 
defined   by 
   \begin{equation} \label{eq.bootsdo}
   {\mathcal F}^\ast_{ \omega,m}(x) = (2\pi)^{-1} \sum_{h \in \Z}C^\ast_h(x) e^{-ih \omega},\ \ \ \  x \in {\mathcal H}.
   \end{equation}
 $ C_h^\ast$ and $ {\mathcal F}^\ast_{\omega,m}$ are essentially finite rank  approximations of the corresponding population operators  $C_h$ and $ {\mathcal F}_{\omega}$ respectively. Thus  and in order for the bootstrap  process ${\bf X}^\ast$ to capture the infinite dimensional structure of the underlying functional process and  the infinite order dependence structure of the vector time series generating the scores, the dimension $m$ as well as  the autoregressive order $p$, used in the functional sieve bootstrap algorithm,  have to increase to infinity (at some appropriate rate)  as the sample size $n$ increases to infinity.  This  rate
 should 
      take into account  the fact that the true scores  and eigenfunctions appearing  in the Karhunen-Lo\`eve
    expansion  are not observed and, therefore,   sample estimates are used instead. Furthermore,     
    % not only the  order $p$ but also the dimension $m$ of the fitted autoregressive process has to increase to infinity with  the sample size which makes the asymptotic analysis quite  involved.  Finally,  
     the lower bound $ \delta_m$ 
    of the  spectral density matrix   of the scores   $ f_\xi $, 
    %corresponding vector autoregressive process, 
    approaches zero as  the sample size $ n $ increases to infinity. This  is due to  the fact that the eigenvalues $ \nu_j(\omega)$ of the spectral density operator $ {\mathcal F}_\omega$ converge to  zero as $j \rightarrow \infty$. These facts make the asymptotic analysis quite involved  and impose several restrictions regarding the  
    behavior of  $m$ and $p$  with respect to the sample size $n$ which  are   summarized   in the following assumption.

{\bf Assumption 4}\   The sequences  $m=m(n)$ and $ p=p(n)$ satisfy 
$m \rightarrow \infty$ and $p \rightarrow \infty$ as $ n \rightarrow \infty$ such that,  
\begin{enumerate}
\item[(i)] \ $ m^{3/2}=O(p^{1/2}) $ 
\item[(ii)] \  ${\displaystyle  \frac{p^7}{n^{1/2}\lambda^2_m }\sqrt{ \sum_{j=1}^{m}\frac{1}{\alpha_j^{2}}} \rightarrow 0 }$, 
where  $ \alpha_1=\lambda_1-\lambda_2$ and $ \alpha_j=  \min\{\lambda_{j-1}
-\lambda_j , \lambda_j-\lambda_{j+1}\}$ for $ j=2,3, ..., m$. 
\item[(iii)] \ $ \delta_m^{-1} \sum_{j=p+1}^{\infty} j^r \| A_{j}(m)\|_F \rightarrow  0$ for some   $ r \geq 0$,  where  $ \delta_m$ is the lower bound of the spectral density matrix $ f_\xi$  given in  (\ref{eq.fbounds}).
\item[(iv)] \  $m^4p^2\|\widetilde{A}_{p,m} - A_{p,m}\|_F
%\sum_{j=1}^{p}\|\widetilde{A}_{j,p}(m) - A_{j,p}(m)\|_F  
= O_P(1)$, where $\widetilde{A}_{p,m}=(\widetilde{A}_{1,p}(m),  \ldots,  \widetilde{A}_{p,p}(m))$, 
 and $ A_{p,m} =(A_{1,p}(m),  \ldots,  A_{p,p}(m))$. Here,  $ \widetilde{A}_{j,p} $, $j=1,2, \ldots, p$ denotes  the same estimator as $ \widehat{A}_{j,p}$, $j=1,2, \ldots, p$,  based on the true vector time series of scores $ \xi_1, \xi_2, \ldots, \xi_n$ instead of their estimates $ \widehat{\xi}_1, \widehat{\xi}_2, \ldots, \widehat{\xi}_n$ and  $ A_{j,p}(m) $, $j=1,2, \ldots, m$ are the coefficient matrices  of the  best  (in the mean square sense) linear predictor of $\xi_t$ based on $ \xi_{t-j}$, $j=1,2, \ldots, p$.  
\end{enumerate}
 
Assumption 4(i) restricts the rate with which the dimension $m$ is allowed to increase to infinity compared with that of $p$. 
Assumption 4(ii)  is   imposed  in order to control the error made by the fact that  the bootstrap procedure  is based on estimated scores and eigenfunctions  instead on the unobserved true quantities in a context where the dimension $m$ and the autoregressive order $p$, both,  increase  to infinity and the lower bound of the spectral density matrix of the $m$-dimensional vector of  scores approaches zero as $ m$ increases to infinity. Part (iii)  relates the rate of increase of the  autoregressive order $p$ to the lower bound of the spectral density matrix $ f_\xi$ and the decay of the norm of the autoregressive matrices to zero. Part (iv) is  essentially  a requirement on the rate at which $m$ and $p$ are allowed to increase to infinity taking into account the convergence rate of the estimator   $ \widetilde{A}_{j,p}$, $j=1,2, \ldots, p$  based on the true scores.  For instance,  calculations similar to that in the proof of Lemma~\ref{le.A3} yield for the Yule-Walker estimator that $\|\widetilde{A}_{p,m} - A_{p,m}\|_F=O_P(mp n^{-1/2}(\sqrt{m} \lambda_m^{-1} + p)^2) $ which, taking into account Assumption 4(i),  implies that Assumption 4(iv) is satisfied if $ m, p \rightarrow \infty$ slowly enough with $n$  such that 
$ m p^6 =O( \sqrt{n}\lambda_m^2) $ and $ p\lambda_m^2 =O(m^2)$.
Notice  that,  for real valued-random variables, such assumptions relating the rate of increase of the autoregressive parameters to the convergence rate of the estimators used,  are common in the autoregressive-sieve bootstrap literature; see  Kreiss et al. (2011) and Meyer and Kreiss (2015).  However, 
 the situation here is  much more involved  since  in our context,  not only the order $p$ but  also the dimension  $m$ of the vector autoregression  has to increase  to infinite with the sample size by  taking into account the fact that  $ \lambda_m$ converges to zero as $m$ increases to infinity. 
 
The following lemma illustrates the rate conditions imposed in Assumption 4 by considering  two particular examples of  the behavior of the difference $ \lambda_j-\lambda_{j+1}$ which is  related to the rate of decrease of the eigenvalues $\lambda_j$. According to this lemma,  $ p$ may increase to infinity  as $ n^a$ for some  $a>0$ while the rate of increase of $m$ depends on  the rate of decrease of $ \lambda_j-\lambda_{j+1}$ respectively of the eigenvalues $\lambda_j$, $j=1,2,\ldots$. If these differences 
 decrease with  a geometric rate, then $m$ may increase at most logarithmically  in the sample size $n$, while if the same differences  decrease with a polynomial rate, then $m$ may increase to infinity faster, like $ n^\zeta$ for some appropriate $\zeta>0$.  

\begin{lemma} \label{le.Ass4}
Assume that  $ \widetilde{A}_{p,m}$ are the Yule-Walker estimators of $A_{p,m} $.
\begin{enumerate}
\item[(i)] \ If $\lambda_j-\lambda_{j+1}  \geq C_\lambda \rho ^{j}$ for $ j=1,2, \ldots $, $ \rho \in (0,1)$ and $ C_\lambda >0$, then  Assumption 4(i), (ii) and (iv)  is satisfied if  
\[ p=O(n^a) \   \ \mbox{and}\ \ m \leq \Big( \frac{1}{6\log(\rho^{-1})}\big(1-14a) - \delta\Big) log(n),\]
%\ \  m = O( log(n)), \]
for  $  a \in (0,1/14)$ and some $\delta >0$.
%and 
%\[ m_n \leq \Big( \frac{1}{6\log(\rho^{-1})}\big(1-14a) - \delta\Big) log(n),\]
%\[ m_n = O( log(n)).\]
%for some $\delta >0$.
\item[(ii)]\  If $\lambda_j-\lambda_{j+1}  \geq C_\lambda j^{-\theta} $  for $ j=1,2, \ldots $ and for some $ \theta >1$ and $ C_\lambda >0$, then Assumption 4(i), (ii) and (iv)  is satisfied if 
\[ p=O(n^a)\  \ \ \mbox{and}\ \  m = O(n^\zeta), \]
%and 
%\[ m_n=O(n^\zeta),\]
for  $  a \in (0,1/14)$ and $ \zeta \in [\zeta_{\min}, \zeta_{\max}]$, where $ \zeta_{\min} = a/(2+2\theta)$ and $ \zeta_{\max} =
 \min\{(1-14a)/(1+6\theta)-\delta, a/3\}$ for some $ \delta >0$.
\end{enumerate}
\end{lemma}

Under the condition that  $m$ and $ p$ increase to infinity at an  appropriate rate with $n$ such that Assumption 4 is satisfied, the following proposition can be established which shows that  the spectral density operator  ${\mathcal F}_{\omega,m}^\ast$
of the bootstrap process $ {\bf X}^\ast$ converges, in  Hilbert-Schmidt norm, to  the spectral density operator  ${\mathcal F}_{\omega}$ of the underlying functional  process $ {\bf X}$. 

\begin{proposition} \label{pr.f} Under Assumptions 1 and 3 and  Assumption 2  and 4 with $r=2$, we have, that, as $ n \rightarrow \infty$, 
\[\sup_{\omega\in [0,\pi]}\|{\mathcal F}^\ast_{ \omega,m} - {\mathcal F}_ \omega\|_{HS} \ \rightarrow \ 0,\]
in probability. 
\end{proposition}

From   the above proposition and  the inversion formulae of Fourier transforms, we immediately  get  for the covariance operators $C^\ast_h$ and $ C_h$ of the bootstrap process ${\bf X}^\ast$  and of the underlying process $ {\bf X}$,  that  $ \sup_{h\in \Z} \|C^\ast_h - C_h\|_{HS} \rightarrow 0$, in probability,    as $n \rightarrow \infty$. Thus  the bootstrap process $ {\bf X}^\ast$,  imitates asymptotically correct  the entire infinite dimensional autocovariance structure of the functional process $ {\bf X}$. This allows  for  the use of the bootstrap functional time $X_1^\ast, X_2^\ast, \ldots, X_n^\ast$  to approximate the distribution of statistics based on the functional time series $ X_1, X_2, \ldots, X_n$. Some examples of such statistics  are discussed in the next section. 

So far we have assumed that the covariance operator $C_0$ has full rank, i.e., that its eigenvalues  $\lambda_j$ are distinct which implies that, for  consistency and in order to capture the entire infinite dimensional dependence structure of the underlying functional process ${\bf X}$, the number $m$ of principal components included, has to increase to infinity with the sample size $n$. The situation is much simpler if we assume that $m_0\in \N$ exists such that $ \lambda_{m_0} >0$ and $\lambda_j=0$ for all $ j >m_0$. In this case only the  finite number of $m_0$ score time series are needed to describe the entire dependence structure of $ {\bf X}$. We are  then essentially in the finite dimensional case with   the $m_0$-dimensional score process $\{ \xi_t=(\langle X_t, v_j\rangle, j=1,\ldots, m_0)^\top, t \in \Z\}$,   possessing a spectral   density matrix which is bounded from bellow by a positive constant independent  of the sample size   $n$. Furthermore, as in the proof of Lemma~\ref{le.A3} and, because in this case  $ \sum_{j=1}^{m_0}\|\widehat{v}_j-v_j\|^2 = O_P(n^{-1/2})$, we get that   $ \| \widehat{A}_{p,m_0}-\widetilde{A}_{p,m_0}\|_F=O_P(p^4/\sqrt{n})$. Standard arguments  applied in  the case of the (finite dimensional) vector autoregressive-sieve bootstrap  can then be used  (see for instance Meyer and Kreiss (2015)), to show that under less restrictive conditions that those stated in Assumption 4, $ \sup_{\omega \in [0,\pi]}\|{\mathcal F}_{\omega,m_0}^\ast - {\mathcal F}_\omega\|_{HS} \stackrel{P}{\rightarrow} 0,$ in probability. 
% This suggests that 
%the functional sieve bootstrap procedure proposed, can be successfully applied to approximate the distribution of  statistics the  distribution of which 
%depends on the second order  structure  of the underlying functional process.  
%Examples of such statistics 
% will be considered in the next section.

\section{Bootstrap Validity}
In this section we investigate the  validity of the functional sieve  bootstrap  applied in order to approximate the distribution of some statistic          
 $T_n=T(X_1, X_2, \ldots, X_n)$ 
of interest,  when 
 the bootstrap analogue  $ T^\ast_n=T(X_1^\ast, X_2^\ast, \ldots, X_n^\ast)$  is used.   Notice that establishing validity of  a bootstrap procedure for time series heavily depends on two issues; see also  Kreiss and Paparoditis (2011). On the dependence structure of the underlying process which  affects the distribution of the statistic of interest and on the capability of the bootstrap procedure used to mimic appropriately this dependence structure. Furthermore, since  proving bootstrap  validity is a case by case matter, we  demonstrate in the following  applications  of the functional sieve bootstrap procedure proposed to some  statistics that have  recently  attracted considerable interest in the functional time series literature.   
     
 \subsection{Functional finite Fourier transform}
 Consider the distribution of 
  the functional Fourier transform  
 \begin{equation} \label{eq.fFT}
  S_n(\omega) = \sum_{t=1}^{n}X_t e^{-i t \omega}, \ \ \ \omega\in [-\pi, \pi].
  \end{equation}
Notice that  the sample mean  $ \overline{X}_n=n^{-1}S_n(0)$ is  just a  special case of  (\ref{eq.fFT}).  In order to elaborate 
on  the limiting distribution of $ S_n(\omega)$ we first  fix some notation. We say that a  random element  $Z \in  {\mathcal H}_\C:={\mathcal H} + i {\mathcal  H}$,  follows a circularly-symmetric complex Gaussian distribution with mean zero  and covariance $ {\mathcal G}$,  we write $ Z \sim CN (0, {\mathcal G}) $,  if 
\[  \left(\begin{array}{cc} Re(Z) \\ Im(Z) \end{array} \right) \sim N_{{\mathcal H} \times {\mathcal H}}\Big(  
\left(\begin{array}{cc} 0 \\ 0 \end{array} \right),  \frac{1}{2}
\left(\begin{array}{cc} Re({\mathcal G}) & - Im({\mathcal G})  \\  Im({\mathcal G}) & Re({\mathcal G}) \end{array} \right)\Big); \] 
see also Cerovecki and H\"ormann (2017) for a general discussion of the complex Gaussian distribution.
%We then assume that the sequence $\{n^{-1/2}S_n(\omega), n \in \N\}$ satisfies the following assumption.

Under a  range of different weak dependence assumptions on the functional process ${\bf X}$, 
 it has been shown  that  
 \begin{equation} \label{eq.snnor}
 n^{-1/2} S_n(\omega) \Rightarrow CN(0, 2\pi {\mathcal F}_\omega)
 \end{equation}
  as $ n \rightarrow \infty$, where $ \Rightarrow$ denotes weak convergence on ${\mathcal H}_{\C}$.
%  where for 
%  if $ is a Gaussian process in $ L^2$ with 
% $ EG=0$ and covariance operator $C_G$, which is  an  integral operator with kernel
% \begin{equation}  \label{eq.kernG} 
%  c_{G}(\tau_1,\tau_2) = E[X_0(\tau_1)X_0(\tau_2)] +\sum_{h=1}^{\infty}E[X_0(\tau_1)X_h(\tau_2)]+\sum_{h=1}^{\infty}E[X_0(\tau_2)X_h(\tau_1)],
% \end{equation} 
% for all $\tau_1,\tau_2 \in [0,1]$.
For $\omega=0$, such a limiting behavior   has been established  
for   linear functional   processes by  Merlev\`ede et al. (1997) and  for $ L^p- {\mathcal M}$ approximable processes  by Horv\`ath et al. (2013).
%and for functional processes   ${\bf X}$ fulfilling  a  variety of weak dependence conditions. To elaborate,  for the case of the sample mean $(\omega=0)$,   established such result for Hilbert-valued linear  while  Horv\'ath et al. (2013) for $ L^2-{\mathcal M}$ approximable processes.
 Panaretros and Tavakoli  (2013) derived the above  limiting distribution  of  
 $n^{-1/2}S_n(\omega)$ for $\omega \in [0,\pi]$,   under  a summability   condition of the functional  cumulants, while      
  more general results for  the same statistic and under weaker conditions,  have been recently obtained  by Cerovecki and H\"ormann (2017). 
  
 We propose to use the bootstrap statistic $ n^{-1/2}S^\ast_n(\omega)=n^{-1/2}\sum_{t=1}^{n} X_t^\ast e^{-i t \omega}$ in order to approximate the distribution of the statistic $ n^{-1/2} S_n(\omega)$.   The following 
   theorem establishes asymptotic validity of  this  functional sieve bootstrap proposal   for the class of functional Fourier transforms
   considered. In this theorem,  $d$ is any metric metrizing  weak convergence on $ {\mathcal H}_{\C}$.
   %Mallows $d_2$ metric is used to metrisize the distance between the distribution of two  random elements in  ${\mathcal H}_\C$. %Notice that convergence in the sense of this metric implies weak convergence and convergence of second moments; see Bickel %and Freedman (1981) for details.  We denote by  $ {\mathcal L}(X)  $   the law of the random element $X$ in $ {\mathcal H}_{\C}$ %and  by $ {\mathcal L}(X|Y)  $ the conditional law of $X$ given  $ Y$.
   %,  when the bootstrap statistic  $ n^{-1/2}S^\ast_n(\omega)=n^{-1/2} \sum_{t=1}^n X_t^\ast e^{-i t \omega}$  is used  to %approximate the distribution of  the statistic $ n^{-1/2}S_n(\omega)$.

%
%The following theorem establishes   the  asymptotic limiting distribution of the bootstrap finite Fourier transform  $ S_n^\ast = \sum_{t=1}^{n} X_t^\ast e^{-i t \omega}$.   
 
 \begin{theorem} \label{th.th_main}  Suppose that  for  $\omega  \in [0, \pi]$, the sequence $ \{n^{-1/2}S_n(\omega), n \in \N\}$ in $ {\mathcal H}_{\C}$ satisfies (\ref{eq.snnor}). Suppose further that  Assumptions 1 and  3  and Assumption 2  and 4 with $r=2$ are satisfied.  Then,   as $n \rightarrow \infty$,
\begin{enumerate}
\item[(i)] \  $d\big({\mathcal L}(n^{-1/2}S_n(\omega)), {\mathcal L}(n^{-1/2}S^\ast_n(\omega)|X_1,X_2, \ldots, X_n) ) \ \rightarrow 0,$   \ and
\item[(ii)] \ $  \|n^{-1}E^\ast S^\ast_n(\omega)\otimes S_n^\ast(\omega) - n^{-1}E S_n (\omega)\otimes S_n(\omega)\|_{HS} \stackrel{P}{\rightarrow} 0$,\\ 
\end{enumerate}
 in probability.
 %, where $ d$ is any metric metrizing weak convergence  on $ {\mathcal H}_{\C}$. 
% \[ 
% \[d_K\big({\mathcal L}(n^{-1/2}S_n^\ast(\omega)),   {\mathcal L}(n^{-1/2}S_n(\omega)) \big) \rightarrow  0,\]
% in probability, where $ {\mathcal L}(Y)$ denotes the law of the $ {\mathcal H}$-valued random element $Y$ and $ d_K$ Kolmogorov's distance.
 \end{theorem}
 
 \begin{remark} {\rm Notice that as a special case of the above theorem we get that,   under the assumptions made,  and as $n \rightarrow \infty$,
 $ \sqrt{n}\overline{X}^\ast_n  \Rightarrow   N\big(0, \sum_{h \in Z} C_h\big)$,
 in probability and $ n E^\ast \overline{X}^\ast_n \otimes  \overline{X}^\ast_n \stackrel{P}{\rightarrow}  2\pi {\mathcal F}_{0}$,   which provides one of the first instances of a central limit theorem for the bootstrap for functional time series under 
 the weak dependence conditions stated  in Assumption 1. 
%  To the best of our knowledge, the only  comparable  result is 
% that  in  Politis and Romano (1994), Theorem 3.1,  where  weak convergence   of   the stationary bootstrap  sample mean
% under strong mixing assumptions and  boundeness conditions has been established.
} 
 \end{remark}  
 
 \subsection{Fully functional testing} \label{sec.test} \  In a variety of functional testing situations one is faced with the problem that the limiting distribution under the null of a fully functional test statistic, depends, in a complicated way,  on difficult to estimate characteristics of the underlying functional process. This makes   the practical implementation of  asymptotic results  derived in order to calculate  critical values  of tests a difficult task. To overcome this problem, a common approach in the literature  is to consider tests based on finite dimensional projections. However, such tests have  non-degenerated power only for  alternatives  which are not orthogonal to the space captured by the particular projections considered; see Horv\'ath et. al (2013) and Horv\'ath et al. (2014) for examples.  Using  as an example the two sample mean problem, we demonstrate in the following  how   the sieve bootstrap procedure proposed in this paper,  can be successfully applied  to approximate the null distribution of a fully functional test.
 
 Let ${\bf X}=\{X_t,t\in \Z\}$ and $ {\bf Y}=\{Y_t,t\in \Z\} $ be two independent,  strictly stationary  functional processes with mean functions $ \mu_X=EX_t $ and $ \mu_Y=EY_t$ respectively and consider the testing problem $ H_0: \mu_X=\mu_Y$  against the alternative $ H_1: \mu_X \neq \mu_Y$. Given two time series $ X_1, X_2, \ldots, X_{n_1}$ and $ Y_1, Y_2, \ldots, Y_{n_2} $  stemming from ${\bf X}$ and $ {\bf Y}$ respectively,  a natural test statistic for these hypotheses  is given by
  \[ U_{n_1,n_2} = \frac{n_1 n_2}{n_1+n_2} \|\overline{X}_{n_1}-\overline{Y}_{n_2}\|^2, \]
where $ \overline{X}_{n_1}=n_1^{-1}\sum_{t=1}^{n_1} X_t$ and $\overline{Y}_{n_2}=n_2^{-1}\sum_{t=1}^{n_2}Y_t$.
% are the corresponding sample means.  
If  both  processes satisfy Assumption 1 and $n_1,n_2 \rightarrow \infty$ such that  $ n_1/(n_1+n_2) \rightarrow \theta \in (0,1)$, it has been shown in Horv\'ath et al. (2013),  that $ U_{n_1,n_2} \stackrel{d}{\rightarrow} \int_0^1 \Gamma^2(\tau)d\tau$, where $ \{ \Gamma(\tau), \tau \in [0,1]\}$ is a mean zero Gaussian process with covariance function $ E(\Gamma(\tau_1)\Gamma(\tau_2)) = (1-\theta) c_X(\tau_1,\tau_2) + \theta c_Y(\tau_1,\tau_2)$ for $ \tau_1, \tau_2 \in [0,1]$ and    $ c_X(\tau_1,\tau_2) = Cov(X_0(\tau_1),X_0(\tau_2)) + \sum_{h\geq 1}  $  $ Cov(X_0(\tau_1), X_h(\tau_2)) + \sum_{h\geq 1} $ $ Cov(X_0(\tau_2), X_h(\tau_1))$ and $ c_Y(\tau_1,\tau_2) =  Cov(Y_0(\tau_1), $ $ Y_0(\tau_2)) + \sum_{h\geq 1} Cov(Y_0(\tau_1),$ $  Y_h(\tau_2))  + \sum_{h\geq 1} $  $Cov(Y_0(\tau_2), Y_h(\tau_1))$.  Notice that  the kernel functions $ c_X$ and $ c_Y$ are unknown, which makes the  calculation of  critical values of the test $U_{n_1,n_2}$  a  difficult task. 

Since the functional sieve bootstrap procedure proposed satisfactory imitates the autocovariance structure of the underlying processes, 
it can be successfully  applied  to estimate the critical values of the test $U_{n_1,n_2}$. To elaborate, the goal is to generate two independent functional pseudo-time series $ X^\ast_1, X^\ast_2, \ldots, X^\ast_{n_1}$ and 
$ Y^\ast_1, Y^\ast_2, \ldots, Y^\ast_{n_2}$,  that mimic the autocovariance structure of the processes ${\bf X}$ and $ {\bf Y}$ respectively and  satisfy, at the same time,  the null hypothesis of interest. For this  let $ X_t^\ast$ and $ Y^\ast_t$  be  generated by means of  equation (\ref{eq.bootX})  of the functional sieve bootstrap algorithm,  where for the generation of the $ X_t^\ast$'s   the sample scores $ \widehat{\xi}_t^{(X)}=(\widehat{\xi}^{(X)}_{j,t}= \langle X_t, \widehat{v}^{(X)}_j\rangle, j=1,2, \ldots, m_1)^{\top}$, $ t=1,2, \ldots, n_1$ and  for the generation of the $ Y_t^\ast$'s, the sample scores  $ \widehat{\xi}_t^{(Y)}=(\widehat{\xi}^{(Y)}_{j,t}= \langle Y_t, \widehat{v}^{(Y)}_j\rangle, j=1,2, \ldots, m_2)^{\top}$, $ t=1,2, \ldots, n_2$ are  used in {\it Step 1} of this algorithm. Here $ \widehat{v}^{(X)}_j$, $j=1, \ldots, m_1$ and $ \widehat{v}^{(Y)}_j$, $j=1, \ldots, m_2$,  denote the orthonormalized eigenfunctions of the $m_1$ respectively $m_2$ largest eigenvalues of the sample covariance operators $ \widehat{C}_{0}^{(X)}=n_1^{-1}\sum_{t=1}^{n_1} (X_t-\overline{X}_{n_1}) \otimes (X_t-\overline{X}_{n_1})$ and $ \widehat{C}_{0}^{(Y)}=n_2^{-1}\sum_{t=1}^{n_2} (Y_t-\overline{Y}_{n_2}) \otimes (Y_t-\overline{Y}_{n_2})$ respectively. Notice that generation of   $ X_t^\ast$ and $Y_t^\ast$  by  using  (\ref{eq.bootX})   ensures that    $E^\ast X_t^\ast=E^\ast Y^\ast_t=0$, that is the generated functional pseudo-time series $X^\ast_1, X^\ast_2, \ldots, X^\ast_{n_1} $ and $Y^\ast_1, Y^\ast_2, \ldots, Y^\ast_{n_2} $ satisfy the null hypothesis $H_0$.  Now, let $ \overline{X}^\ast_{n_1}=n_1^{-1}\sum_{t=1}^{n_1}X^\ast_t$ and $ \overline{Y}^\ast_{n_2}=n_2^{-1}\sum_{t=1}^{n_2}Y^\ast_t$   and define the bootstrap analogue of $ U_{n_1,n_2}$ as
 \[ U^\ast_{n_1,n_2} = \frac{n_1n_2}{n_1+n_2}  \|\overline{X}^\ast_{n_1}-\overline{Y}^\ast_{n_2}\|^2. \]
The following theorem establishes validity of the  sieve bootstrap applied to the  functional testing problem considered.

\begin{theorem} \label{th.testing}
Let the conditions of Theorem~\ref{th.th_main} be satisfied and assume that $ n_1,n_2\rightarrow \infty$ such that $ n_1/(n_1+n_2) \rightarrow \theta \in (0,1)$. Then,
\[ \sup_{x\in \R} \big| P(U_{n_1,n_2} \leq x) - P(U^\ast_{n_1,n_2} \leq x| {\bf X}_{n_1}, {\bf Y}_{n_2}) \big|  \rightarrow 0,\]
in probability, where $P(U^\ast_{n_1,n_2} \leq \cdot | {\bf X}_{n_1}, {\bf Y}_{n_2}) $ denotes the distribution function of  \, $ U^\ast_{n_1,n_2}$ conditional on $ {\bf X}_{n_1}=(X_1, X_2, \ldots, X_{n_1}) $ and ${\bf Y}_{n_2}=(Y_1,Y_2, $ $ \ldots, Y_{n_2})$.
\end{theorem}
%
% Theorem~\ref{th.testing} justifies  the  use of percentage points of the  distribution of $ U_{n_1,n_2}^\ast$ in order to obtain bootstrap critical values of  the test $U_{n_1,n_2}$. Furthermore, if $ H_1$ is true, that is if $ \|\mu_X-\mu_Y\| >0$  and $ U_{n_1,n_2} \stackrel{p}{\rightarrow} \infty$ as $ n_1,n_2 \rightarrow \infty$, see for instance  Theorem 4 of Horv\'ath et al. (2013),  then  the  consistency 
% of the  test $U_{n_1,n_2}$ based  on sieve bootstrap estimated critical values,  follows.    

% The assumption of finite fourth moments seems inavoidable since autoregressive parameter estimates for the Fou
% rier coefficients ...
% We stress here the fact that 
% the assumptions imposed on the underlying process ${\bf X}$ for our bootstrap procedure to work are quite minimal, i.e.,  it is only required  that ${\bf X}$  has finite fourth moments and possesses a continuous and positive finite spectral density operator $ {\mathcal F}_\omega$.
  
%\subsection{Long-run variance}

\section{Choice of Parameters  and Numerical Results}
\subsection{Choice  of the sieve  bootstrap parameters}  Implementation of the functional sieve bootstrap requires the choice of two tuning parameters: the  order $p$  and  the dimension $m$. By choosing  these parameters, the problem of overfitting caused  
 by  selecting   a large dimension and/or a high  order  vector autoregressive model,
 should be   seriously taken into account.  
%In this section we discuss some  proposals  on how to choose these parameters in practice,  where 
%we focus   on simple data driven   rules  which lead to  
%%choices of the bootstrap parameters of interest that result in  
%parsimonious  models 
%that   work well in practice and especially in small sample situations.

Several approaches for selecting the number of principal components in functional data analysis have been proposed in the literature; see among others  Yao et al. (2005)  and Li et al. (2013) for the use of  information type criteria. For our purpose, one useful and simple  criterion for selecting the dimension $m$ 
%and  which is commonly used in the functional set-up, 
 is  based  on the ratio of the  total variance explained by the number $m$ of principal components   included,   to the  variance of $X_t$.  
According to this rule,   $m$ is selected  as the smallest positive integer for which  the empirical variance ratio  ($VR_n$) satisfies $ VR_n(m) = \sum_{j=1}^m\widehat{\lambda}_j/ \sum_{j= 1}^n\widehat{\lambda}_j \geq Q$, with $Q$  a predetermined value  %representing the desired proportion  of variability  that should be explained by the number of principal components selected 
and   $ Q=0.80$ or  $ Q=0.85$ two   common choices; 
cf. H\'orvath and Kokoszka (2012).  One drawback  of   the VR-rule  applied to   functional time series,  is 
 that this criterion  does not take into account dependence.
 %,  that is,   the lost on information  on the dependence structure  of ${\bf X}$  
 %of the underlying process 
 %caused by  the restriction to the $m$-dimensional space.  

To overcome  this drawback we  introduce in the following a  generalized variance ratio criterion.  Measuring the total variability of the underlying functional process ${\bf ×}$  by the quantity $ \int_{(-\pi,\pi]}\|{\mathcal F}_\omega\|^2_{HS} d\omega$,   yields by straightforward calculations and evaluating the Hilbert-Schmidt norm using the orthonormal basis  $ \{v_j, j=1,2, \ldots\}$, the expression 
\[   \int_{(-\pi,\pi]}\|{\mathcal F}_\omega\|^2_{HS} d\omega=  
\sum_{l=1}^\infty \sum_{r=1}^\infty \int_{(-\pi,\pi]}\big|f_{\xi_l,\xi_r}(\omega)\big|^2d\omega,\]
where $f_{\xi_l,\xi_r} $  denotes the cross spectral density of the score processes $ \{\xi_{l,t} \}$ and $ \{\xi_{r,t} \}$. Define  next  a  functional process $ {\bf X}^+_m=\{X^+_t, \in \Z\}$,    where $ X^+_t=X^+_{t,m} + U^+_{t,m}$, 
 $ X^+_{t,m} = \sum_{j=1}^{m} \xi_{j,t} v_j$,  $ U^+_{t,m}= \sum_{j=m+1}^\infty \zeta_{j,t}v_j$ and $ \{\zeta_{j,t},t\in \Z\}, j=m+1, m+2, \ldots,  $ are  independent,  i.i.d. processes  which are  independent from $ X_{t,m}^+$ and have mean zero and  $ Var(\zeta_{j,t}) = \lambda_j$.   Observe  that  for any $m$ fixed and ignoring estimation  errors, 
it is the dependence structure of  ${\bf X}_m^+$   which is  essentially mimicked by  the functional sieve bootstrap process ${\bf X}^\ast$. This is so since  in the bootstrap world,  $ U_{t,m} = X_t - \sum_{ j=1}^m\xi_{j,t} v_j$ is treated as an i.i.d. process  and   the  (possible)  correlation between the processes $ \{X_{t,m} =\sum_{j=1}^m \xi_{j,t} v_j\}$ and $\{U_{t,m}\}$ is ignored.
% for  any fixed $m$ after setting $p=\infty$ and replacing the estimated parameters by their true values. 
Let $ {\mathcal F}^+_{\omega,m}$  be the spectral density operator of $ {\bf X}^+_m$.
Using the same measure of  total variability as  for the process ${\bf X}$, we get 
\[  \int_{(-\pi,\pi]}\|{\mathcal F}^+_{\omega,m}\|^2_{HS} d\omega = \sum_{l=1}^m \sum_{r=1}^m \int_{(-\pi,\pi]}\big|f_{\xi_l,\xi_r}(\omega)\big|^2d\omega  + (2\pi)^{-1}\sum_{l=m+1}^\infty \lambda_l^2.\]
Notice that  the term $ (2\pi)^{-1}\sum_{l=m+1}^\infty \lambda_l^2$ is  due to integrating the  squared Hilbert-Schmidt norm of the spectral density operator of the process $ \{U^+_{t,m}\}$. This process  is included in the definition of $ {\bf X}^+_m$ because of  the functional i.i.d.  innovations  $ U^\ast_t$ used in {\it Step 6} of  the sieve bootstrap algorithm to generate the $ X^\ast_t$'s.

The ratio 
\[ GVR(m) =  \int_{(-\pi,\pi]}\|{\mathcal F}^+_{\omega,m}\|^2_{HS} d\omega\Big/ \int_{(-\pi,\pi]}\|{\mathcal F}_\omega\|^2_{HS} d\omega,\]
can then be considered as the proportion of  total variability of the process $ {\bf X}$ captured by that of  the process 
$ {\bf X}^+_m$. Recall that   $ \big|f_{\xi_l,\xi_r}(\omega)\big|^2=\kappa_{l,r}^2(\omega) f_{\xi_l,\xi_l}(\omega)f_{\xi_r,\xi_r}(\omega) $ with $ \kappa_{l,r}$ the squared coherency between the score processes $\{\xi_{l,t}\} $ and $ \{\xi_{r,t}\}$. That is,  $GVR$ explicitly takes into account the entire autocovariance structure of the   processes ${\bf X} $ and ${\bf X}^+_m $.   $GVR(m)$  can then be interpreted as a measure of  the  los on information on the dependence structure of ${\bf X}$  caused by the functional sieve bootstrap procedure based on  $m$ principal components. Note that  if ${\bf X}$ is a white noise process, then  $ GRV(m)=1$ for every value of $m$. In this case we  set  $m=0$ as the most parsimonious choice, i.e.,  no vector autoregression is fitted,  which implies that the functional sieve bootstrap (correctly) reduces to  an i.i.d. bootstrap. 
%applied to the observed  functional, white noise time series $ X_1,X_2, \ldots, X_n$. 
% Finally, we mention   that if the Hilbert-Schmidt norm is replaced by the trace norm of the spectral density operators involved and the additional term $(2\pi)^{-1}\sum_{l=m+1}^\infty \lambda_l^2$  is ignored, then the GVR ratio   reduces to the $VR$ ratio given by   $VR(m)=\sum_{l=1}^m\lambda_l/\sum_{l=1}^{\infty}\lambda_l$.  
%It is worth mentioning here that by ignoring this term and replacing  $ |f_{\xi_l,\xi_r}(\omega)|^2$  by its squared root, the  following VR criterion for dependent functional observations is obtained,
%\[  DVR(m) =\frac{\displaystyle \sum_{l=1}^m \sum_{r=1}^m \int_{(-\pi,\pi]}\big|f_{\xi_l,\xi_r}(\omega)\big| d\omega}{\displaystyle \sum_{l=1}^\infty \sum_{r=1}^\infty \int_{(-\pi,\pi]}\big|f_{\xi_l,\xi_r}(\omega)\big|d\omega} \]
%In particular, if $ {\bf X}$ is an i.i.d. process then, as  it is easily seen,  DVR(m)=VR(m) for all values of $m$.

Now, observe that $ \lambda_j$,  $ \int_{(-\pi,\pi]}\big|f_{\xi_l,\xi_r}(\omega)\big|^2d\omega$  and $ \int_{(-\pi,\pi]}\|{\mathcal F}\omega\|^2_{HS}d\omega $ can be consistently estimated by $\widehat{\lambda}_j $,  $2\pi n^{-1}\sum_{j\in {F_n}} | I_{\xi_l,\xi_r}(\omega_j) |^2$ and $  2\pi n^{-1}\sum_{j\in F_n} $ $ \| I_{n,\omega_j}\|^2_{HS} $,  respectively, where 
$ I_{\xi_l,\xi_r}(\omega)= J_{\xi_l}(\omega)J_{\xi_r}(-\omega)$ and $ J_{\xi_s}(\omega)= $  $ (2\pi n)^{-1/2} $ $ \sum_{t=1}^{n}\xi_{s,t}$ $ e^{-i \omega t}$ for any $ s \geq 1$. Furthermore,   $ I_{n,\omega} $ is the periodogram operator with kernel $I_{n,\omega}(\tau_1,\tau_2)=J_{n,\omega}(\tau_1)\overline{J}_{n,\omega}(\tau_2)$, $J_{n,\omega}(\tau)= (2\pi n)^{-1/2} \sum_{t=1}^{n} X_t(\tau)e^{-i \omega t}$, $ \omega_j=2\pi j/n$, $F_n =\{-N,\ldots, -1,1,\ldots, N\}$ and $ N=[n/2]$. This suggests to select the  dimension  $m$  as the smallest  positive integer for which the   empirical  generalized variance ratio ($GVR_n$)  satisfies
\[ GVR_n(m) =\frac{\displaystyle \sum_{l=1}^m \sum_{r=1}^m \frac{\displaystyle 2\pi}{\displaystyle n}\sum_{j\in F_n} \big|\widehat{I}_{\xi_l,\xi_r}(\omega_j)\big|^2   + \frac{1}{2\pi}\sum_{l=m+1}^n \widehat{\lambda}_l^2 }{\displaystyle  \frac{ \displaystyle 2\pi}{\displaystyle n}\sum_{j\in F_n} \big\| I_{n,\omega_j}\big\|^2_{HS}} \ \geq \ Q.\]
 Here $ \widehat{I}_{\xi_l,\xi_r}(\omega)= \widehat{J}_{\xi_l}(\omega)\widehat{J}_{\xi_r}(-\omega)$ with  $ \widehat{J}_{\xi_s}(\omega)=(2\pi n)^{-1/2}\sum_{t=1}^{n}\widehat{\xi}_{s,t}e^{-i \omega t}$   the finite Fourier transform of  the time series of estimated scores. 

\begin{remark}  {\rm 
$GVR_n$ has been  developed  for the functional sieve  bootstrap situation considered in this paper.  However, a simple modification of this criterion leads to  an alternative  to the $VR_n$ rule   which is appropriate for  dependent functional data and   which is of interest on its own. In particular,   ignoring the second term of  the nominator of $ GVR_n$,  the following dependent variance ratio  ($DVR_n$) criterion is obtained,
\[ DVR_n(m)=  \sum_{l=1}^m \sum_{r=1}^m  \sum_{j\in F_n} \big|\widehat{I}_{\xi_l,\xi_r}(\omega_j)\big|^2\Big/ \sum_{j\in F_n} \big\| I_{n,\omega_j}\big\|_{HS}^2.\]
%\[ DVR_n=  \sum_{l=1}^m \sum_{r=1}^m \frac{2\pi}{n}\sum_{j=-N}^{N} \big|\widehat{I}_{\xi_l,\xi_r}(\omega_j)\big|^2\Big/ \sum_{l=1}^K \sum_{r=1}^K \frac{2\pi}{n}\sum_{j=-N}^{N} \big|\widehat{I}_{\xi_l,\xi_r}(\omega_j)\big|^2 \ \geq \ Q.\]
 $DVR_n$  delivers  an empirical   measure  of the lost on information on the dependence structure of ${\bf X}$ associated with the use of the $m$-dimensional space and can be therefore,  used as a simple criterion  to select the number  $m$ of principal components   in a  functional time series setting.  Notice that  if the Hilbert-Schmidt norm in GVR is replaced by the trace norm of the spectral density operators involved and the  additional term $(2\pi)^{-1}\sum_{l=m+1}^\infty \lambda_l^2$  is ignored, then the corresponding $DVR(m)$  ratio given by $$ DVR(m)= \sum_{l=1}^m \sum_{r=1}^m \int_{-\pi}^{\pi}\big|f_{\xi_l,\xi_r}(\omega)\big|^2d\omega \Big/ \sum_{l=1}^\infty \sum_{r=1}^\infty \int_{-\pi}^{\pi}\big|f_{\xi_l,\xi_r}(\omega)\big|^2d\omega,$$ reduces   to the   $VR(m)=\sum_{l=1}^m\lambda_l/\sum_{l=1}^{\infty}\lambda_l$ ratio.  
}
\end{remark}

  Notice that both,  the VR and the GVR criterion,  refer to  a fixed sample size $n$ and the purpose is to select the number of principal components in a way which ensures that a desired fraction  $Q$ of the variance of the process is  captured by the number of principal components included in the analysis. This is  important for  our bootstrap proposal where the objective 
  is  to appropriately mimic   the dependence structure of the  functional time series  at hand.   However,   consistency requires that $m$  increases to infinity with $n$ which is not the case   if $ Q$ remains fixed with $n$.  At  the same time and  as we have seen, the  rate at which $m$  has to increase to infinity  should  take into account    the rate of decrease of the eigenvalues $\lambda_j$ respectively  of the  differences $ \lambda_j-\lambda_{j+1}$ to zero.  One way to accommodate  such  aspects  in our practical selection of $m$,  is to combine the discussed VR respectively  GVR criterion with  an approach  for selecting   $m$ proposed  by  H\"ormann and Kidzi\'nski (2015) and  which  explicitly takes into account the behavior of  the eigenvalues  $\widehat{\lambda}_j$.   To elaborate,   denote by $ m_{n,E}$ the number of principal components selected   by the rule 
  \[ m_{n,E}={\rm argmax}\Big\{ j \geq 1: \frac{\widehat{\lambda}_1}{\widehat{\lambda}_j} \leq \sqrt{n}/\log(n) \Big\}.\]
Notice that $m_{n,E}$ allows for   the $j$-th  principal component to be included in the analysis if the corresponding estimated eigenvalue $\widehat{\lambda}_j$  is big enough, i.e., if the ratio $ 1/\widehat{\lambda}_j$     does not exceed the threshold $ \sqrt{n}/\log(n)$.   The  nominator  $\widehat{\lambda}_1$ acts solely as   a normalization   to adapt for scaling; for this and for   the choice the  particular threshold  see   H\"ormann and Kidzi\'nski (2015).   Denote now by    $ m_{n,Q}$  the number of principal components selected using, the VR  or  the GVR criterion for some given $Q$. 
 The practical  selection  of $m$ we then propose is  to set  this parameter equal to   
\[  \widehat{m}_n =  \max\{ m_{n,Q}, m_{n,E}\}.\]
According to  this proposal,  only those principal directions  are included in the analysis the eigenvalues of which can be estimated with a reasonably accuracy  ensuring at the same time that the number of principal components selected  explains at least a desired  portion of the variability  of the time series at hand. 
%As our numerical results in Section~\ref{sec.sim} and in  the Supplementary file show, for small to moderate sample sizes,  $ \widehat{m}%_n$ is dominated by $m_{n,Q}$ while as $ n $ increases $ m_{n,E}$ takes over allowing for $m$ to increase to infinity with $n$ , i.e.,  for  %the proportion of variance explained  to converge to unity.   
We  remark that although   for functional time series the   GVR criterion is theoretically more appealing,  for  short  time series of  $n \leq 100$ observations, we  still recommend the use  the VR-criterion 
since   it leads to selections of $m$ with  a smaller variability avoiding, therefore,  the potential  fit of vector autoregressions of large dimensions and/or of high orders which is an important issue for  small samples sizes; see also  Section~\ref{sec.sim} for details.
%and numerical results on using the approach to select $m$ discussed so far.

Once  the dimension $m$ has been selected, the order $p$ of the vector autoregression fitted can be chosen using  the AICC criterion; see  Hurvich and Tsai (1993). This criterion is preferred  because it is  based on an approximately  unbiased   estimator of the expected Kullback-Leibler information of the fitted model and, more importantly,  avoids overfitting.  The order    $p$ is then selected by minimizing    
%\[ AICC(p) = n\log|\widehat{\Sigma}_{e,p} | + \frac{n(nm+pm^2)}{n-m(p+1)-1} ,\]
$ AICC(p) = n\log|\widehat{\Sigma}_{e,p} | +n(nm+pm^2)\big/(n-m(p+1)-1)$,
over a range of possible values of $p$, where $ \widehat{\Sigma}_{e,p}=n^{-1}\sum_{t=p+1}^n\widehat{e}_{t,p}\widehat{e}^T_{t,p}$ and $\widehat{e}_{t,p}$  is defined in {\it Step 4} of the functional sieve bootstrap algorithm.

\subsection{Simulations} \label{sec.sim}
To investigate  the finite sample behavior  of the functional sieve bootstrap (FSP)  we have performed  simulations using time series stemming from a first order functional moving average  process  given by    
\begin{equation} \label{eq.ma1sim}
X_t = \varepsilon_t + \Theta(\varepsilon_{t-1}).
\end{equation}
 as in Aue et al. (2015). To elaborate, $ \Theta$ is specified as  $ \Theta=0.8 \Psi$, where  $\Psi$ is a  linear operator, $ \Psi:{\mathcal H}_D \rightarrow {\mathcal H}_D$, $ {\mathcal H}_D =\overline{sp}\{f_1, f_2, \ldots, f_D\}$,  $D=21$ and  $ f_j$, $j=1,2, \ldots, D$ are  Fourier basis functions on the interval $[0,1]$.    Notice that for 
$ x \in {\mathcal H}_D$, $ x=\sum_{j=1}^D c_j f_j$ with  $ c_j =\langle x, f_j \rangle$, 
the operator $\Psi$  acts  as  
%\sum_{j=1}^Dc_j \Psi(f_j)=
$ \Psi(x)=\sum_{j=1}^D\sum_{l=1}^D c_j \langle \Psi(f_j),f_l\rangle f_l = (B_\Psi  c)^{'} v$,  where $ c=(c_1, \ldots, c_D)^{'}$ and $ v=(f_1,\ldots, f_D){'}$ and the matrix $ B_\Psi$ has element in the $j$th column and $l$th row given  by $\langle \Psi(f_j),f_l \rangle$.
Following  Aue et al. (2015) the operator   $\Psi$ was chosen at random.  For this a    $ D\times D$ matrix of independent, normal random variables 
with mean zero was  first generated   where   its  $ (j_1,j_2)$th element has standard deviation  $  \sigma_{j_1,j_2} =j_1^{-1}j_{2}^{-1}$.  This matrix was then scaled so that the resulting matrix $B_\Psi$ has induced norm equal to 1 and in  every iteration of the simulation runs $B_\Psi$ was newly generated.  The  corresponding i.i.d. innovations $\varepsilon_t$ in (\ref{eq.ma1sim}) were generated   as $ \varepsilon_t =\sum_{j=1}^D Z_{t,j} f_j$,  where  $ Z_{t,j}$ are i.i.d. Gaussian  with mean zero and standard deviation  equal to $ j^{-1}$.

We first consider the performance of the     VR and GVR criteria  in selecting the number $m$ of principal components, when  $Q=0.85$. Table 1 shows the frequencies of selected  dimensions $m$ over $R=1000$ replications of the considered FMA(1) model for  different sample sizes. As it is seen from this table, the VR criterion is quite stable over the different sample sizes leading to the  selections  $ m=3$ or $m=4$ in almost all situations. The GVR criterion exhibits a greater variability for small sample sizes (n$\leq$ 100) and becomes more concentrated around the dimensions $m=4$ and $m=5$ as  $n$ increases.  Observe that  because the GVR criterion explicitly takes into account the dependence structure of the processes involved, it selects more frequently the  larger dimension   $m=4$  compared to the dimension  $m=3$ which is  more frequently selected by  the VR criterion.  Notice  further that the smaller variability of the VR rule for small sample sizes, prohibits  the selection of vector autoregressions of  large dimension which is particularly important in our set-up. Thus  for $n\leq 100$ observations we recommend to apply  the $\widehat{m}_n$ rule using the VR criterion to calculate $ m_{n,Q}$ and  the AICC criterion in order to select the values of $m$ and $p$.

\begin{table*}
\caption{Frequency  of selected values of $m$ by the VR and the GVR criterion ($ R=1000$ replications).} 
\begin{tabular}{llllllllll}
\hline
%& & &   &  & & &  & & \\
 &&m =&  1 & 2 & 3& 4& 5& 6& 7\\
\hline
n=100 & $ VR_n    $  & & 0 & 0.3 & 67.1 & 32.6 & 0 & 0 & 0 \\
           &  $GVR_n $  & & 0 & 0.9 & 19.6 & 55.2 & 22.8 & 1.5 & 0\\
%           &  & &  &&  &  &  & & &  & \\ 
 n=200 &   $VR_n $   &  & 0 & 0  & 62.7 & 37.3 & 0 & 0 & 0 \\
           &  $GVR_n  $&  & 0 & 0.1 & 10.4 & 68.7 & 20.6 & 0.2 &  0 \\
%           &  & &  &  && &  &  & & &  & \\   
  n=300 &  $ VR_n $   &   & 0 & 0  & 66.2 & 33.8 & 0 & 0 & 0 \\
           &  $GVR_n $ &  & 0 & 0 & 4.3 & 75.4 & 20.3 & 0 &  0 \\ 
  n=500 &  $ VR_n $   &   & 0 & 0  & 64.7 & 35.3 & 0 & 0 & 0 \\
           &  $GVR_n $ &  & 0 & 0 & 0.9 & 83.0 & 16.1 & 0 &  0 \\ 
  n=1000 &  $ VR_n $   &   & 0 & 0  & 64.1 & 35.9 & 0 & 0 & 0 \\
           &  $GVR_n $ &  & 0 & 0 & 0.2 & 89.3 & 10.5 & 0 &  0 \\          
           \hline                 
\end{tabular}
\end{table*}

To  investigated the  behavior of $ \widehat{m}_n$ for  the FMA(1) model considered,  we use  a range of sample sizes with   $ m_{n,Q}$  chosen according to  the $VR$ $(n\leq 100)$ respectively $GVR$ criterion with $ Q=0.85$. Table 1 of   the supplementary material    shows the results obtained over $R=1000$ repetitions for each of the sample sizes considered. As it is seen from this table,   the behavior of   $ \widehat{m}_n$  is dominated for small to moderate sample sizes by $ m_{n,Q}$  ensuring, therefore,  the 
desired description of the variability   of the  functional time series by the number $m$ of principal components selected.  However, as  $n$ increases the behavior of $ \widehat{m}_n$ 
becomes  dominated by $ m_{n,E}$  which allows  for  the number of principal components selected   as well as for the part of the variance explained,  to increase with $n$.  
%In particualr, as it is seen from this table, for small to moderate sample sizes, the number of principal components selected is dominated by %$ m_{Q,n}$. Furthermore,  as $ n$  increases, the number of principal components selected  increases too due to the fact that the %behavior of $ \widehat{m}_n$ is dominated by that of $ m_{E,n}$ XXXXXX  

We next consider the behavior  of the FSB procedure in  estimating the standard deviation  of the sample  mean  $\sqrt{n}\overline{X}_n(\tau_j)= n^{-1/2}\sum_{t=1}^n X_t(\tau_j)$, calculated for time series of length $n=100$  observations  and for  $ \tau_j$, $j=1,2, \ldots, T$, $T=21$, equidistant time points in the interval $[0,1]$.  The exact standard deviation of the sample mean is 
estimated using 20,000 replications of the moving average model (\ref{eq.ma1sim}).  
All estimates presented  are based on $R=1,000$ replications and $ B=1,000$ bootstrap repetitions. 
Table 2 of the supplementary material  shows the FSB estimates obtained using  some different values of the bootstrap parameters  $m$ and $p$ as well as for the values of  these parameters chosen by means of the $\widehat{m}_n$ and $AICC$  rule and which are denoted by $ (\widehat{m},\widehat{p})$.
Note that  $(m,p)=(3,3)$ is the most frequently chosen pair using this data driven selection rule.  As this table  shows the FSB estimates  are  quite good even for  the short  functional time series  of $n=100$ observations. These estimates  also   
seem   not to be very sensitive with respect to the different choices of the parameter  $m$ used to truncate  the Karhunen-Lo\'eve expansion.  
%Observe  that as $m$ increases, the variability  of the FSB estimates  increases too. This is expected since   in this case the %number of parameters to be estimated increases too. Recall that  the number of autoregressive coefficients  of   the vector %autoregression fitted to the time series of scores equals $m^2 p$. 

Table 2  compares the   results using the FSB procedure  with those of three different block bootstrap
methods, the moving block bootstrap (MBB), the tapered block bootstrap (TBB) and the stationary bootstrap (SB).  To asses the overall behavior of the different bootstrap estimates, we use the averaged  absolute bias (ABias),   $ T^{-1}\sum_{j=1}^T|\overline{\sigma}^\ast(\tau_j)-\sigma(\tau_j)|$, the averaged relative  bias (RBias),   $ T^{-1}\sum_{j=1}^T|\overline{\sigma}^\ast(\tau_j)/\sigma(\tau_j)-1|$ and the averaged standard deviation of the bootstrap estimates (AStd), calculated as 
  $T^{-1}\sum_{j=1}^T\sqrt{\widehat{Var}(\sigma^\ast(\tau_j))} $, where 
  %and the averaged mean square error (AMSE) calculated as $ T^{-1}\sum_{j=1}^T\{(\overline{\sigma}^\ast(\tau_j)-\sigma(\tau_j))^2+\widehat{Var}(\sigma^\ast(\tau_j))\}$ of the bootstrap estimator  $ \sigma^\ast(\tau_j)$, where  
  $ \sigma(\tau_j)$  is the  estimated exact standard deviation,  $\widehat{Var}(\sigma^\ast(\tau_j))=(R-1)^{-1}\sum_{r=1}^R(\sigma^\ast_r(\tau_j)-\overline{\sigma}^\ast(\tau_j) )^2 $, with $ \sigma^\ast_r(\tau_j)$ denoting the bootstrap estimate of $ \sigma(\tau_j)$ obtained in the $r$th replication, $r=1,2, \ldots, R$, and $\overline{\sigma}^\ast(\tau_j)=R^{-1}\sum_{r=1}^R  \sigma^\ast_r(\tau_j)$. For the three block bootstrap methods considered we report the results for  two block sizes denoted by $b_1$ and $b_2$, for which the corresponding methods achieve the two lowest ABias  respectively   RBias values. 
  %which deliver the lowest AMSE  denoted by $ b_{M}$ or the lowest averaged squared  bias,  denoted by $ b_{B}$, when  $b_M$ %and $b_B$ are different blocksizes.
Thus  the results presented for the three block bootstrap methods in Table 2 are  those having  the  overall lowest bias. 
%that can be obtained by  these  methods. 
Finally,  for the FSB procedure we report the results for the values $(m,p)=(2,3)$, $(m,p)=(3,3)$ and for  the values of these parameters chosen by the $ \widehat{m}_n
$ and $ AICC$ rule denoted  by $(\widehat{m},\widehat{p})$.

\begin{table*}
\caption{Averaged absolute  bias (ABias), Averaged relative bias (RBias)  and Averaged standard deviation (AStd)  of the moving block bootstrap  (MBB), the tapered block bootstrap (TBB), the stationary bootstrap (SB) and the functional sieve bootstrap (FSB)  estimates of  the standard deviation of the sample mean  $ \overline{X}_n$.   
%We report results for b$b_1$ and $b_2$ denote the block sizes for which the lowest AMSE respectively ABias  is obtained by each %one of the three block bootstrap methods provided $b_1\neq b_2$.
% obtained for  $ R=1000$ replications
} 
\begin{tabular}{llllllllll}
\hline
            &  \multicolumn{2}{l}{MBB}& \multicolumn{2}{l}{TBB}  &  \multicolumn{2}{l}{SB}  &  \multicolumn{3}{l}{FSB}    \\
           &  $b_1=5$ & $b_2$=9 & $b_1=7$&  $b_2=6$ &$b_1=5$ & $b_2=6$ & (2,3) & (3,3) & ($\widehat{m},\widehat{p}$) 
           \\          
\hline
ABias  & 0.206&   0.208 & 0.139 & 0.153 &0.255 & 0.256 & 0.037 & 0.054 & 0.121 \\ 
RBias  & 0.091 & 0.092 & 0.061 &0.068 & 0.112 & 0.113 & 0.016 & 0.024 &0.053\\ 
AStd  &  0.321& 0.406 & 0.350 & 0.312 &0.341 & 0.371 & 0.445 & 0.462 & 0.484 \\ 
\hline
\end{tabular}
\end{table*}

As it is seen from Table 2,   between the three block bootstrap estimators considered, the MBB estimator 
 seems to behave  better that the SB estimator, while  both  estimators are outperformed  by   the  TBB estimator. 
 %The later method performs   best among the three different block bootstrap methods considered. \
 However, 
compared to the FSB estimates,  all block bootstrap estimates  are  quite  biased and they  are  clearly outperformed 
%in terms of bias  
by the FSB estimates. This is true even for the case where the parameters of the FSB procedure are chosen data dependent, where the bias of the FSB  estimates  is smaller that the lowest bias achieved by the block bootstrap methods. 
The FSB estimates have a larger standard deviation which, however,  is not surprising taking into account the fact that 
this bootstrap method  requires the estimation  of $ m^2p$ autoregressive coefficients. It  is worth investigating whether the  standard deviation of the FSB estimates  can be reduced by using  sparse methods 
to fit  the vector autoregression involved in the bootstrap procedure. 

The results of a small simulation study investigating the finite sample size and power  behavior of the bootstrap based, fully functional test for the two-sample mean problem considered in Section~\ref{sec.test},  are presented   in the supplementary material. 
   
\section{Auxiliary Results and Proofs}

\begin{lemma} \label{le.A1} Let Assumption 1, 2 and 3 be satisfied. Denote by $ \Psi_j(m)$, $j=1,2, \ldots $,  the coefficients matrices  of the power series $ A^{-1}_m(z)$,  where $ A_m(z) = I_m-\sum_{j=1}^{\infty} A_{j}(m) z^j$, $ |z| \leq 1$, and  let $ \Sigma_e(m)=E(e_{t}(m)e^{\top}_{t}(m))$.  Then,  
\begin{enumerate}
\item[(i)] \ $ \sum_{j=1}^\infty (1 + j)^{r} \| A_j(m)\|_F  =O(1)$,
\item[(ii)] \ $ \sum_{j=1}^\infty (1 + j)^{r} \| \Psi_j(m)\|_F =O(1)$, \ and
\item[(iii)] \  $ 0 < c_e \leq   \| \Sigma_e(m) \|_F =O(1)$,
\end{enumerate}
where all bounds on the right hand side are valid uniformly in $m$.
\end{lemma}

The following version of Baxter's inequality is very useful in our  setting because it relates the  approximation error  of the coefficient matrices of  the  finite predictor and of  the autoregressive-representation of the $m$-dimensional process of scores to  the lower bound of the spectral density matrix $f_\xi(\cdot)$.  It is an immediate consequence of  Lemma~\ref{le.eigen} and of Theorem 3.2   in  Meyer et al. (2016).

\begin{lemma}\label{le.A2} Let Assumption 1, 2 and 3 be satisfied. Then there exists a constant  $C>0$ which does not depend on $m$,  such that  for all $ 0 \leq s  \leq r-1$,  
\[ \sum_{j=1}^{p}(1 +j)^{s} \| A_{j,p}(m) - A_j(m) \|_F \leq  C\delta_m^{-1} \sum_{j=p+1}^{\infty}(1+j)^{s+1} \| A_j(m)\|_F, \]
where $ \delta_m$ is given in Lemma~\ref{le.eigen}.
\end{lemma}

The following lemma provides a useful bound between the estimated  matrices  of the autoregressive parameters based on the vector of scores $ \xi_t$ and on the vector of their estimates $ \widehat{\xi}_t$, $t=1,2, \ldots, n$. It deals with the case of the  Yule-Walker estimators but similar bounds    can be established along the same lines  for  other estimators, like  for instance for least squares estimators.

\begin{lemma} \label{le.A3} Let  Assumption 1 be satisfied, let  $\widehat{A}_{p,m} =(\widehat{A}_{j,p}(m), j=1,2, \ldots,p)$ and  let  $\widetilde{A}_{p,m} =(\widetilde{A}_{j,p}(m), j=1,2, \ldots,p)$ be  the Yule-Walker
estimators of $ A_{j,p}(m)$, $ j=1,2, \ldots, p$, based on the time series of true scores $ \xi_1, \xi_2, \ldots, \xi_n$.
Then, 
\[ \big\| \widehat{A}_{p,m} - \widetilde{A}_{p,m}\big\|_F =
O_P\Big(\Big(\frac{p\sqrt{m}}{\lambda_m} + p^2\Big)^2\Big\{\frac{1}{n}\sum_{j=1}^{m}\frac{1}{\alpha^{2}_j}\Big\}^{1/2}\Big).\]
\end{lemma}

\begin{lemma} \label{le.A4} Let Assumption 1 and 2 (with $r=0$) be satisfied and $ A_{p,m}(z) = I-\sum_{j=1}^{p}A_{j,p}(m) z^j$, $ z \in \C$.  There exists $ p_m \in \N$ and a positive constant $C$ which does not depend on $m$ such that  for $m\in \N$ and all $ p > p_m$, 
\[ \inf_{|z| \leq 1 + 1/p} \Big| det(A_{p,m}(z))\Big| \geq Cm^{-1/2}.\]
\end{lemma}

%{\bf Proof:} \   We first show that the assertion is true for $ |z| \leq 1$.  Since  $ |det (A_{p,m}(z))| \neq 0$ for $ |z| \leq 1$ it follows by the minimum modulus principle for holomorphic functions that  $|det(A_{p,m}(z)| \geq \inf_{|\overline{z}|=1}|det A_{p,m}(\overline{z})|$.  Now, recall that for $\omega \in [-\pi, \pi]$, 
%$ 2\pi f_\xi (\omega)= A_{m,p}^{-1}(e^{-i\omega}) \Sigma_e(m) \overline{A}_{m,p}^{-1}(e^{-i\omega})$.  Let  $\mu_1(\omega)$ be the largest eigenvalue of $f_\xi(\omega)$,  we then have     
%\begin{align*}
% |det (A_{p,m}(e^{-i\omega}))|^2 &=det (\Sigma_e(m)) / (2\pi |det(f_\xi (\omega))|)\\
% & \geq c_e /(2\pi\, m\, \mu_1(\omega))\\
% & \geq \widetilde C m^{-1},
% \end{align*}
%for some constant $\widetilde C>0$ independent of $m$. Notice that  
% the first  inequality follows by  Lemma~\ref{le.A1}(iii) and  the last by the fact that  $ \mu_1(\omega)$  is  bounded uniformly in $m$; see Lemma~\ref{le.eigen}. Thus $\inf_{\omega\in[-\pi,\pi]} |det (A_{p,m}(e^{-i\omega})|^2 \geq  C m^{-1}$ which implies that  
% $ \inf_{|z|\leq 1}|det A_{p,m}(z)| \geq Cm^{-1/2}$ with some constant $C>0$ independent of $m$. Extension of this lower bound 
%to the slightly larger region  $ |z| \leq 1 + 1/p$ and  for  all $ p > p_m$ for some  $ p_m \in \N$,  follows then exactly along the same lines as the proof of Lemma 3.2 of Meyer and Kreiss (2015); see also Lemma 2.3 of Kreiss et al. (2011). 
%\hfill $\Box$

To state the next lemma we first  fix the following notation.  $\Psi_{j}(m)$, $\Psi_{j,p}(m)$,  $\widetilde{\Psi}_{j,p}(m)$ and $ \widehat{\Psi}_{j,p}(m)$
$j=1,2, \ldots $ denote the coefficient matrices in the power series expansions of $ A^{-1}_m(z)$, $ A_{p,m}^{-1}(z)$, $ \widetilde{A}_{p,m}^{-1}(z)$ and $ \widehat{A}_{p,m}^{-1}(z)$,  respectively, $ |z| \leq 1$. We set $ \Psi_{0}(m)=\Psi_{0,p}(m)=\widetilde{\Psi}_{0,p}(m)=\widehat{\Psi}_{0,p}(m)=I_m$. Furthermore, $e_{t}(m)=\xi_t -\sum_{j=1}^{\infty}A_{j}(m)\xi_{t-j}$, $e_{t,p}(m)=\xi_t -\sum_{j=1}^{p}A_{j,p}(m)\xi_{t-j}$, $\widetilde{e}_{t,p}(m)=\xi_t -\sum_{j=1}^{p}\widetilde{A}_{j,p}(m)\xi_{t-j}$ and $\widehat{e}_{t,p}(m)=\widehat{\xi}_t -\sum_{j=1}^{p}\widehat{A}_{j,p}(m)\widehat{\xi}_{t-j}$, while $\widetilde{\Sigma}_{e,p}(m) =E^+( \widetilde{e}_{t,p}(m)-\overline{\widetilde{e}}_{n,p}(m)) ( \widetilde{e}_{t,p}(m)-\overline{\widetilde{e}}_{n,p}(m))^{\top}$ and $\widehat{\Sigma}_{e,p}(m) =E^\ast( \widehat{e}_{t,p}(m)-\overline{\widehat{e}}_{n,p}(m)) ( \widehat{e}_{t,p}(m)-\overline{\widehat{e}}_{n,p}(m))^{\top}$ with 
$\overline{\widetilde{e}}_{n,p}(m) =(n-p)^{-1}\sum_{t=p+1}^n\widetilde{e}_{t,p}(m)  $ and $\overline{\widehat{e}}_{n,p}(m) =(n-p)^{-1}\sum_{t=p+1}^n\widehat{e}_{t,p}(m)$,  where  $ E^+$ denotes expectation with respect to the measure assigning  probability $ (n-p)^{-1}$ to each $ 
\widetilde{e}_{t,p}(m)$, $ t=p+1, p+2, \ldots, n$.

\begin{lemma} \label{le.A5} Let  Assumptions 1 and 3 and Assumption 2 and 4 (r=2) be satisfied. Then, as $ n \rightarrow \infty$, 
\begin{enumerate}
\item[(i)] \ $ \sum_{j=1}^\infty \| \widetilde{\Psi}_{j,p}(m)- \Psi_{j,p}(m)\|_F \stackrel{P}{\rightarrow}  0$,
\item[(ii)] \ $ \|\widetilde{\Sigma}_{e,p}(m) - \Sigma_{e,p}(m)\|_F \stackrel{P}{\rightarrow} 0$,
\item[(iii)] \ $ \sum_{j=1}^\infty \| \widehat{\Psi}_{j,p}(m)- \Psi_{j,p}(m)\|_F \stackrel{P}{\rightarrow}  0$,
\item[(iv)] \ $ \|\widehat{\Sigma}_{e,p}(m) - \Sigma_{e,p}(m)\|_F\stackrel{P}{\rightarrow} 0$,
\item[(v)] \ $ \sum_{j=1}^\infty \| \Psi_{j,p}(m)- \Psi_{j}(m)\|_F \stackrel{}{\rightarrow}  0$,
\item[(vi)] \ $ \|\Sigma_{e,p}(m) - \Sigma_{e}(m)\|_F \stackrel{}{\rightarrow} 0$.
\end{enumerate}
\end{lemma}

{\bf Proof of Lemma  \ref{le.eigen}:} \  Expression (\ref{eq.gammasum}) imediately leads, for all $\omega \in [0,\pi]$, to an upper bound of  $ f_{{\bf \xi}^{(M)}} (\omega)$. 
To derive a lower bound, recall that $ \Gamma_{\xi^{(M)}}(h)= \big(\langle C_h(v_{j_r}),v_{j_s}\rangle \big)_{r,s=1,2, \ldots, m} $ and observe that 
\[ f_{{\bf \xi}^{(M)}}(\omega)
%=\Big(\langle (2\pi)^{-1}\sum_h C_h(v_{j_r})e^{-ih\omega},v_{j_s}\rangle \Big)_{r,s=1,2, \ldots, m} 
= \Big(\langle {\mathcal F}_\omega(v_{j_r}),v_{j_s}\rangle \Big)_{r,s=1,2, \ldots, m}.  \]
Let $  \mu_j(\omega)$, $j=1,2, \ldots,m$,   be the eigenvalues of $ f_{{\bf \xi}^{(M)}}(\omega)$ (including multiplicity). It suffices to show that $ \min_{1\leq j\leq m}  \mu_j(\omega) \geq \delta_M>0$ for all frequencies $\omega\in [0,\pi]$.
%, i.e., that the eigenvalues of the spectral density matrix $f_{{\bf \xi}^{(M)}}(\omega)$ are uniformly (in $\omega$) 
%bounded away from zero. 
For this let 
$ c_j(\omega)=(c_{j,1}(\omega), c_{j,2}(\omega), \ldots, c_{j,m}(\omega))^\top \in \C^m$, $ j=1,2, \ldots,m$,
 be the corresponding normalized eigenvectors. Then for every $ j \in \{1,2, \ldots, m\}$, we have
\begin{align*}
 \mu_j(\omega)  
 %& = c_j^\top(\omega) f_{{\bf \xi}^{(M)}}(\omega) c_j(\omega)\\
&= c_j^\top(\omega) \Big(\langle {\mathcal F}_\omega(v_{j_r}), v_{j_s}\rangle \Big)_{r,s=1,2, \ldots, m} c_j(\omega)\\
& =\langle {\mathcal F}_\omega(y_j(\omega)),y_j(\omega) \rangle \, >0,
\end{align*}   
by the positivity of $ {\mathcal F}_\omega$, where $ y_j(\omega)=\sum_{r=1}^{m} c_{j,r}(\omega) v_{j_r} \in \overline{V}_M=\overline{sp}\{v_{j_1},v_{j_2},$ $  \ldots, v_{j_m}\}$  and $ \|y_j\|=1$.   Because of the norm summability of the autocovariance matrix function $ \Gamma_{\xi^{M}}(h) $, the spectral density $ f_{\xi^{(M)}} (\omega)$ and consequently the eigenvalues $ \mu_j(\omega)$, $j=1,2, \ldots, m$,  are continuous functions of $ \omega$. Let $ \delta_M(\omega) =\min_{1\leq j\leq m} \mu_j(\omega) $ and notice that  $\delta_M(\omega)$ is continuous in $\omega$ and $\delta_M(\omega)>0$ for all $\omega \in [0,\pi]$. Define $ \delta_M=\min_{\omega \in [0,\pi]}\delta_M(\omega)$ which is positive by the continuity of $ \delta_M(\cdot)$ in the compact interval $[0,\pi]$. Hence   $  \min_{1\leq j\leq m}  \mu_j(\omega) \geq \delta_M>0$ for all  $\omega\in [0,\pi]$. 
 \hfill $\Box$  
 
{\bf Proof of Proposition~\ref{pr.Lp_m}:}
Recall the definition of $ X_t^\ast=\sum_{j=1}^m{\bf 1}_j^\top \xi_t^\ast \widehat{v}_j + U_t^\ast $ and observe  that  $ \xi_t^\ast=\sum_{l=0}^\infty \widehat{\Psi}_{l,p}(m)e^\ast_{t-l}$, where $ \widehat{\Psi}_{0,p}(m)=I_m$ and  the power series $ \widehat{\Psi}_{m,p}(z) = I_m + \sum_{l=1}^{\infty}\widehat{\Psi}_{l,p}(m) z^l= (I_m-\sum_{j=1}^p \widehat{A}_{j,p}(m)z^j)^{-1}$  converges for $ |z| \leq 1$. Write  $ X_t^\ast = \sum_{l=0}^{\infty} \sum_{j=1}^m{\bf 1}_j^\top\widehat{\Psi}_{l,p}(m)e^\ast_{t-l}\widehat{v}_j  + U^\ast_t$ and define
$X_{t,M}^{\ast}=\sum_{l=0}^{M-1}\sum_{j=1}^m {\bf 1}_j^\top \widehat{\Psi}_{l,p}(m) e^\ast_{t-l}\widehat{v}_j+ \sum_{l=M}^\infty\sum_{j=1}^m {\bf 1}_j^\top \widehat{\Psi}_{l,p}(m) e^{\ast}_{t-l,t}\widehat{v}_j +U^\ast_t$, 
where for each $ t \in \Z$, $ \{e_{s,t}^\ast, s \in \Z\}$ is an independent copy of $ \{e^\ast_s, s \in \Z\}$. Notice that 
$ X^\ast_M - X^\ast_{M,M}= \sum_{l=M}^\infty \sum_{j=1}^m {\bf 1}_j^\top \widehat{\Psi}_{l,p}(m) ( e^\ast_{M-l}-e^\ast_{M-l,M} )\widehat{v}_j$.
By Minkowski's inequality we have
\begin{align} \label{eq.l2mapp} 
\sqrt{E\|X_M^\ast - X_{M,M}^\ast\|^2} & \leq \sqrt{E\|\sum_{l=M}^\infty \sum_{j=1}^m {\bf 1}_j^\top \widehat{\Psi}_{l,p}(m) e^\ast_{M-l}\widehat{v}_j  \|^2}  \nonumber \\
& \ \ \ \ + \sqrt{E\|\sum_{l=M}^\infty \sum_{j=1}^m {\bf 1}_j^\top \widehat{\Psi}_{l,p}(m) e^\ast_{M-l,M} \widehat{v}_j \|^2}.
\end{align}
Evaluating the first expectation term we get  using $ \|A\|^2_F = tr(A A^\top)$ and the submultiplicative property of the Frobenius matrix norm,  that 
\begin{align*}
E\|\sum_{l=M}^\infty \sum_{j=1}^m {\bf 1}_j^\top \widehat{\Psi}_{l,p}(m) e^\ast_{M-l}\widehat{v}_j  \|^2 &  = \sum_{l=M}^{\infty} tr\big(\widehat{\Psi}_{l,p}(m) \Sigma^\ast(m)\widehat{\Psi}_{l,p}^\top (m)\big)\\
& \leq  \|\widehat{\Sigma}^{1/2}_{e,p}(m)\|^2_F \sum_{l=M}^\infty \|\widehat{\Psi}_{l,p} (m)\|^2_F,
\end{align*}
where $\widehat{\Sigma}_{e,p}(m) =  \widehat{\Sigma}_{e,p}^{1/2}(m)  \widehat{\Sigma}_{e,p}^{1/2}(m)$.  An identical expression appears for the second expectation term on the right hand side of (\ref{eq.l2mapp}). Applying Minkowski's inequality again we get by the exponential decay of  $ \|\widehat{\Psi}_{l,p}(m)\|_F$, that 
\begin{align*}
\sum_{M=1}^{\infty}\sqrt{E\|X_M^\ast - X_{M,M}^\ast\|^2}  & \leq 2  \|\widehat{\Sigma}^{1/2}_{e,p}(m)\|_F\sum_{M=1}^\infty\sum_{l=M}^\infty \| \widehat{\Psi}_{l,p}(m)\|_F\\
& = 2  \|\widehat{\Sigma}^{1/2}_{e,p}(m)\|_F\sum_{l=1}^\infty  l  \| \widehat{\Psi}_{l,p}(m)\|_F = O_P(1).
\end{align*} 
\hfill $\Box$

{\bf Proof of Theorem~\ref{th.th_main}} Let 
\[ L_{n,m}^+= \frac{1}{\sqrt{n}} \sum_{t=1}^n \sum_{j=1}^m\xi_{j,t}^+ v_j e^{-it\omega},\]
where $ \xi_t^+=(\xi^+_{1,t}, \xi_{2,t}^+, \ldots, \xi_{m,t}^+ )^\top$, $ t=1,2, \ldots, n$ with $ \xi_t^+=\sum_{j=1}^p \widetilde{A}_{j,p}(m) \xi_{t-j}^+  + e^+_t$, where $\widetilde{A}_{j,p}(m)$, $j=1,2, \ldots, p$ are the estimators of the autoregressive parameter matrices based on  the vector time series of true scores  $\xi_t$, $t=1,2, \ldots, n$ and $ e^+_t$ are obtained by i.i.d. resampling from the centered residuals $ \widehat{e}_t=\xi_t - \sum_{j=1}^{p}\widetilde{A}_{j,p}(m) \xi_{t-j}$, $ t=p+1, p+2, \ldots, n$.   That is,  the pseudo-variable  $ L_{n,m}^+$ is obtained  using  the true eingefunctions $v_j$ and  the true 
scores $\xi_{j,t}$ instead of their estimates  $\widehat{v}_j$ and $\widehat{\xi}_{j,t}$ respectively.
Decompose $n^{-1/2}S^\ast_n(\omega)$ as 
\begin{align*}
n^{-1/2}S^\ast_n(\omega) &
%= \frac{1}{\sqrt{n}}\sum_{t=1}^{n}X^\ast_t e^{-i \omega t}\\
%&
= \frac{1}{\sqrt{n}}\sum_{t=1}^{n}\sum_{j=1}^m \xi^+_{j,t} v_j e^{-i t \omega} +  
\frac{1}{\sqrt{n}}\sum_{t=1}^{n}\sum_{j=1}^m \xi^\ast_{j,t} (\widehat{v}_j-v_j) e^{-i t \omega}  \\
& \ \ \  \ +  
\frac{1}{\sqrt{n}}\sum_{t=1}^{n}\sum_{j=1}^m (\xi_{j,t}^\ast -\xi^+_{j,t})v_j e^{-i t \omega} 
+ \frac{1}{\sqrt{n}}\sum_{t=1}^{n}U^\ast_{t,m}e^{-i t \omega} \\
& = L^+_{n,m} +  V_{n,m}^\ast + D^\ast_{n,m} + R^\ast_{n,m}
\end{align*}
with an obvious notation for $L^+_{n,m} $, $ V_{n,m}^\ast$,  $D^\ast_{n,m}$ and $R^\ast_{n,m} $.   Notice that  the terms 
$ V_{n,m}^\ast$ and $ D^\ast_{n,m}$ are due to the fact that, in the bootstrap procedure, the unknown  scores and eigenfunctions are replaced by their sample estimates,  while $R^\ast_{n,m}$ is due to the $m$-dimensional  approximation of the infinite dimensional structure of the underlying process.   Assertion  (i) of the theorem  follows then 
%since   $ R^\ast_{n,m} \stackrel{P}{\rightarrow} 0 $, $ V_{n,m}^\ast \stackrel{P}{\rightarrow} 0 $ and  $D^\ast_{n,m}\stackrel{P}{\rightarrow} 0 $ 
from Lemma~\ref{le.le1Th_main},  ~\ref{le.le2Th_main},  ~\ref{le.le3Th_main}  and  ~\ref{le.le4Th_main} and Slutsky's theorem.

Consider assertion (ii). Since 
\begin{align*}
n^{-1}\|E^\ast S_n^\ast & (\omega) \otimes S^\ast_n(\omega)  -E S_n(\omega) \otimes S(\omega)   \|_{HS} \\
&  \leq \|n^{-1}E^\ast S_n^\ast(\omega) \otimes S^\ast_n(\omega)- 2\pi {\mathcal F}^\ast_{\omega,m} \|_{HS} \\
&\ \ \ \  + 2\pi \|{\mathcal F}^\ast_{\omega,m} -  {\mathcal F}_\omega\|_{HS} + \|n^{-1} ES_n(\omega) \otimes S_n - 2\pi {\mathcal F}_{\omega} \|_{HS},
\end{align*}
it suffices in view of Proposition~\ref{pr.f} and Theorem 2 of  Cerovecki and H\"ormann (2015),  to show that the first term on the right hand side of the above inequality converges to zero in probability.  For this we have
using  $ n^{-1}E^\ast S_n^\ast(\omega) \otimes S_n^\ast(\omega)= n^{-1}\sum_{-n+1}^{n-1}(1-|h|/n)C^\ast_h$, that this term is bounded by 
\[ \sum_{|h| \geq n} \|C^\ast_h\|_{HS} + n^{-1}\sum_{h=-n+1}^{n-1}|h| \|C^\ast_h\|_{HS}.\]
Now,  since $\sum_{h \in \Z}\|C^\ast_{h}\|_{HS} = O_P(1)$ uniformly in $p$ and $m$,  we get that
$  \sum_{|h| \geq n} \|C^\ast_h\|_{HS} =o_P(1)$ and by Kronecker's lemma that $ n^{-1}\sum_{h=-n+1}^{n-1}|h|$ $ \|C^\ast_h\|_{HS}=o_P(1)$.  To verify the uniform boundeness of    $\sum_{h \in \Z}\|C^\ast_{h}\|_{HS} $,  notice first 
that from the expression of  $ C_h^\ast$ given in Section 3.2  we get that
$\sum_{h \in \Z}\|C^\ast_{h}\|_{HS}  \leq   \sum_{h\in \Z}\|\Gamma^\ast_h\|_F +  \|C^\ast_U\|_{HS}$. The square of the second term on the right hand side of the last inequality  equals $ \|E^\ast U^\ast_t\otimes U^\ast_t\|^2_{HS}$ which converges to zero in probability, see the proof of Proposition~\ref{pr.f}.
For the first term we have that $ \sum_{h\in \Z}\|\Gamma^\ast_h\|_F \leq \big(\sum_{j=0}^\infty \|\widehat{\Psi}_{j,p}(m)\|_F\big)^2\|\widehat{\Sigma}_{e,p}(m)\|_F=O_P(1)$ uniformly in $p$ and $m$ by Lemma ~\ref{le.A1} and Lemma~\ref{le.A5}.  \hfill $\Box$

\begin{lemma} \label{le.le1Th_main} Under the assumptions of Theorem~\ref{th.th_main} it holds true that, $  R^\ast_{n,m} \stackrel{P}{\rightarrow} 0$, as $ n \rightarrow \infty$.
\end{lemma}

\begin{lemma} \label{le.le2Th_main}  Under the assumptions of Theorem~\ref{th.th_main} it holds true that,  $ D^\ast_{n,m} \stackrel{P}{\rightarrow} 0$, as $ n \rightarrow \infty$.
\end{lemma}

\begin{lemma} \label{le.le3Th_main}
Under the assumptions of Theorem~\ref{th.th_main} it holds true that, $  V^\ast_{n,m} \stackrel{P}{\rightarrow} 0$, as $ n \rightarrow \infty$.
\end{lemma}

\begin{lemma}   \label{le.le4Th_main} 
Under the assumptions of Theorem~\ref{th.th_main} it holds true that, for all $ \omega \in [-\pi,\pi]$ and as $n \rightarrow \infty$, 
 \[ L^+_{n,m}(\omega) \Rightarrow NC(0,2\pi {\mathcal F}_\omega),\]
 in probability.
\end{lemma}

\section*{Acknowledgements} The author thanks the Editor, the Associate Editor and two referees for their careful reading  and thoughtful comments and questions that led to an improved version of the paper.

\begin{center}
SUPPLEMENTARY MATERIAL
\end{center}

{\bf Online Supplement: ``Sieve Bootstrap for Functional Time Series''}.  The online supplement contains the proofs that were omitted in this paper and additional numerical results.

\newpage

\begin{frontmatter}
\title{Supplement to \\ ``Sieve Bootstrap for Functional Time Series"}
\runtitle{Functional Sieve Bootstrap}
%\thankstext{T1}{Footnote to the title with the ``thankstext'' command.}

%\begin{aug}
%\author{\fnms{Efstathios} \snm{Paparoditis}\thanksref{t1}\ead[label=e1]{stathisp@ucy.ac.cy}}
%%\author{\fnms{Second} \snm{Author}\thanksref{t3,m1,m2}\ead[label=e2]{second@somewhere.com}}
%%\and
%%\author{\fnms{Third} \snm{Author}\thanksref{t1,m2}
%%\ead[label=e3]{third@somewhere.com}
%%\ead[label=u1,url]{http://www.foo.com}}
%
%%\thankstext{t1}{Some comment}
%\thankstext{t1}{Supported in part   by a University of Cyprus Research Grant.}
%%\thankstext{t3}{Second supporter of the project}
%\runauthor{Efstathios Paparoditis}
%
%\affiliation{University of Cyprus
%%\thanksmark{m1} 
%%and Another University\thanksmark{m2}
%}
%
%\address{UNIVERISTY OF CYPRUS\\
%DEPARTMENT OF MATHEMATICS AND STATISTICS\\
%P.O.Box 20537\\
%CY-1678 NICOSIA\\
%CYPRUS\\
%\printead{e1}\\
%\phantom{E-mail:\ }%\printead*{e2}
%}
%%\address{Address of the Third author\\
%%Usually a few lines long\\
%%Usually a few lines long\\
%%\printead{e3}\\
%%\printead{u1}}
%\end{aug}
%
%\begin{keyword}[class=MSC]
%\kwd[Primary ]{62M10, 62M15}
%%\kwd{60K35}
%\kwd[; secondary ]{62G09}
%\end{keyword}
%
%\begin{keyword}
%\kwd{Bootstrap, Fourier transform, Principal components, Karhunen-Lo\`eve expansion, Spectral density operator}
%%\kwd{\LaTeXe}
%\end{keyword}
%
\end{frontmatter}

\setcounter{section}{0}

This supplement contains technical proofs of the results   presented in the main paper  Paparoditis (2016) as well as  some additional numerical results. 
In particular, Section 1 contains the proofs of the  auxiliary lemmas presented in the mentioned paper, Section 2 the proof of Lemma 3.1,  Section 3 the proof of Proposition 3.2, Section 4 the proofs of the lemmas related to Theorem 4.1,  Section 5 the proof of Theorem 4.2 and Section 6 discusses some implementation issues and presents some additional numerical results.  

\section{Proofs of auxiliary lemmas}  {~} \\
%\ In this section we prove some auxiliary lemmas.

{\bf Proof of Lemma 6.1:} \  Consider  (i) and (ii).   Let  ${\mathcal  C}_v$ be   the class of all $m\times m$ matrix-valued functions on $[-\pi,\pi]$ with $ \C^{m\times m} $-valued Fourier coefficient matrices $(F_k, k \in \Z)$ satisfying  the condition  $ \sum_{h\in \Z} (1+|h|)^r \|F_k\|_F < C < \infty$, where 
 $C$ is  independent of $m$.  Then,   $f_{\xi}\in{\mathcal C}_v $  since the autocovariance matrix function of $ f_{\xi}$ satisfies  
  $ \sum_{h\in \Z}(1+|h|)^r \|\Gamma_{\xi}(h)\|_F <  C < \infty$, see (2.2).  Furthermore,  $ f_{\xi}(\omega) = \phi(\omega)\overline{\phi}(\omega)$, with $ \phi$  the optimal factor of $f_{\xi} $; see Cheng and Pourahmadi (1993), p. 116.  From the boundeness conditions it follows  that   $ det(f_{\xi}(\omega)) \geq \delta_m >0$ for all $ m \in \N$, and, therefore,
  $ \phi(\omega)$ is invertible with inverse denoted by $ \phi^{-1}$. Notice that  $ \phi, \phi^{-1} \in {\mathcal C}_v$.   According to Wiener and Masani (1958), Theorem 5.5 and Theorem 5.7, there exist sequences $ \{C_n(m),n\in\N\}$ and $\{D_n(m),n \in \N\}$ which are independent of $t$ such that for all $ t \in \Z$, 
  $ {\mathcal P}_{\overline{sp}\{\xi_{t-j}, j \geq 1\}}(\xi_t)=\sum_{j=1}^{\infty} C_j(m)e_{t-k}(m)$ and $ 
  e_{t}(m)=\sum_{j=0}^{\infty} D_j(m) \xi_{t-j}$, where $ D_0(m)=\Sigma_e^{-1/2}(m)$ and the infinite sums are $L_2$-convergent. The coefficients in the autoregressive and the Wold representation are obtained by setting  $A_0(m)=\Sigma^{1/2}_e(m) D_0(m)=I_m$, 
   $ A_{j}(m)= - \Sigma_e^{1/2}(m) D_j(m)$ and  $ \Psi_j(m) = C_j(m)\Sigma_e^{-1/2}(m)$, where $ C_j(m)$, $j=1,2, \ldots$  and $ D_j(m)$, $j=0,1,2,\ldots$,  are the Fourier coefficients of $ \phi$   and $ \phi^{-1}$, respectively. Since $ \phi, \phi^{-1} \in {\mathcal C}_v$,  we get that  
  $ \sum_{j\in \N} (1+j)^r\|A_j(m)\|_F $ and  $ \sum_{j\in \N} (1+j)^r\|\Psi_j(m)\|_F  $ are bounded uniformly in $m$.  In  (iii) 
 the lower bound   follows from the regularity of  the infinite dimensional process of scores $ \{\xi_t=(\xi_{j,t}, j=1,2, \ldots, )^\top, t\in \Z\}$ which in turn follows   from the regularity of ${\bf X}$. For the upper bound,  let $\sigma_j^{(m)}$, $j=1,2, \ldots, m$, be the (positive) eigenvalues of $ \Sigma_e(m)$. Then, since $ 0 \leq \sum_{j=1}^{m} \sigma_j^{(m)} = \sum_{j=1}^{m} Var(e_{j,t}) \leq \sum_{j=1}^m Var (\xi_{j,t}) = \sum_{j=1}^{m}\lambda_j$ we have  $ \| \Sigma_e(m) \|_F = \sqrt{\sum_{j=1}^m \sigma_j^2} \leq \sqrt{\sum_{j=1}^m \lambda_j^2}
  \leq \| C_0\|_{HS} < \infty$. \hfill $\Box$

{\bf Proof of Lemma 6.3:} \  We first show that 
\begin{equation} \label{eq.difGamma}
 \sup_{-p \leq h \leq p} \| \widehat{\Gamma}_h -\widetilde{\Gamma}_h\|_F  = O_P(\{n^{-1}\sum_{j=1}^{m}\alpha^{-2}_j\}^{1/2}),
 \end{equation}
where  $\widehat{\Gamma}_h =n^{-1}\sum_{t=1}^{n-h}\widehat{\xi}_t \widehat{\xi}_{t+h}^{\top}$,   $\widetilde{\Gamma}_h =n^{-1}\sum_{t=1}^{n-h}\xi_t \xi_{t+h}^{\top}$, $ h=0,1, \ldots, n-1$,   $ \alpha_1 =\lambda_1-\lambda_2$ and $\alpha_j=\min\{\lambda_{j-1}-\lambda_j, \lambda_j-\lambda_{j+1}\}$, $ j=2,3, \ldots, m$. To simplify notation we also write $ \Gamma_h$ for $ \Gamma_\xi (h)$ in what follows. Recall that the covariance matrices introduced refer to the $m$-dimensional vector of scores $ \xi_t=( \xi_{j,t}=\langle X_t,v_j\rangle, j=1,2, \ldots, m)^\top$ or to  its estimator $ \widehat{\xi}_t=(\widehat{\xi}_{j,t}=\langle X_t,\widehat{v}_j\rangle, j=1,2, \ldots, m)^\top$.
Since $ \| \widehat{\Gamma}_h -\widetilde\Gamma_h \|_F \leq \|n^{-1}\sum_{t=1}^{n-h} (\widehat{\xi}_{t+h}-\xi_{t+h})\widehat{\xi}_t^\top \|_F + \| n^{-1}\sum_{t=1}^{n-h} \xi_{t+h}(\widehat{\xi}_{t}-\xi_{t})^\top \|_F$ it  suffices to consider only one of the two terms  on the right hand side of the last bound. By the triangular and the  Cauchy-Schwarz inequality  we have 
\begin{align*}
\|n^{-1}\sum_{t=1}^{n-h}(\widehat{\xi}_{t+h}- & \xi_{t+h})\widehat{\xi}^\top_t\|_F  \leq n^{-1}\sum_{t=1}^{n-h} \|\big(\langle X_{t+h},\widehat{v}_j-v_j\rangle, j=1, \ldots,m\big)^\top\|\\
& \hspace*{2cm} \times \|\big(\langle X_{t+h},\widehat{v}_j\rangle, j=1, \ldots,m\big)^\top \|\\
& \leq \Big(\sum_{j=1}^{m}\|\widehat{v}_j-v_j \|^2 \Big)^{1/2}\frac{1}{n}\sum_{t=1}^{n} \|X_{t}\| \Big(\sum_{j=1}^m\langle X_t,  \widehat{v}_j\rangle^2\big)^{1/2}\\
& =O_P\Big((\sum_{j=1}^{m}\|\widehat{v}_j-v_j \|^2 )^{1/2}\Big),
\end{align*}
with the $O_P$ term  uniformly in $h$.  The assertion follows because by Assumption 1, 
$ \sum_{j=1}^m\|\widehat{v}_j-v_j \|^2=O_P(n^{-1}\sum_{j=1}^m \alpha_{j}^{-2})$;  see H\"ormann and Kokoszka (2010). 

We next proof the assertion of the lemma. First notice that for invertible matrices $ A_n$ and $ B$ such that $ \|A_n -B\|_F \rightarrow 0$ as $ n \rightarrow \infty$, we have  the bound
\begin{align*}
\|A_n^{-1}-B^{-1}\|_F & = \|A_n^{-1}(B-A)B^{-1}\|_F \\
& \leq \|A_n^{-1}-B^{-1}\|_F \|B-A_n\|_F \|B^{-1}\|_F + \|B^{-1}\|^2_F \|B-A_n\|_F,
\end{align*}
from which we get,  for $n$ large enough such that $1 - \|A_n-B\|_F\|B^{-1}\|_F >0 $, the inequality
\begin{equation}\label{eq.ineqAB}
\|A_n^{-1}-B^{-1}\|_F \leq \frac{\|B^{-1}\|^2_F\|A_n-B\|_F}{1 - \|A_n-B\|_F\|B^{-1}\|_F}.
\end{equation}
Then recall the solution of the Yule-Walker equations, 
\[ A_{p,m}=(A_{1,p}(m), A_{2,p}(m), \ldots, A_{p,p}(m)) =  G_1G_{0,p}^{-1},\]
where the $ G_{0,p}\in \R^{mp \times mp}$ matrix is given by 
\[ G_{0,p}=\left( \begin{array}{cccc} \Gamma_0 & \Gamma_1&  \ldots & \Gamma_{p-1} \\
 \Gamma_{-1} & \Gamma_0 & \ldots & \Gamma_{p-2} \\
 \vdots & \vdots &  & \vdots\\
 \Gamma_{-p+1} & \Gamma_{-p+2} & \ldots & \Gamma_0\end{array}\right) \ \ \mbox{and} \ \   G_1=(\Gamma_1, \Gamma_2, \ldots, \Gamma_p).\]
 Let $ \widehat{A}_{p,m}= \widehat{G}_1 \widehat{G}_{0,p}^{-1}$ where $ \widehat{G}_{0,p}$ and $\widehat{G}_1 $ are the same matrices as $ G_{0,p}$ and $ G_1$ with 
  $ \Gamma_h$ replaced by $\widehat{\Gamma}_h$ and let   $ \widetilde{A}_{p,m}= \widetilde{G}_1 \widetilde{G}_{0,p}^{-1}$, where $ \widetilde{G}_{0,p}$ and $\widetilde{G}_1 $ are the same matrices as $ G_{0,p}$ and $ G_1$ with 
  $ \Gamma_h$ replaced by $\widetilde{\Gamma}_h$. We then have
  \begin{align} \label{eq.boundA}
 \|\widehat{A}_{p,m}- \widetilde{A}_{p,m}\|_F & \leq \|\widehat{G}_{0,p}^{-1}-\widetilde{G}_{0,p}^{-1}\|_F\|\widehat{G}_1\|_F + \| \widetilde{G}_{0,p}^{-1}\|_F \|\widehat{G}_1-\widetilde{G}_1 \|_F.
  \end{align}
 We first show that
 \begin{equation} \label{eq.G1}
 \|G_{0,p}^{-1}\|_F = O_P\big(\sqrt{m}\lambda_m^{-1} + p\big).
 \end{equation}
 Toward this notice first  the recursive relation
 \begin{equation}  \label{eq.recG} G_{0,p+1}^{-1} = \left( \begin{array}{cc} G_{0,p}^{-1} & 0 \\ 0 & 0 \end{array}\right) + R_p,
 \end{equation}
 where 
 \[ R_p= \left( \begin{array}{cc} 0 & -J_p\overline{A}_{p,m}\overline{V}_p^{-1/2} \\ 0 &  \overline{V}_p^{-1/2}  \end{array}\right)  \left( \begin{array}{cc} 0 & 0 \\ -\overline{V}_p^{-1/2}\overline{A}_{p,m}^\top J_p & \overline{V}_p^{-1/2}  \end{array}\right), \]
see Brockwell and Davis (1991), Ch. 11.4 and Sowell (1989), where $J_p= \overline{I}_p \otimes I_m$ with $ I_m$ the $m\times m$ unity matrix and $  \overline{I}_p$ the matrix with ones on the diagonal from the bottom left to the top right and zero elsewhere, $ \overline{V}_p =E(\xi_t-\sum_{j=1}^p\overline{A}_{j,p}(m) \xi_{t+j})(\xi_t-\sum_{j=1}^p\overline{A}_{j,p}(m) \xi_{t+j})^\top$ and $\overline{A}_{p,m}=(\overline{A}_{1,p}(m)^\top, \overline{A}_{2,p}(m)^\top, $ $  \ldots, \overline{A}_{p,p}(m)^\top)$ the  coefficient matrices that minimize the ``forward prediction variance'' $ E(\xi_t-\sum_{j=1}^pD_{j,p}(m) \xi_{t+j})(\xi_t-\sum_{j=1}^pD_{j,p}(m) \xi_{t+j})^\top$. We then get from  the recursive relation (\ref{eq.recG}) that 
\[  G_{0,p}^{-1} = \left( \begin{array}{cc} \Gamma_{0}^{-1} & 0 \\ 0 & 0 \end{array}\right) + \sum_{j=1}^p R_j.\]
Using the definition of the matrix $ R_s$  and because $\|\overline{V}_s^{-1/2} \|_F=O(1)$ uniformly in $ s$ and $ m$,  we get 
\begin{align*}
\|R_s\|_F  & \leq  \|\overline{V}_s^{-1/2} \|^2_F\big( 1 + \|J_s\overline{A}_{s,m}\|_F\big)^2 \\
 &  \leq O(1)\big(1 + \sum_{j=1}^s \|\overline{A}_{j,s}(m)\|_F\big)^2 \\
 & = O(1),
 \end{align*}
 since, as in Lemma 6.1,  $  \sum_{j=1}^s \|\overline{A}_{j,s}(m)\|_F =O(1)$ uniformly in $ s$ and $ m$. Thus $ \|\sum_{j=1}^p R_j\|_F \leq \sum_{j=1}^p\|R_j\|_F=O(p)$ and  using the bound 
 $ \|\Gamma_0^{-1}\|_F=\sqrt{\sum_{j=1}^m \lambda_j^{-2}} \leq \sqrt{m}\lambda_m^{-1}$  we conclude that 
 $\|G_{0,p}^{-1}\|_F = O\big(\sqrt{m}\lambda_m^{-1} + p\big)$.
 
We next show that 
\begin{equation} \label{eq.G2}
\|\widehat{G}_{0,p}^{-1}-\widetilde{G}_{0,p}^{-1}\|_F = O_P\Big((p\sqrt{m}\lambda_m^{-1} +p^2)^2 \sqrt{\frac{1}{n}\sum_{j=1}^m \alpha_j^{-2}}\Big).
\end{equation} 
For this notice that using (\ref{eq.ineqAB}) we get 
 \begin{align*}
 \| \widehat{G}_{0,p}^{-1}-\widetilde{G}_{0,p}^{-1}\|_F & \leq \frac{\displaystyle  \|\widetilde{G}_{0,p}^{-1}  \|^2_F\|\widehat{G}_{0,p} - \widetilde{G}_{0,p} \|_F  
 }{\displaystyle  1 - \|\widetilde{G}_{0,p}^{-1}  \|_F\|\widehat{G}_{0,p} - \widetilde{G}_{0,p} \|_F  }\\
& =\frac{\displaystyle  \|G_{0,p}^{-1}  \|^2_F\|\widehat{G}_{0,p} - \widetilde{G}_{0,p} \|_F  
 }{\displaystyle  1 - \|\widetilde{G}_{0,p}^{-1}  \|_F\|\widehat{G}_{0,p} - \widetilde{G}_{0,p} \|_F  }
 +
  \mbox{a lower order term}\\
 & = O_P\big( \| G_{0,p}^{-1}\|^2_F \|\widehat{G}_{0,p} -\widetilde{G}_{0,p} \|_F\big)\\
 %& =O_P\big(\lambda_m^{-1} \big\{m\,n^{-1}\sum_{j=1}^{m}\alpha^{-2}_j\big\}^{1/2}\big),
\end{align*}
and since  by (\ref{eq.difGamma}),
\begin{align*}
 \| \widehat{G}_{0,p} -\widetilde{G}_{0,p} \|_F &  \leq \sum_{i=1}^p\sum_{j=1}^p \|\widehat{\Gamma}_{i-j}-\widetilde{\Gamma}_{i-j}\|_F\\
 & \leq p^2\max_{-p+1 \leq h \leq p-1} \|\widehat{\Gamma}_{h}-\widetilde{\Gamma}_{h}\|_F \\
 & = O_P\Big(p^2\sqrt{n^{-1}\sum_{j=1}^{m}\alpha_j^{-2}}\Big),
 \end{align*} 
 we get using (\ref{eq.G1}) the  assertion (\ref{eq.G2}). 
 
 Furthermore, 
 \begin{align} \label{eq.G3}
 \|\widehat{G}_1\|_F & \leq  \sum_{j=1}^p\|\widehat{\Gamma}_{j}\|_F \nonumber \\
 & \leq \sum_{j=1}^p\|\Gamma_{j}\|_F  +\sum_{ j=1}^p\|\widehat{\Gamma}_{j}-\Gamma_j\|_F \nonumber \\
 & =O(1) +  O_P\Big( p\sqrt{n^{-1}\sum_{j=1}^{m}\alpha_j^{-2}}\Big),
  \end{align}
 where the $O(1)$ term is uniformly in $m$, and, 
 \begin{equation} \label{eq.G4}
 \|\widehat{G}_1-\widetilde{G}_1 \|_F \leq \sum_{j=1}^p\|\widehat{\Gamma}_{j}-\widetilde{\Gamma}_{j}\|_F = O_P\Big(p\sqrt{n^{-1}\sum_{j=1}^{m}\alpha_j^{-2}}\Big).
 \end{equation}
 Thus from (\ref{eq.boundA}) and using the bounds (\ref{eq.G1})-(\ref{eq.G4}) we get that 
 \[  \|\widehat{A}_{p,m}- \widetilde{A}_{p,m}\|_F =O_P\Big((p\sqrt{m}\lambda_m^{-1} +p^2)^2 \sqrt{\frac{1}{n}\sum_{j=1}^m \alpha_j^{-2}} \Big) .\]
 \hfill $\Box$

{\bf Proof of Lemma 6.4:} \   We first show that the assertion is true for $ |z| \leq 1$.  Since  $ |det (A_{p,m}(z))| \neq 0$ for $ |z| \leq 1$ it follows by the minimum modulus principle for holomorphic functions that  $|det(A_{p,m}(z)| \geq \inf_{|\overline{z}|=1}|det A_{p,m}(\overline{z})|$.  Now, recall that for $\omega \in [-\pi, \pi]$, 
$ 2\pi f_\xi (\omega)= A_{m,p}^{-1}(e^{-i\omega}) \Sigma_e(m) \overline{A}_{m,p}^{-1}(e^{-i\omega})$.  Let  $\mu_1(\omega)$ be the largest eigenvalue of $f_\xi(\omega)$. We then have     
\begin{align*}
 |det (A_{p,m}(e^{-i\omega}))|^2 &=det (\Sigma_e(m)) / (2\pi |det(f_\xi (\omega))|)\\
 & \geq c_e /(2\pi\, m\, \mu_1(\omega))\\
 & \geq \widetilde C m^{-1},
 \end{align*}
for some constant $\widetilde C>0$ independent of $m$. Notice that  
 the first  inequality follows by  Lemma 6.1(iii) and  the last by the fact that  $ \mu_1(\omega)$  is  bounded uniformly in $m$; see Lemma 2.1. Thus $\inf_{\omega\in[-\pi,\pi]} |det (A_{p,m}(e^{-i\omega})|^2 \geq  C m^{-1}$ which implies that  
 $ \inf_{|z|\leq 1}|det A_{p,m}(z)| \geq Cm^{-1/2}$ with some constant $C>0$ independent of $m$. Extension of this lower bound 
to the slightly larger region  $ |z| \leq 1 + 1/p$ and  for  all $ p > p_m$ for some  $ p_m \in \N$,  follows exactly along the same lines as the proof of Lemma 3.2 of Meyer and Kreiss (2015); see also Lemma 2.3 of Kreiss et al. (2011). 
\hfill $\Box$

{\bf Proof of Lemma 6.5:} \  To see (i) let $ A^{(r,s)}$  be the $(r,s)$th element of a  matrix $A$ and notice that by  Cauchy's inequality for holomorphic  functions we have 
\begin{equation} \label{eq.Cauchy}
| \widetilde{\Psi}_{j,p}^{(r,s)}(m)- \Psi_{j,p}^{(r,s)}(m)| \leq \big(1+\frac{1}{p}\big)^{-j}\max_{|z|=1+1/p}\|\widetilde{A}^{-1}_{p,m}(z) - A^{-1}_{p,m}(z)\|_F
\end{equation}
and 
\begin{align*}
\max_{|z|=1+1/p} \| \widetilde{A}^{-1}_{p,m}(z)- & A^{-1}_{p,m}(z)\|_F  \leq \max_{|z|=1+1/p}\frac{1}{|det(\widetilde{A}_{p,m}(z)|}\| \widetilde{A}^{Adj}_{p,m}(z) - A^{Adj}_{p,m}(z)\|_F\\
& \ + \max_{|z|=1+1/p}\big|\frac{1}{det(\widetilde{A}_{p,m}(z))} -\frac{1}{det(A_{p,m}(z))} \big|\|A_{p,m}^{Adj}(z)\|_F\\
& = R_{1,n}(z) + R_{2,n}(z),
\end{align*} 
with an obvious notation for $ R_{1,n}(z)$ and  $ R_{2,n}(z)$.   By Theorem 2.12 of Ipsen and Rehman (2008) we get  that 
\[ |det(\widetilde{A}_{p,m}(z) )- det(A_{p,m}(z))| \leq m\|\widetilde{A}_{p,m}(z) - A_{p,m}(z)\|_2 B^{m-1}_{\max}(z),\]
where $ B_{\max}(z)=max\{ \| \widetilde{A}_{p,m}(z)\|_2, \| A_{p,m}(z)\|_2\}$ and $ \|A\|_2$ denotes the spectral norm  (i.e.,  the largest singular value) of the matrix $A$.  Since   $\|A\|_2 \leq \|A\|_F$ and using the bound 
\[ \|A\|_2 \leq \frac{2}{|det(A)|} \Big(\frac{\|A\|_F}{\sqrt{m+1}}\Big)^{m+1},\]
for the largest singular value of  a non-singular matrix $A$, see Merikoski and Kumar (2005), p. 373, we get by straightforward calculations and in view of Lemma 6.4  and the constant $C$ appearing there, that 
\[ \max_{|z|=1+1/p}B_{\max}^{m-1}(z) \leq \frac{2^{m-1}m^{(m-1)/2}\max_{|z|=1+1/p}\|B_{\max}(z)\|_F^{m^2-1} }{C^{m-1}(m+1)^{(m^2-1)/2}}=o_P(1),\]
since $ \max_{|z|=1+1/p} \|A_{p,m}(z)\|_F=O(1) $ and $ \max_{|z|=1+1/p} \|\widetilde{A}_{p,m}(z)\|_F=O_P(1)$ uniformly in $m$. 
Thus 
\begin{equation} \label{eq.det-bound}
|det(\widetilde{A}_{p,m}(z) )- det(A_{p,m}(z))| \leq m\|\widetilde{A}_{p,m}(z) - A_{p,m}(z)\|_F\, o_P(1).
\end{equation}
 Lemma 6.1  and  Assumption 4(iv) lead by  Cauchy-Schwarz's inequality to  the bound 
\begin{align*}
\sup_{|z| \leq  1+1/p}\| \widetilde{A}_{p,m}(z) - A_{p,m}(z) \|  & \leq (1+1/p)^p\sum_{j=1}^p\| \widetilde{A}_{p,m}(z) - A_{p,m}(z) \| \\
& 
\leq O_P( \sqrt{p}\|\widetilde{A}_{p,m} -A_{p,m}\|_F),
\end{align*}
from which we derive, using  (\ref{eq.det-bound}), that 
\begin{align*}
 \sup_{|z| \leq  1+1/p}|det(\widetilde{A}_{p,m}(z) )- det(A_{p,m}(z)) | &  \leq o_P(1) m \sup_{|z| \leq  1+1/p}\| \widetilde{A}_{p,m}(z) - A_{p,m}(z) \|_F\\
% & = o_P(1)O_P(m \sqrt{p}\|\widetilde{A}_{p,m} -A_{p,m}\|_F)\\
 &=o_P(m\sqrt{p}\|\widetilde{A}_{p,m} -A_{p,m}\|_F)
\end{align*}
and by Lemma 6.4  that  
\begin{align*}
R_{1,n}(z) &  \leq \delta_m^{-1}\sum_{r=1}^m\sum_{s=1}^m \sup_{|z|\leq 1+1/p}|det(\widetilde{A}_{p,m}^{(-r,-s)}(z) )- det(A_{p,m}^{(-r,-s)}(z)) | \\
& = O_P(m^{1/2})O_P(m^2)o_P(m\sqrt{p}\|\widetilde{A}_{p,m} -A_{p,m}\|_F) \\
& = o_P(m^{7/2} p^{1/2} \|\widetilde{A}_{p,m} -A_{p,m}\|_F).
\end{align*}
Furthermore,  by Lemma 6.4  and the bound (\ref{eq.det-bound}) we get     
\begin{align*}
R_2(z) & \leq \delta_m^{-2}\max_{|z|=1+1/p}\|A_{p,m}^{Adj}(z)\|_F\max_{|z|=1+1/p} |det(\widetilde{A}_{p,m}(z) )- det(A_{p,m}(z))|\\
 & =o_P(\delta_m^{-2} mp^{1/2}\|\widetilde{A}_{p,m} -A_{p,m}\|_F)\\
 &= o_P(m^{2} p^{1/2}\|\widetilde{A}_{p,m} -A_{p,m}\|_F).
\end{align*}
Thus and using equation (\ref{eq.Cauchy}),  we conclude that
\begin{align*}
\sum_{j=1}^\infty \| \widetilde{\Psi}_{j,p}(m)- \Psi_{j,p}(m)\|_F   & \leq \sum_{j=1}^\infty\sum_{r=1}^m\sum_{s=1}^m | \widetilde{\Psi}_{j,p}^{(r,s)}(m)- \Psi_{j,p}^{(r,s)}(m)| \\
& =  o_P(m^{11/2}p^{3/2} \|\widetilde{A}_{p,m} -A_{p,m}\|_F) \\
& \ \ \ \ \ \ + o_P(m^4p^{3/2}\|\widetilde{A}_{p,m} -A_{p,m}\|_F)\\
&=o_P(m^{3/2}p^{-1/2}) + o_P(p^{-1/2}) \rightarrow 0,
\end{align*}
by Assumption 4.

Consider   (ii) so we have,
\begin{align*}
\|\widetilde{\Sigma}_{e,p}(m) &  - \Sigma_{e,p}(m)\|_F  \leq \|\frac{1}{n-p}\sum_{t=p+1}^n\big(\widetilde{e}_{t,p}(m)\widetilde{e}_{t,p}^{\top} (m) - e_{t,p}(m)e_{t,p}^{\top}(m)\big)  \|_F\\
& \ \ + \|\frac{1}{n-p}\sum_{t=p+1}^ne_{t,p}(m)e_{t,p}^{\top} (m) - Ee_{t,p}(m)e_{t,p}^{\top}(m)  \|_F \\
& \ \ \ \ \ \ + \|\overline{\widetilde{e}}_{n,p}(m)\overline{\widetilde{e}}_{n,p}^{\top}(m) \|_F\\
& = E_{1,n}+E_{2,n} + E_{3,n},
\end{align*} 
with an obvious notation for $ E_{j,n}$, $j=1,2,3$. We show that all three terms converge to zero in probability.   By the triangular inequality and in order  to show  $E_{1,n}\stackrel{P}{\rightarrow} 0$,   it suffices to show that  $E^{(1)}_{1,n} =\| (n-p)^{-1}\sum_{t=p+1}^n( \widetilde{e}_{t,p}(m)- e_{t,p}(m)) \widetilde{e}_{t,p}^{'}(m)\|_F  \stackrel{P}{\rightarrow} 0$.  For this  we use  the bound  
\begin{align} \label{eq.bound1}
E^{(1)}_{1,n} & \leq \|\frac{1}{n-p}\sum_{t=p+1}^n \sum_{j=1}^{p}(\widetilde{A}_{j,p}(m)-A_{j,p}(m)) \xi_{t-j}   \widetilde{e}_{t,p}^{\top}(m) \|_F  \nonumber \\
& + \| \frac{1}{n-p}\sum_{t=p+1}^n \sum_{j=1}^{p}(A_{j,p}(m)-A_{j}(m)) \xi_{t-j}   \widetilde{e}_{t,p}^{\top}(m)   \|_F \nonumber  \\
& + \|\frac{1}{n-p}\sum_{t=p+1}^n \sum_{j=p+1}^{\infty}A_{j}(m) \xi_{t-j}   \widetilde{e}_{t,p}^{\top}(m)  \|_F.
\end{align}
Since by straightforward calculations it yields that $ \sum_{j=1}^p\|\xi_{t-j} \widetilde{e}_{t,p}^{\top}(m) \|^2_F=O_P(m^2p)$,  we get by Assumption 4(iv) and Cauchy-Schwarz's inequality that   
\begin{align*}
\|\frac{1}{n-p}\sum_{t=p+1}^n \sum_{j=1}^{p}(\widetilde{A}_{j,p}(m)  & -A_{j,p}(m)) \xi_{t-j}   \widetilde{e}_{t,p}^{\top}(m) \|_F
 \\
 &=O_P\Big(\|\widetilde{A}_{p,m} -A_{p,m}\|_F\sqrt{ \sum_{j=1}^p\|\xi_{t-j} \widetilde{e}_{t,p}^{\top}(m) \|^2_F}\Big)\\
& =O_P(m^{-3}p^{-3/2}) \rightarrow 0.
 \end{align*}
For the second term on the right hand side of (\ref{eq.bound1}) we get  by replacing $ \widetilde{e}_{t,p}(m) $ by $ e_{t,p}(m)$  and 
using  Lemma 6.2, that 
\begin{align*}
E\| \sum_{j=1}^{p}(A_{j,p}(m)-  A_j(m))\xi_{t-j}e^{\top}_{t,p}(m)& \|_F   \leq \sum_{j=1}^{p}\|A_{j,p}(m)-A_j(m)\|_F O(m)\\
& =O(m \delta^{-1}_m \sum_{j=p+1}^{\infty}j \|A_{j}(m)\|_F)\\
& = O(mp^{-1} \delta^{-1}_m\sum_{j=p+1}^\infty j^2 \|A_j(m)\|_F) \rightarrow 0,
\end{align*}
by Assumption 4. Finally and by the same assumption, we get for the third term of (\ref{eq.bound1})    using 
\begin{align*}
E \sum_{j=p+1}^{\infty} \|A_j(m) \|_F\|\xi_{t-j}e_{t,p}^{\top}(m)\|  = &  O(m)\sum_{j=p+1}^\infty \|A_j(m)\|_F\\
 = & O(m p^{-1} \sum_{j=p+1}^\infty j \|A_j(m)\|_F),
\end{align*}
which  converges to zero in probability. 

Since the  term  $ E_{2,n}$ is easier to deal with using similar arguments as for the term $ E_{1,n}$, we    consider 
 the term $ E_{3,n}$. Using  $ \overline{e}_{n}(m) =(n-p)^{-1}\sum_{t=p+1}^n e_{t}(m)$ we have  that 
\begin{align*}
E_{3,n} & \leq \|(\overline{\widetilde{e}}_{n,p}(m)-\overline{e}_{n}(m) )(\overline{\widetilde{e}}_{n,p}(m) -\overline{e}_{n}(m)  )^{\top} \|_F + \|\overline{e}_{n}(m) \overline{e}_{n}^{\top}(m)   \|_F \\
&\ \ \ \  + 2\|(\overline{\widetilde{e}}_{n,p}(m)-\overline{e}_{n}(m) )\overline{e}_{n}^{\top}(m)  \|_F.
\end{align*}
Since $ \overline{e}_{n}(m) =O_P((n-p)^{-1/2})$ uniformly in $m$ and by similar arguments as above, we have  
\begin{align*}
\|\overline{\widetilde{e}}_{n,p}(m)-\overline{e}_{n}(m) \| & \leq  \frac{1}{n-p}\sum_{t=p+1}^n\sqrt{\sum_{j=1}^p
 \|\widetilde{A}_{j,p}(m)-A_{j,p}(m)\|^2_F }\sqrt{\sum_{j=1}^p\|\xi_{t-j}\|^2} \\
& \ \ \ \  + \frac{1}{n-p}\sum_{t=p+1}^n\sum_{j=1}^p \|A_{j,p}(m)-A_{j}(m)\|_F\|\xi_{t-j}\| \\
& \ \ \ \  + \frac{1}{n-p}\sum_{t=p+1}^n\sum_{j=p+1}^\infty \|A_{j}(m)\|_F\|\xi_{t-j}\| \\
& = O_P(m^{-7/2}p^{-3/2})  + O_P(m^{1/2}\delta_m^{-1} \sum_{j=p+1}^\infty j \|A_j(m)\|_F).
%+ O_P(m^{1/2} \sum_{j=p+1}^\infty\|A_j(m)\|_F).
\end{align*}
Thus we  conclude using  Assumption 4,   that $ E_{3,n} \stackrel{P}{\rightarrow} 0$.

Consider (iii). By (i) it suffices to show that $  \sum_{j=1}^\infty \| \widehat{\Psi}_{j,p}(m)-\widetilde{ \Psi}_{j,p}(m)\|_F \stackrel{P}{\rightarrow}  0$. For this
%, let $ \widehat{A}_{p,m}(z)=I_m -\sum_{j=1}^p \widehat{A}_{j,p}(m)z^j$ and 
% $ \widetilde{A}_{p,m}(z)=I_m -\sum_{j=1}^p \widetilde{A}_{j,p}(m)z^j$ and 
notice  that  by Cauchy-Schwarz's inequality and Lemma 6.3, that
\begin{align*}
\sup_{|z| \leq 1 + 1/p}\|\widehat{A}_{p,m}(z)- \widetilde{A}_{p,m}(z)\|_F  & \leq \big(1 + \frac{1}{p})^p\sum_{j=1}^{p}
\|\widehat{A}_{j,p}(m) -\widetilde{A}_{j,p}(m) \|_F\\
& \leq O(1)O_P(\sqrt{p}\|\widehat{A}_{p,m} - \widetilde{A}_{p,m}\|_F)\\
& = O_P\Big( \Big(\frac{p\sqrt{m}}{\lambda_m} +p^2\Big)^2\sqrt{\frac{p}{n}\sum_{j=1}^{m}\alpha^{-2}_j}\Big). 
\end{align*}
By Cauchy's inequality for holomorphic  functions we get for the $(r,s)$th element of the matrices $ \widehat{\Psi}_{j,p}(m)$ and $ \widetilde{\Psi}_{j,p}(m)$,  that 
\[
\big|\widehat{\Psi}^{(r,s)}_{j,p}(m) - \widetilde{\Psi}^{(r,s)}_{j,p}(m)\big|  \leq \big(1 + \frac{1}{p}\big)^{-j} \max_{|z|=1+1/p}\|\widehat{A}^{-1}_{p,m}(z)- \widetilde{A}^{-1}_{p,m}(z)\|_F \]
and 
\begin{align*}
\max_{|z|=1+1/p} \| \widehat{A}^{-1}_{p,m}(z) & - \widetilde{A}^{-1}_{p,m}(z)\|_F  \leq \max_{|z|=1+1/p}\frac{1}{|det(\widehat{A}_{p,m}(z)|}\| \widehat{A}^{Adj}_{p,m}(z) - \widetilde{A}^{Adj}_{p,m}(z)\|_F\\
& + \max_{|z|=1+1/p}\big|\frac{1}{det(\widehat{A}_{p,m}(z)} -\frac{1}{det(\widetilde{A}_{p,m}(z)} \big|\|\widetilde{A}_{p,m}^{Adj}(z)\|_F.
\end{align*} 
%where $ A^{Adj}$ denotes the adjoint of  the matrix $A$. 
From the above bound and by   Lemma 6.3 and Lemma 6.4, we get by the same  arguments as those leading to the bounds of $R_{1,n}(z) $ and $ R_{2,n}(z)$, that   uniformly in $j$, 
\[ \big|\widehat{\Psi}^{(r,s)}_{j,p}(m) - \widetilde{\Psi}^{(r,s)}_{j,p}(m)\big|  \leq \big(1 + \frac{1}{p}\big)^{-j}O_P\Big( m^{7/2}\Big(\frac{p\sqrt{m}}{\lambda_m} +p^2\Big)^2\sqrt{\frac{p}{n}\sum_{j=1}^{m}\alpha^{-2}_j}\Big),\]
that is 
\begin{align*}
 \| \widehat{\Psi}_{j,p}(m)-\widetilde{ \Psi}_{j,p}(m)\|_F  & \leq \sum_{r,s=1}^{m} \big|\widehat{\Psi}^{(r,s)}_{j,p}(m) - \widetilde{\Psi}^{(r,s)}_{j,p}(m)\big| \\
 & =\big(1 + \frac{1}{p}\big)^{-j}O_P\Big( m^{11/2}\Big(\frac{p\sqrt{m}}{\lambda_m} +p^2\Big)^2\sqrt{\frac{p}{n}\sum_{j=1}^{m}\alpha^{-2}_j}\Big),
 \end{align*}
 from which we get 
 \[ \sum_{j=1}^{\infty}\| \widehat{\Psi}_{j,p}(m)-\widetilde{ \Psi}_{j,p}(m)\|_F = O_P\Big( pm^{11/2}\Big(\frac{p\sqrt{m}}{\lambda_m} +p^2\Big)^2\sqrt{\frac{p}{n}\sum_{j=1}^{m}\alpha^{-2}_j}\Big) \rightarrow 0,\]
by Assumption 4.
 
To establish (iv) notice that using (ii) it suffices to show that  $ \|\widehat{\Sigma}_{e,p}(m) - \widetilde{\Sigma}_{e,p}(m)\| \stackrel{P}{\rightarrow} 0$.  By the triangular inequality it suffices to show that
\[ \|\frac{1}{n-p}\sum_{t=p+1}^n\big[(\widehat{e}_{t,p}(m) - \overline{\widehat{e}}_{n,p}(m)) - (\widetilde{e}_{t,p}(m)) - 
\overline{\widetilde{e}}_{n,p}(m)\big]\big(\widehat{e}_{t,p}(m)-\widetilde{e}_{t,p}(m)\big)\| \stackrel{P}{\rightarrow} 0.\]
Since the above term can be bounded by 
\begin{align*}
\frac{1}{n-p}\sum_{t=p+1}^n\|\widehat{e}_{t,p}(m)-\widetilde{e}_{t,p}(m)\|^2 + \|\overline{\widehat{e}}_{n,p}(m)-\overline{\widetilde{e}}_{n,p}(m)\|\frac{1}{n-p} \| \widehat{e}_{t,p}(m)-\widetilde{e}_{t,p}(m)\|,
\end{align*}
we show that both terms above converge to zero in probability. We use the bound  
\begin{align*}
\frac{1}{n-p}& \sum_{t=p+1}^n\|  \widehat{e}_{t,p}(m)  -\widetilde{e}_{t,p}(m)\|^2  \leq 4\sum_{j=1}^p \|\widehat{A}_{j,p}(m)-\widetilde{A}_{j,p}(m)\|^2_F\frac{1}{n-p}\sum_{t=p+1}^n \|\widehat{\xi}_{t-j}\|^2 \\
& +  \frac{2}{n-p}\sum_{t=p+1}^n \|\widehat{\xi}_t - \xi_t\|^2 + 
 4 \sum_{j=1}^p\|\widetilde{A}_{j,p}(m)\|^2_F \frac{1}{n-p}\sum_{t=p+1}^n\|\widehat{\xi}_{t-j}-\xi_{t-j}\|^2. 
\end{align*}
From  Lemma 6.3 we get by  straightforward calculations  that,  $ (n-p)^{-1}\sum_{t=p+1}^n \|\widehat{\xi}_t - \xi_t\|^2 =O_P(n^{-1}\sum_{j=1}^m\alpha_j^{-2})$ and because $ (n-p)^{-1}\sum_{t=p+1}^n\|\widehat{\xi}_{t-j}\|^2=O_P(1) $, we get 
\begin{align*}\sum_{j=1}^p \|\widehat{A}_{j,p}(m)-& \widetilde{A}_{j,p}(m)\|^2_F\frac{1}{n-p}\sum_{t=p+1}^n \|\widehat{\xi}_{t-j}\|^2 \\
& = \|\widehat{A}_{p,m}-\widetilde{A}_{p,m}\|^2_F\frac{1}{n-p}\sum_{t=p+1}^{n}\|\widehat{\xi}_{t-j}\|^2 \\
& =o_P(1),
 \end{align*}
 by Assumption 4. Furthermore, since $ \sum_{j=1}^p\|\widetilde{A}_{j,p}(m)\|^2_F = O_P(1)$,  we get 
\[ \sum_{j=1}^p\|\widetilde{A}_{j,p}(m)\|^2_F \frac{1}{n-p}\sum_{t=p+1}^n\|\widehat{\xi}_{t-j}-\xi_{t-j}\|^2= O_P(1)O_P\Big(\frac{1}{n} \sum_{j=1}^m \alpha_j^{-2} \Big)=o_P(1).\]
Similar arguments yield
\begin{align*}
\|\overline{\widehat{e}}_{n,p}(m) \|^2  & \leq  2\|\frac{1}{n-p}\sum_{t=p+1}^n \widehat{\xi}_t\|^2 + 2 \|\sum_{j=1}^p\widehat{A}_{j,p}(m)\frac{1}{n-p}\sum_{t=p+1}^n \widehat{\xi}_{t-j}\|^2\\
& = O_P\Big( m (n-p)^{-1} + n^{-1}\sum_{j=1}^m\alpha_j^{-2}\Big),
\end{align*}
and $  \| \overline{\widetilde{e}}_{n,p}(m)\|^2 = O_P(m (n-p)^{-1}) \rightarrow 0$.  

The proof of (v) and (vi) is straightforward and uses  Lemma 6.1 and   Lemma 6.2.

\section{Proof of Lemma 3.1}  {~}\\
(i) Notice that $ \lambda_j \geq C_\lambda \rho^{j}$. Since $ p = O(n^a)$ with $ a \in (0,1/14)$, Assumption 4(i) is satisfied 
because  $ m^{3/2}/p^{1/2} = O(log^{3/2}(n)/n^{a/2})$. For Assumption 4(ii) we have 
\begin{align*}
\frac{p^7}{\sqrt{n}\lambda^2_m}\sqrt{\sum_{j=1}^m\alpha_j^{-2}} & \leq  \frac{p^7}{\sqrt{n}C_\lambda^3\rho^{3m}}\sqrt{\sum_{j=1}^m\rho^{2(m-j)}}\\
& \leq  \frac{p^7}{\sqrt{n}C_\lambda^3\rho^{3m}}\times\frac{1}{\sqrt{1-\rho^2}} \rightarrow 0,
\end{align*}
if  $ n^{1/2 - 7a} \rho^{3m}  \rightarrow \infty$, which is satisfied for $$ m \leq \Big( \frac{\displaystyle 1}{\displaystyle 6 log(\rho^{-1})}\big(1 - 14 a\big) - \delta\Big)log(n)$$ for some $ \delta >0$. 
Finally, straightforward calculations as for Assumption 4(ii) show that  
$ mp^6=O({\sqrt{n}\lambda_m^2}) $  and $ p\lambda_m^2 = O(m^2) $  which imply that Assumption 4(iv) is satisfied.

(ii) Notice that $\lambda_j^{-1} \leq C_\lambda^{-1} j^{\theta}$ and recall that $ p = O(n^a)$ with $ a \in (0,1/14)$.  Then $ m^{3/2}/p^{1/2} = O(n^{3\zeta/2 - a/2}) =O(1)$ for $ 0< \zeta \leq \alpha/3$. Consider Assumption 4(ii) and observe that 
\[ \sum_{j=1}^{m} \alpha_{j}^{-2} \leq  C_\lambda^{-2} \sum_{j=1}^{m} j^{2\theta} \leq C_\lambda^{-2} \frac{1}{2\theta +1}(m+1)^{2\theta+1}.\] 
Thus  
\begin{align*}
\frac{p^7}{\sqrt{n}\lambda^2_m}\sqrt{\sum_{j=1}^m\alpha_j^{-2}}  & = O\big(n^{-(1/2 -7a-3\theta\zeta -\zeta/2)}\big)  \rightarrow 0,
\end{align*}  
 for $\zeta \in (0, \zeta_{\max}]$ and $ \zeta_{\max} =\min\Big\{\frac{\displaystyle 1-14a}{\displaystyle 1+6\theta} - \delta, a/3\Big\}$ and for some $\delta >0$. Finally, verify that for $ \zeta < (1-14a)(1+6\theta)^{-1} $ we have that $ mp^6=O(\sqrt{n}\lambda^2_m)$ and that if  $ \zeta \geq a/(2+2\theta)$,  then  $ p\lambda_m^2 =O(m^2)$.

\section{Proof of Proposition 3.2} {~}\\
   Recall that the spectral density operator ${\mathcal F}_\omega$ can be expressed as 
$2\pi {\mathcal F}_\omega = \sum_{j=1}^\infty \sum_{l=1}^\infty \sum_{h=-\infty}^\infty \langle C_h(v_j),v_l\rangle e^{-ih \omega}(v_j 
\otimes v_l)$.           
Define for $m\in \N$, 
$ 2\pi  {\mathcal F}_{\omega,m} =  \sum_{j=1}^m \sum_{l=1}^m \sum_{h=-\infty}^\infty \langle C_h(v_j),v_l\rangle e^{-ih \omega}(v_j 
\otimes v_l)$ and verify  that since $ \langle C_h(v_j),v_l\rangle =E(\xi_{j,t} \xi_{l,t+h})$,   
%$2\pi  {\mathcal F}_{\omega,m}$ can be also expressed as 
the following expression is also valid, 
$$2\pi {\mathcal F}_{\omega,m}(\cdot)=\sum_{j=1}^m\sum_{l=1}^m {\bf 1}_j^{\top} \Psi_{m}(e^{-i\omega}) \Sigma_{e}(m)\overline{\Psi}_{m}(e^{-i\omega}){\bf 1}_l \langle v_j, \cdot \rangle v_l, $$
where $ \Psi_m(z) = I_m + \sum_{j=1}^\infty \Psi_{j}(m)z^j = (I_m-\sum_{j=1}^\infty A_{j}(m) z^j)^{-1}$,  $ |z| \leq 1$.
% and 
%notice that $ 2\pi F_{\omega,m}=\Psi_{m}(e^{-i\omega})\Sigma_{e}(m)\overline{\Psi}_{m}(e^{-i\omega}) $, where
%$ \Psi_{m}(e^{-i\omega})=\sum_{j=1}^{\infty}\Psi_{j}(m)e^{-i\omega j}$. 
Let $ 2\pi  \widetilde F_{\omega,m}= \sum_{j=1}^m\sum_{l=1}^m {\bf 1}_j^{\top}\Psi_{p,m}(e^{-i\omega}) \Sigma_{e,p}(m)\overline{\Psi}_{p,m}(e^{-i\omega}) {\bf 1}_l \langle v_j, \cdot \rangle v_l$ where
$ \Psi_{p,m}(z)=\sum_{j=1}^{\infty}\Psi_{j,p}(m)z^{j} $, $ |z| \leq 1$.  Finally recall that 
$$ 2\pi  F^\ast_{\omega,m} =  \sum_{j=1}^m\sum_{l=1}^m  {\bf 1}_j^{\top}\widehat{\Psi}_{p,m}(e^{-i\omega}) \widehat{\Sigma}_{e,p}(m)\overline{\widehat{\Psi}}_{p,m}(e^{-i\omega}) {\bf 1}_l \langle \widehat{v}_j, \cdot \rangle \widehat{v}_l + E^\ast U^\ast_t \otimes U^\ast_t.$$ 
Then, 
\begin{equation} \label{eq.speczer}
 \| {\mathcal F}_{\omega,m}^\ast-{\mathcal F}_\omega \|_{HS}\leq \|{\mathcal F}^\ast_{\omega,m}- \widetilde{\mathcal F}_{\omega,m}  \|_{HS} + \|\widetilde{\mathcal F}_{\omega,m}- {\mathcal F}_{\omega,m}  \|_{HS} +\|{\mathcal F}_{\omega,m}- {\mathcal F}_{\omega}\|_{HS}.
 \end{equation}
The first term on the right hand side above is bounded by 
\begin{align*}
\|&{\mathcal F}^\ast_{\omega,m}  - \widetilde{\mathcal F}_{\omega,m}  \|_{HS}  \leq  \|E^\ast U^\ast_t \otimes U^\ast_t\|_{HS}\\
 + & \|\sum_{j,l} {\bf 1}^{'}_j \widehat{\Psi}_{p,m}(e^{-i\omega}) \widehat{\Sigma}_{e,p}(m)\overline{\widehat{\Psi}}_{p,m}(e^{-i\omega}) {\bf 1}_l \big(\langle \widehat{v}_j, \cdot \rangle \widehat{v}_l- \langle v_j, \cdot \rangle v_l\big)\|_{HS}\\
 +&  \|\sum_{j,l} {\bf 1}^{'}_j \big(\widehat{\Psi}_{p,m}(e^{-i\omega}) \widehat{\Sigma}_{e,p}(m)\overline{\widehat{\Psi}}_{p,m}(e^{-i\omega})- \Psi_{p,m}(e^{-i\omega}) \Sigma_{e,p}(m)\overline{\Psi}_{p,m}(e^{-i\omega})\big)\\ & \ \ \ \ \  \times  {\bf 1}_l \langle v_j, \cdot \rangle v_l\|_{HS}.
\end{align*}
Furthermore, 
\[  \|E^\ast U^\ast_t \otimes U^\ast_t\|_{HS} \leq
 \|E^+ U^+_t \otimes U^+_t\|_{HS}+ \|E^\ast U^\ast_t \otimes U^\ast_t -E^+ U^+_t \otimes U^+_t \|_{HS} ,\]
 where $ U^+_t$ are i.i.d. random variables taking values with probability $n^{-1}$ in the set $ \{U_t^c=U_t-\overline{U}_n, t=1,2, \ldots, n\}$ and  $ \overline{U}_n=n^{-1}\sum_{t=1}^{n}U_t$.  Then
 \begin{align*}
  \|E^+ U^+_t \otimes U^+_t\|_{HS} & 
  %\leq \|\frac{1}{n}\sum_{t=1}^n\langle U_t,\cdot\rangle U_t\|_{HS} + \|\langle \overline{U}_n,\cdot
 % \rangle \overline{U}_n\|_{HS}\\
   \leq  \|\sum_{j,l=m+1}^\infty \langle \widehat{C}_0(v_j)v_l\rangle \langle v_j,\cdot  \rangle v_l\|_{HS} + \|\overline{U}_n\| \rightarrow 0,
  \end{align*}
  in probability, since $\| \widehat{C}_0-C_0\|_{HS} \rightarrow 0$ and the operator $ C_0$ is  Hilbert-Schmidt.  
  Furthermore, 
  \begin{align*}
  \|E^\ast U^\ast_t \otimes U^\ast_t -E^+ U^+_t \otimes U^+_t \|_{HS}  & \leq  \|\frac{1}{n}\sum_{t=1}^n \big(\langle\widehat{U}_t ,\cdot\rangle \widehat{U}_t  - \langle U_t,\cdot \rangle U_t\big)  \|_{HS}\\
  & \ \ \ + \| \langle\overline{\widehat{U}}_n,\cdot \rangle\overline{\widehat{U}}_n - \langle \overline{U}_n,\cdot\rangle\overline{U}_n\|_{HS}\\
  & = O_P\Big(\sqrt{\sum_{j=1}^m\|\widehat{v}_j-v_j\|^2}\Big) \rightarrow 0,
  \end{align*}
 in probability, where the last equality follows  by straightforward calculations  and using  $ \widehat{U}_t-U_t=\sum_{j=1}^m(\langle  X_t, \widehat{v}_j \rangle \widehat{v}_j - \langle  X_t, v_j \rangle v_j) $.   
 Similarly, and by the same arguments as above and using Lemma 6.5, we get  
 \begin{align*}
 \|\sum_{j,l} {\bf 1}^{'}_j \widehat{\Psi}_{p,m}(e^{-i\omega}) &  \widehat{\Sigma}_{e,p}(m)\overline{\widehat{\Psi}}_{p,m}(e^{-i\omega}) {\bf 1}_l \big(\langle \widehat{v}_j, \cdot \rangle \widehat{v}_l- \langle v_j, \cdot \rangle v_l\big)\|_{HS}\\
 & = O_P\Big(\sqrt{\sum_{j=1}^m\|\widehat{v}_j-v_j\|^2}\Big) \rightarrow 0,
 \end{align*}
   
  Finally, straightforward calculations yield 
  \begin{align*}
 \| &\sum_{j,l} {\bf 1}^{'}_j \big(\widehat{\Psi}_{p,m}(e^{-i\omega}) \widehat{\Sigma}_{e,p}(m)\overline{\widehat{\Psi}}_{p,m}(e^{-i\omega})- \Psi_{p,m}(e^{-i\omega}) \Sigma_{e,p}(m)\overline{\Psi}_{p,m}(e^{-i\omega})\big) \\ & \ \ \ \ \ \times{\bf 1}_l \langle v_j, \cdot \rangle v_l\|_{HS}\\
 & =O_P(\sum_{j=1}^\infty \|\widehat{\Psi}_{j,p}(m) -\Psi_{j,p}(m) \|_F + \|\Sigma_{e,p}(m)-\widehat{\Sigma}_{e,p}(m) \|_F ) \rightarrow 0,
 \end{align*}
 by Lemma 6.5(iii) and (iv). This concludes the proof that $\|{\mathcal F}^\ast_{\omega,m}  - \widetilde{\mathcal F}_{\omega,m}  \|_{HS}\stackrel{P}{\rightarrow} 0 $.
 
%Now, $\|{\mathcal F}^\ast_{\omega,m}- \widetilde{\mathcal F}_{\omega,m}  \|_{HS}=O_P(\sum_{j=1}^{\infty}\|\widehat{\Psi}_{j,p}(m)-\Psi_{j,p}(m)\|_F ) + O_P(\|\widehat{\Sigma}_{e,p}^\ast(m)-\Sigma_{e,p}(m)\|_F )\rightarrow 0$, in probablity,  by Lemma~\ref{le.A5}(iii) and (iv).
%Similarly,  
% $ 2\pi {\mathcal F}_{\omega,m}=\Psi_{m}(e^{-i\omega})\Sigma_{e}(m)\overline{\Psi}_{m}(e^{-i\omega}) $, where
%$ \Psi_{m}(e^{-i\omega})=\sum_{j=1}^{\infty}\Psi_{j}(m)e^{-i\omega j}$,
Consider next  the second term on the right hand side of (\ref{eq.speczer}). For this term we get   $\|\widetilde{\mathcal F}_{\omega,m}-  {\mathcal F}_{\omega,m}  \|_{HS}=
O_P(\sum_{j=1}^{\infty}\|\Psi_{j,p}(m)-\Psi_{j}(m)\|_F ) + O_P(\|\Sigma_{e,p}(m)-\Sigma_{e}(m)\|_F )$,  i.e., $\|\widetilde{\mathcal F}_{\omega,m}-  {\mathcal F}_{\omega,m}  \|_{HS}$ converges to zero in probability  by Lemma 6.5(v) and (vi). For the third and last term  on the right hand side of  (\ref{eq.speczer}) we obtain
\begin{align*}
\|{\mathcal F}_{\omega,m}- {\mathcal F}_{\omega}\|_{HS} &\leq \|\sum_{j=m+1}^{\infty}\sum_{l=1}^{m}\langle {\mathcal F}_\omega(v_j),v_l \rangle (v_j\otimes v_l)\|_{HS}\\
 & \ \ +  \|\sum_{j=1}^{m}\sum_{l=m+1}^{\infty}\langle{\mathcal F}_\omega(v_j),v_l \rangle (v_j\otimes v_l)\|_{HS}\\
 & \ \ +  \|\sum_{j=m+1}^{\infty}\sum_{l=m+1}^{\infty}\langle {\mathcal F}_\omega(v_j),v_l \rangle (v_j\otimes v_l)\|_{HS} \rightarrow 0,
\end{align*}
as $m \rightarrow \infty$, since 
 $ \{(v_j\otimes v_l), j=1,2, \ldots, l=1,2, \ldots  \}$ is a complete  orthonormal basis of $ {\mathcal H} \otimes  {\mathcal H}$.
\hfill $\Box$

\section{Proofs of the lemmas used for the proof of   Theorem 4.1}  {~} \\

{\bf Proof of Lemma 6.6:} \  Note that 
\begin{align*}
E^\ast\|R^\ast_{n,m} \|^2& 
%= E^\ast \| U^\ast_{1,m}, U^\ast_{1,m}\|^2 \\
 %&
  = \frac{1}{n}\sum_{t=1}^{n}\|\widehat{U}_{t,m}-\overline{\widehat{U}}_{n}\|^2\\
 & \leq   \frac{2}{n}\sum_{t=1}^{n}\|\widehat{U}_{t,m}\|^2 + 2\|\overline{\widehat{U}}_{n}\|^2.
\end{align*}
Using $\|\widehat{v}_j-v_j\|\leq 2\sqrt{2}\alpha_j^{-1}\|\widehat{C}_0-C_0\|_{HS}$, see  H\"ormann and Kokoszka (2010),  we get  
\begin{align*}
\frac{1}{n}\sum_{t=1}^{n}\|\widehat{U}_{t,m}\|^2 & \leq \frac{4}{n}\sum_{t=1}^{n}\|\sum_{j=1}^{m}\langle X_t,v_j\rangle(v_j-\widehat{v}_j)\|^2 + \frac{4}{n}\sum_{t=1}^{n}\|\sum_{j=1}^{m}\langle X_t,v_j-\widehat{v}_j\rangle\widehat{v}_j\|^2  \\
& =  4\sum_{j=1}^m\sum_{l=1}^m \langle \widehat{C}_0(v_j),v_l\rangle \langle v_j-\widehat{v}_j,   v_l-\widehat{v}_l\rangle + 4\sum_{j=1}^m\frac{1}{n}\sum_{t=1}^n\langle X_t, v_j-\widehat{v}_j\rangle^2\\
& \leq  4 \|\widehat{C}_0\|_{HS}\big(\sum_{j=1}^m \|\widehat{v}_j-v_j\|\big)^2 +4 \|\widehat{C}_0\|_{HS}\sum_{j=1}^m \|\widehat{v}_j-v_j\|^2 \\
& \leq 32\|\widehat{C}_0\|_{HS}  \|\widehat{C}_0-C_0\|^2_{HS}  \Big( \big(\sum_{j=1}^m\alpha_j^{-1} \big)^2  +\sum_{j=1}^m \alpha_j^{-2}\Big)\\
& = O_P\big(n^{-1/2}\sum_{j=1}^{m}\alpha_j^{-1}\big),
\end{align*} 
where the last equality follows because  $   \|\widehat{C}_0-C_0\|_{HS}=O_P(n^{-1/2})$. Furthermore,  $\|\overline{\widehat{U}}_{n}\|^2\stackrel{P}{\rightarrow} 0$ follows using similar arguments and since 
$ \overline{\widehat{U}}_{n}=\overline{U}_n + n^{-1}\sum_{t=1}^n \sum_{j=1}^{m}
\big(\langle X_t,v_j \rangle v_j - \langle X_t,\widehat{v}_j \rangle \widehat{v}_j\big)$, where  $ \overline{U}_n=n^{-1}\sum_{t=1}^n U_{t,m}$.
 \hfill $\Box$

{\bf Proof of Lemma 6.7:} \  We have 
\begin{align*}
E\|\frac{1}{\sqrt{n}}\sum_{t=1}^n  \sum_{j=1}^m (\xi_{j,t}^\ast- \xi_{j,t}^+)v_j\|^2= & \frac{1}{n}\sum_{t,s=1}^n\sum_{j=1}^m{\bf 1}_j^{\top} E\xi^\ast_t(\xi_s^\ast-\xi^+_s)^{\top}{\bf 1}_j  \\
& \ \ + \frac{1}{n}\sum_{t,s=1}^n\sum_{j=1}^m{\bf 1}_j^{\top} E\xi^+_t(\xi_s^+ -\xi^\ast_s)^\top{\bf 1}_j \\ &   =  
D^{(1)}_{n,m} +D^{(2)}_{n,m},
\end{align*}
with an obvious notation for $D^{(1)}_{n,m}  $ and $D^{(2)}_{n,m}  $. We   consider  $D^{(1)}_{n,m} $ only since $D^{(2)}_{n,m} $ can be  handled similarly. For this term we have
\begin{align} \label{eq.Dnm1}
D^{(1)}_{n,m}& = \frac{1}{n}\sum_{t,s=1}^n\sum_{j=1}^{m}\sum_{l=0}^{\infty}{\bf 1}_j^\top \widetilde{\Psi}_{l,p}(m)\Sigma^\ast_{e,p}(m)(\widehat{\Psi}_{l+s-t,p}(m)-\widetilde{\Psi}_{l+s-t,p}(m) )^\top {\bf 1}_j  \nonumber \\
& + \frac{1}{n}\sum_{t,s=1}^n\sum_{j=1}^{m}\sum_{l=0}^{\infty}{\bf 1}_j^\top \widetilde{\Psi}_{l,p}(m)E\big[ e^\ast_{t,p}(m)(e_{t,p}^\ast(m)- e^+_{t,p}(m))\big]\widetilde{\Psi}_{l+s-t,p}(m)^\top {\bf 1}_j 
\end{align}
and, using Lemma 6.1 and 6.5 we get for the first term on the right hand side of (\ref{eq.Dnm1}), that, this term 
   is bounded by
\[ \|\Sigma^\ast_{e,p}(m)\|_F\sum_{l=0}^\infty\|\sum_{j=1}^{m}{\bf 1}_j^\top\widetilde{\Psi}_{l,p}(m)\|_F \sum_{l=0}^\infty\|\sum_{j=1}^{m}{\bf 1}_j^\top(\widehat{\Psi}_{l,p}(m)-\widetilde{\Psi}_{l,p}(m))\|_F \rightarrow 0, \]
in probability. The second term of (\ref{eq.Dnm1})   is bounded by
\begin{align*}
 \sqrt{E\|e^\ast_{t,p}(m)\|^2}\sqrt{E\| e^\ast_{t,p}(m) -  e^\ast_{t,p}(m)\|^2} &  \sum_{l=0}^\infty \|\sum_{j=1}^{m}{\bf 1}_j^\top\widetilde{\Psi}_{l,p}(m)\|_F  \\ & \times \sum_{l=0}^\infty\|\sum_{j=1}^{m}{\bf 1}_j^\top\widehat{\Psi}_{l,p}(m)\|_F ,
 \end{align*}
which converges to zero in probability,  because $E\| e^\ast_{t,p}(m) -  e^+_{t,p}(m)\|^2 \rightarrow 0$ in probability. This follows since  
\begin{align*}
E \| e^\ast_{t,p}(m) -  e^+_{t,p}(m)\|^2  \leq &  \frac{2}{n-p}\sum_{t=p+1}^n\|\widehat{e}_{t,p}(m)-\widetilde{e}_{t,p}(m)\|^2 + 4 \big(\|\overline{\widehat{e}}_{n}\|^2 + \|\overline{\widetilde{e}}_{n}\|^2\big)\\
\leq & \frac{4}{n-p}\sum_{t=p+1}^n\|\widehat{\xi}_t -\xi_t\|^2  \\
& +  \frac{4}{n-p}\sum_{t=p+1}^n\|\sum_{j=1}^p(
\widehat{A}_{j,p}(m)\widehat{\xi}_{t-j} -\widetilde{A}_{j,p}(m)\xi_{t-j}) \|^2 \\
& + 4 \big(\|\overline{\widehat{e}}_{n}\|^2 + \|\overline{\widetilde{e}}_{n}\|^2\big),
\end{align*}
and 
\[ \frac{1}{n-p}\sum_{t=p+1}^n\|\widehat{\xi}_t -\xi_t\|^2  \leq \frac{1}{n-p}\sum_{t=p+1}^n\|X_t\|^2 \sum_{j=1}^m \|\widehat{v}_j - v_j\|^2 = O_P\big(n^{-1}\sum_{j=1}^m \alpha_j^{-2}\big)\rightarrow 0.\]
Furthermore,
\begin{align*}
 \frac{1}{n-p}\sum_{t=p+1}^n\|\sum_{j=1}^p\big(
\widehat{A}_{j,p}(m)& \widehat{\xi}_{t-j}  - \widetilde{A}_{j,p}(m)\xi_{t-j}\big) \|^2 \\
 & \leq  2\sum_{j=1}^p\|\widehat{A}_{j,p}(m)\|^2_F\frac{1}{n-p}\sum_{t=p+1}^n\|\widehat{\xi}_{t-j}-\xi_{t-j}\|_F\\
& +  2\sum_{j=1}^p\|\widehat{A}_{j,p}(m)-\widetilde{A}_{j,p}(m)\|_F^2\frac{1}{n-p}\sum_{t=p+1}^n\|\xi_{t-j}\|^2\\
&=O_P\big(n^{-1}\sum_{j=1}^m\alpha_j^{-2}\big) + O_P\big(\lambda_m^{-2}n^{-1}mp\sum_{j=1}^m \alpha_j^2\big) \rightarrow 0,
\end{align*}
where the last equality follows using Lemma 6.1 and Lemma 6.3.
Finally,  
\begin{align*}
\|\overline{\widehat{e}}_n\|^2 & \leq 2\|\frac{1}{n-p}\sum_{t=p+1}^n \widehat{\xi}_{t}\|^2 + 2 \big(\sum_{j=1}^p\|\widehat{A}_{j,p}(m)\|_F\big)^2\big(\|\frac{1}{n-p}\sum_{t=p+1}^{n}\widehat{\xi}_{t-j}\| \big)^2 \rightarrow 0\,
\end{align*}
in probability, since 
\begin{align*}
 \|\frac{1}{n-p}\sum_{t=p+1}^n \widehat{\xi}_{t}\|^2  & \leq 2 \|\frac{1}{n-p}\sum_{t=p+1}^n \xi_{t}\|^2  + O_P(n^{-1}\sum_{j=1}^m\alpha_j^{-2})\\
& = O_P(m/(n-p)) +  O_P(n^{-1}\sum_{j=1}^m\alpha_j^{-2}) \rightarrow 0,
\end{align*}
and
\begin{align*}
\sum_{j=1}^p\|\widehat{A}_{j,p}(m)\|_F \|\frac{1}{n-p}\sum_{t=p+1}^{n}& \widehat{\xi}_{t-j}\|  =  O_P(1)\\
& \ \ \times O_P\Big(\sqrt{m/(n-p)} + \sqrt{n^{-1}\sum_{j=1}^m\alpha_j^{-2}}\Big) \rightarrow 0.
\end{align*}
By similar arguments we get $\|\overline{\widetilde{e}}_n\|^2 \rightarrow 0$, in probability.
\hfill $ \Box$

{\bf Proof of Lemma 6.8:} \ 
\begin{align*}
E\big\|\frac{1}{\sqrt{n}}\sum_{t=1}^n \sum_{j=1}^m \xi_{j,t}^\ast(\widehat{v}_j - v_j)\big\|^2 & =\sum_{j=1}^m\sum_{l=1}^m\frac{1}{n}\sum_{t=1}^n\sum_{s=1}^n {\bf 1}^\top_j\Gamma^\ast_{t-s}{\bf 1}_l \langle\widehat{v}_j-v_j , \widehat{v}_l-v_l\rangle \\
& \leq \big(\sum_{j=1}^{m}\|\widehat{v}_j-v_j\|\big)^2\frac{1}{n}\sum_{t=1}^n\sum_{s=1}^n \|\Gamma^\ast_{t-s}\|_F\\
& =O_P\big(\big(n^{-1/2}\sum_{j=1}^{m} \alpha_j^{-1}\big)^2\big) \rightarrow 0.
\end{align*}
\hfill $ \Box$

{\bf Proof of Lemma 6.9:} \ Write $ L^+_{n,m}(\omega) =n^{-1/2}\sum_{t=1}^n W^+_t e^{-i t \omega}$ where $ W^+_{t} = \sum_{j=1}^m{\bf 1}_j^\top \xi_{t}^+ v_j$ with $ \xi_t^+=\sum_{l=0}^\infty \widetilde{\Psi}_{l,p}(m)e^+_{t-l}$, $\widetilde{\Psi}_{0,p}(m)=I_m$,    a random element in $ {\mathcal H}$.
Notice that  $ E^+(W^+_t)=0$,  while  using $ \xi_{t} = \sum_{l=0}^\infty \Psi_l(m) e_{t-l}$, $\Psi_{0}(m)=I_m$ , we get 
\begin{align*}
 E^+ W^+_t \otimes  W^+_{t+h} & = \sum_{l=0}^\infty\sum_{j=1}^m\sum_{s=1}^m {\bf 1}_j^\top \widetilde{\Psi}_{l,p}(m)\widetilde{\Sigma}_{e,p}(m)  \widetilde{\Psi}_{l+h,p}^\top(m) {\bf 1}_s \langle v_j,\cdot \rangle v_s\\
 & = \sum_{l=0}^\infty\sum_{j=1}^m\sum_{s=1}^m {\bf 1}_j^\top \Psi_{l}(m)\Sigma_{e}(m) \Psi_{l+h,p}^\top(m) {\bf 1}_s \langle v_j,\cdot \rangle v_s + \widetilde{D}_n\\
% & =  E\langle \sum_{j=1}^m \xi_{j,t} v_j, \cdot \rangle \sum_{s=1}^m\xi_{s,t+h} v_s + \widetilde{D}_n\\
 & = E\langle X_t-U_{t,m} ,\cdot  \rangle(X_{t+h}-U_{t+h,m})   + \widetilde{D}_n\\
 & = C_h(\cdot) - E\langle U_{t,m},\cdot  \rangle X_{t+h} - E\langle X_t , \cdot   \rangle U_{t+h,m}    \\
 & \ \ \ \ + E\langle U_{t,m} , \cdot  \rangle U_{t+h,m} + \widetilde{D}_n,
\end{align*}
with an obvious notation for $\widetilde{D}_n$. It is easily seen that  $ \widetilde{D}_n=O_P(\sum_{l=0}^\infty \|\widetilde{\Psi}_{l,p}(m) $ $-\Psi_{l}(m)\|_F + \| \widetilde{\Sigma}_{e,p}(m) -\Sigma_{e}(m)\|_F)$ and  therefore $ \|\widetilde{D}_n\|_{HS} \rightarrow 0$ in probability,   by Lemma 6.5. Hence and 
using  $ E\|U_{t,m}\|^2 \rightarrow 0$ as $ m \rightarrow \infty$, 
we get  that $ \|  E^+W^+_t\otimes W^+_{t+h} -C_h\|_{HS} \rightarrow 0$ in probability, as $ n \rightarrow \infty$.

Let $ \xi_{t}^o =   \sum_{l=0}^\infty \Psi_l(m) e^+_{t-l}$ and  define $ W^o_{t} = \sum_{j=1}^m{\bf 1}_j^\top \xi_{t}^o v_j$  and $ L^o_{n,m}(\omega) = n^{-1/2}\sum_{t=1}^n W_t^o e^{-it\omega}$. It easily follows by simple algebra and using   Lemma 6.5
that $ E^+\|L_{n,m}^+(\omega) - L_{n,m}^o(\omega)\|=O_P(\sum_{l=0}^\infty \|\widetilde{\Psi}_{l,p}(m)-\Psi_l(m) \|_F)\rightarrow 0$, in probability, that is 
 $ L_{n,m}^+(\omega) = L_{n,m}^o(\omega) + o_P(1)$. Thus to prove the assertion of the lemma it suffices to  show that $ L_{n,m}^o(\omega) \Rightarrow NC(0,2\pi {\mathcal F}_\omega)$. For this  we  show that Assumption 2  of Cerovecki and H\"ormann (2015) is satisfied, that is, using the notation  $ S_{n,m}^o(\omega)= \sum_{t=1}^n W^o_t e^{-it\omega}$ , we show that the following two conditions are fulfilled, in probability.
 \begin{equation} \label{th2.cond1}
 Z^o_n(\omega) \equiv \sum_{t=0}^n {\mathcal P}_0(W_t^o)e^{-it\omega}  \ \ \ \ 
 \mbox{is a Cauchy sequence in ${\mathcal H}$}, 
 \end{equation}
 and
 \begin{equation} \label{th2.cond2}
E\|E(S^o_{n,m}(\omega) | {\mathcal G}_0)\|^2 = o(n),
 \end{equation} 
where  the operator $ {\mathcal P}_0 $  is defined as  $ {\mathcal P}_0(\cdot) = E(\cdot | {\mathcal G}_0) - E(\cdot | {\mathcal G}_{-1})$ and $ {\mathcal G}_s=\sigma(W_s^o, W_{s-1}^o, W^o_{s-2}, \ldots)$.  Toward this we first define  
\[ W^o_{s,s} =\sum_{j=1}^m{\bf 1}_j^\top\sum_{l=0}^{\infty} \Psi_{l}(m)e^+_{s-l,s}v_j,\]
 where $ e^+_{t,s}=e^+_t$ if $ t >0$ and 
$ e^+_{t,s}=\widetilde{e}_{t}$ if $ t \leq 0$  with $ \widetilde{e}_{t}^+$  a copy of $ e^+_t$ which is independent of $ e^+_t $ for $ t<0$. We show that  
\begin{equation} \label{eq.uo}
\sum_{s=1}^\infty \sqrt{E^+\|W_0^o-W^o_{0,s}\|^2}  =O_P(1),
\end{equation}
where the $ O_P(1) $ term is independent of $m$ and $p$.  Notice first that by Minkowski's inequality
\begin{align*}
\sqrt{E^+\|W_0^o-W^o_{0,s}\|^2}  
% &= \sqrt{E^+\|\sum_{j=1}^{m}{\bf 1}_j^\top \sum_{l=s}^{\infty} \Psi_{l}(m) (e^+_{-l}-e^+_{-l,s})v_j\|^2 } \\
&  \leq \sqrt{E^+\| \sum_{j=1}^{m}{\bf 1}_j^\top \sum_{l=s}^{\infty} \Psi_{l}(m) e^+_{-l}v_j\|^2 } \\
& \ \ \ \ + \sqrt{E^+\|\sum_{j=1}^{m}{\bf 1}_j^\top \sum_{l=s}^{\infty} \Psi_{l}(m) e^+_{-l,s}v_j \|^2 }\\
%&= 2 \sqrt{\sum_{l=s}^\infty {\mbox tr}\big(\Psi_{l}(m) \widetilde{\Sigma}_{e,p}(m)\Psi^\top_{l}(m)\big)}\\
& \leq 2\|\widetilde{\Sigma}_{e,p}(m)\|_F\sum_{l=s}^{\infty}\|\Psi_l(m)\|_F.
\end{align*}
Thus
\begin{align*}
\sum_{s=1}^\infty \sqrt{E^+\|W_0^o-W^o_{0,s}\|^2} & \leq 2 \sum_{s=1}^{\infty} \|\widetilde{\Sigma}_{e,p}(m)\|_F\sum_{l=s}^{\infty}\|\Psi_l(m)\|_F\\
& \leq 2  \|\widetilde{\Sigma}_{e,p}(m)\|_F \sum_{l=1}^{\infty} l \|\Psi_l(m)\|_F.
\end{align*}
Now,  since  by  Lemma 6.1(ii),  $ \sum_{l=1}^{\infty} l \|\Psi_l(m)\|_F$ is bounded uniformly in $m$,  and,  
by Lemma 6.5(ii) and  (vi) and Lemma 6.1(iii),  $ \|\widetilde{\Sigma}_{e,p}(m)\|_F $ is bounded in probability,
where  the bound is independent of $p$ and $m$, assertion (\ref{eq.uo}) follows.  

Consider next  condition (\ref{th2.cond1}). For  positive integers  $ n_2 > n_1$ we have that 
\begin{align*}
E^+\|Z^o_{n_2}(\omega)-Z^o_{n_1}(\omega)\|^2 & \leq  \sum_{t_1=n_1+1}^{n_2}\sum_{t_2=n_1+1}^{n_2} |E^+\langle {\mathcal P}_0(W^o_{t_1}), {\mathcal P}_0(W^o_{t_1}) \rangle |\\
&\leq  \Big(\sum_{t=n_1+1}^{n_2} \sqrt{E^+\|{\mathcal P}_0(W^o_{t})\|^2}\Big)^2.
\end{align*}
Recall the definition of   $ W^o_{s,s}$. Then we have,  
% =\sum_{j=1}^m{\bf 1}_j^\top\sum_{l=0}^{\infty} \Psi_{l}(m)e^+_{t-l,s}v_j$, where $ e^+_{t,s}=e^+_t$ if $ t >0$ and 
%$ e^+_{t,s}=\widetilde{e}_{t}$ if $ t \leq 0$  with $ \widetilde{e}_{t}^+$ is a copy of $ e^+_t$ which is independent of $ e^+_t $ for $ %t<0$.
 since $E(W_{s,s}^o | {\mathcal G}_0) =E(W_{s,s}^o  | {\mathcal G}_{-1})=0$, that 
\begin{align*}
E^+\|{\mathcal P}_0(W^o_{t})\|^2 & = E^+\|{\mathcal P}_0(W_s^o) -{\mathcal P}(W^o_{s,s})\|^2\\
& =E^+\|E^+(W^o_s-W^o_{s,s} |{\mathcal G}_0 ) - E(W_s^o-W_{s,s}^o|{\mathcal G}_{-1}) \|^2 \\
& \leq 2E^+\|E^+(W^o_s-W^o_{s,s} |{\mathcal G}_0 )\|^2 + 2 E^+ \|E(W_s^o-W_{s,s}^o|{\mathcal G}_{-1}) \|^2\\
& \leq 4 E^+\| W^o_0-W^o_{0,s}\|^2
\end{align*} 
Hence 
\begin{align*}
E^+\|Z^o_{n_2}(\omega)-Z^o_{n_1}(\omega)\|^2 & \leq 4\Big( \sum_{s=n_1}^{\infty}\sqrt{E^+\|W_0^o-W^o_{0,s}\|^2}\Big)^2 \rightarrow 0,
\end{align*}
as $ n_1\rightarrow \infty$  because of (\ref{eq.uo}).

To establish  condition (\ref{th2.cond2}) notice that 
\begin{align*}
E\|E(S^o_{n,m}(\omega) | {\mathcal G}_0)\|^2 & 
\leq \sum_{t_1=1}^n\sum_{t_2=1}^n | E^+\langle E^+(W^o_{t_1}| {\mathcal G}_0),  E^+(W^o_{t_2}| {\mathcal G}_0)\rangle  |\\
& = \sum_{t_1=1}^n\sum_{t_2=1}^n | E^+\langle E^+(W^o_{t_1}-W_{t_1,t_1}^o| {\mathcal G}_0),  E^+(W^o_{t_2}-W_{t_2,t_2}^o| {\mathcal G}_0)\rangle  |\\
& \leq \Big( \sum_{t=1}^n \sqrt{E^+\|W^o_{0} - W^o_{0,t}\|^2}\Big)^2\\
& \leq \Big( \sum_{t=1}^\infty \sqrt{E^+\|W^o_{0} - W^o_{0,t}\|^2}\Big)^2,
 \end{align*}
 which is bounded because of (\ref{eq.uo}). \hfill $ \Box$

\section{Proof of   Theorem 4.2}  {~} \\ 

 In view of Theorem 4.1  and  Remark 4.1 of Paparoditis (2016), we get that  
$\sqrt{n}\overline{X}_n \Rightarrow N(0, C_X)$ and  $\sqrt{n}\overline{Y}_n \Rightarrow N(0, C_Y)$, where 
$ C_X=\sum_{h=-\infty}^{\infty} C_{h,X}$ and $ C_Y=\sum_{h=-\infty}^{\infty} C_{h,Y}$, with $ C_{h,X}$ and $ C_{h,Y}$ the  autocovariance operators at lag $h$ of  the processes $ {\bf X}$ and $ {\bf Y}$ respectively. Since $ {\bf X}^\ast$ amd $ {\bf Y}^\ast$ are independent we get, taking into account that $n_1/(n_1+n_2)\rightarrow \theta$,  the following convergence  on $ {\mathcal H}$ as $ n \rightarrow \infty$, 
\[ G_{n_1,n_2} =\sqrt{\frac{n_2n_1}{n_1+n_2}} \overline{X}_{n_1}^\ast + \sqrt{\frac{n_1n_2}{n_1+n_2}} \overline{Y}_{n_2}^\ast \Rightarrow N(0, (1-\theta)C_X + \theta C_Y).\]
By the continuous mapping theorem we then have $ U^\ast_{n_1,n_2} =\|G_{n_1,n_2}\|^2 \Rightarrow \int_0^1\Gamma^2(\tau)\tau$, where $\{\Gamma(\tau), \tau \in [0,1]\} $ is  a Gaussian process in $ {\mathcal H}$ with mean zero and covariance $ Cov(\Gamma(\tau_1),\Gamma(\tau_2))=(1-\theta)  c_X(\tau_1,\tau_2) + \theta c_Y(\tau_1,\tau_2)$,  $\tau_1, \tau_2\in [0,1]$, where $ c_X$ and $c_Y$  denote the covariance kernels of the operators $ C_X $ and $ C_Y $, respectively.
 \hfill $\Box$
  
\section{Additional numerical results}  
\subsection{Choice of starting values} 
To generate the $m$-dimensional time series of pseudo scores $ \xi_1^\ast, \xi_2^\ast, \ldots, \xi_n^\ast$, a set of $p$ starting values have to be chosen. Different  alternatives  can be used. In order to   obtain a  time series of length $n$, we generated time series of length $ n+L$ using as starting values the observed values $ \widehat{\xi}_1, \widehat{\xi}_2, \ldots, \widehat{\xi}_p$  and then discarded  the first $L$ observations to eliminate the effects of these starting values.  The number $L$ has  be chosen  adapting to the multivariate case a proposal made for the univariate case by McLeod and Hipel (1978). To elaborate, we first calculated  $ \widehat{\Gamma}_{\xi}(0)$ given by   
\[ \widehat{\Gamma}_\xi(0) = \sum_{j=0}^{\infty} \widehat{\Psi}_{j,p}(m) \widehat{\Sigma}_e(m) \widehat{\Psi}^\top_{j,p}(m)=\int_{-\pi}^\pi f_{\widehat{\xi}}(\omega)d\omega, \ \ \ \ \widehat{\Psi}_{0,p}(m)=I_m,\]
where $f_{\widehat{\xi}}(\omega) $ denotes  the spectral density of the  VAR(p) model fitted to the $m$-dimensional time series of estimated scores.    
We then selected  a natural number  $S$ such that 
\[ \|\widehat{\Gamma}_\xi(0)-\widetilde{\Gamma}_\xi(0)\|_F < \delta, \]
where  
\[ \widetilde{\Gamma}_\xi(0) = \sum_{j=0}^S \widehat{\Psi}_{j,p}(m) \widehat{\Sigma}_e(m) \widehat{\Psi}^\top_{j,p}(m)\]
and  $\delta $ has been set equal to a very small number, i.e., $ \delta=10^{-5}$.  This essentially implies that observations $ X_{t-j}$ for  $j \geq S$   have practically no effect on the current value $X_t$. 
For instance, for  the model (5.1)  used in the simulations we found in a number of 20 preliminary  runs,    that  the values of $S$ obtained (which depend on the estimates $ \widehat{A}_{j,p}(m)$ and $ \widehat{\Sigma}_e(m)$), vary  between $10$ and $18$.  To be on the safe side, we have  
set    for this model $L=30$ to eliminate the effects of the starting values $ \widehat{\xi}_1, \widehat{\xi}_2, \ldots, \widehat{\xi}_p$. 

\subsection{Additional simulations for the FMA(1) model (5.1)}
Table 1 shows the results obtained for selecting the number $m$ of principal components  according to the rule $ \widehat{m}_n=\max\{m_{n,Q}, m_{n,E}\} $, $Q=0.85$ and for different sample sizes.  Note that for $ n \leq200$ the VR  while for $n>200$ the GVR  criterion is used to calculate $m_{n,Q}$.
 \begin{table*}
\caption{Frequency  of selected values of $m_{n,Q}$, $ m_{n,E}$ and of $ \widehat{m}_n$ ($ R=1000$ replications).} 
\begin{tabular}{llllllllll}
\hline
%& & &   &  & & &  & & \\
 &&m =&  1 & 2 & 3& 4& 5& 6& 7\\
\hline
n=100 & $ m_{n,Q}    $  & & 0 & 0.3 & 67.1 & 32.6 & 0 & 0&0  \\
           &  $m_{n,E} $  & &20.1 & 79.3 & 0.6 & 0 & 0 & 0 & 0   \\
           &  $\widehat{m}_n $  & & 0 & 0.3 & 67.1 & 32.6 & 0 & 0 & 0 \\
 &&&   &  & & & & & \\
 n=500 &   $m_{n,Q} $        &  & 0 & 0  & 0.9 & 83.0 & 16.1 & 0 & 0 \\
           &  $m_{n,E}  $          &  & 0 & 77.3 & 22.7 & 0 & 0 & 0 &  0 \\
           &  $\widehat{m}_n$  &  & 0  & 0 & 0.9 & 83.0 & 16.1& 0  &  0\\
 &&&   &  & & & & & \\ 
  n=1000 &  $ m_{n,Q} $      &   & 0 & 0  & 0.2 & 89.3 & 10.5 & 0 & 0 \\
           &  $m_{n,E} $            &  & 0 & 1.8  & 98.2 & 0 & 0 & 0 &  0 \\ 
           &  $\widehat{m}_n $  & &  0  & 0  & 0.2 & 89.3 & 10.5 & 0 & 0 \\
            &&&   &  & & & & & \\
  n=5000 &  $m_{n,Q} $       &   & 0 & 0  & 0 & 99.4 & 0.6 & 0 & 0 \\
           &  $m_{n,E} $           &    & 0 & 0 & 0 & 100.0 & 0 & 0 &  0 \\ 
           &  $\widehat{m}_n$  & &  0 & 0  & 0 &  99.4 & 0.6 &  0 & 0 \\
            &&&   &  & & & & & \\
   n=10000 &  $ m_{n,Q} $   &   & 0 & 0  & 0 &  99.8& 0.2 & 0 & 0 \\
           &  $m_{n,E} $           &   & 0 & 0 & 0 & 50.8 & 49.2 & 0 &  0 \\ 
           &  $\widehat{m}_n$  &  & 0  & 0 & 0 & 50.8 &  49.2 & 0 & 0 \\
            &&&   &  & & & & & \\
    n=20000 &  $ m_{n,Q} $   &   & 0 & 0  & 0& 100.0 & 0 & 0 & 0 \\
           &  $m_{n,E} $            &  & 0 & 0  & 0 & 0  & 99.5 & 0.5 &  0 \\
           &  $\widehat{m}_n$   & & 0 & 0 &  0 & 0 &  99.5 & 0.5 & 0\\ 
           \hline                 
\end{tabular}
\end{table*}

Table 2 shows the FSB estimates obtained using  some different values of the bootstrap parameters  $m$ and $p$ as well as for the values of  these parameters chosen by means of the $\widehat{m}_n$ and $AICC$  rule and which are denoted by $ (\widehat{m},\widehat{p})$.
Note that  $(m,p)=(3,3)$ is the most frequently chosen pair using this data driven selection rule.  
%As this table  shows the FSB estimates  are  quite good even for  the short  functional time series  of $n=100$ observations. These estimates   
%seem  also  not to be very sensitive with respect to the different choices of the parameter  $m$ used to truncate  the Karhunen-Lo\'eve expansion.  

\begin{table*}
\caption{Estimated exact ($ \sigma_{EE}(\tau_j)$)  and functional sieve bootstrap (FSB) estimates of the standard deviation of the sample mean  $ \overline{X}_n(\tau_j)$    for different values 
of $ \tau_j \in [0,1]$   and for different parameters $ m$ and $p$.   $\overline{\widehat{\sigma}}(\tau_j)$ refers to the mean,  while $S(\widehat{\sigma}(\tau_j) ) $ to the standard deviation  of the  FSB estimates.
% obtained for  $ R=1000$ replications
} 
\begin{tabular}{lcllllll}
\hline
           & & \multicolumn{2}{l}{m=2, p=3}  &  \multicolumn{2}{l}{m=3,  p=3}  &  \multicolumn{2}{l}{$\widehat{m}, \widehat{p}$}    \\
$\tau_j$ & $ \sigma_{EE}(\tau_j)$   & $\overline{\widehat{\sigma}}(\tau_j)$  & $S(\widehat{\sigma}(\tau_j) )$&  $\overline{\widehat{\sigma}}(\tau_j)$  & $S(\widehat{\sigma}(\tau_j) )$ &  $\overline{\widehat{\sigma}}(\tau_j)$  & $S(\widehat{\sigma}(\tau_j) )$   \\
\hline
      0.00 &  2.149 &  2.124 & 0.392 &  2.188 & 0.440 &  2.025  & 0.462   \\
    0.05 &  2.203 &  2.172 & 0.404 &  2.227 & 0.440 &  2.072  & 0.473   \\
    0.10 &  2.272 &  2.262 & 0.441 &  2.305 & 0.458 &  2.141  & 0.480   \\
    0.15 &  2.325 &  2.362 & 0.466 & 2.385& 0.477&  2.196  & 0.501  \\
    0.20 &  2.358 &  2.429 & 0.484 &  2.434 & 0.492 &  2.227  & 0.510   \\
    0.25 &  2.370 &  2.457 & 0.488 &  2.452 & 0.488 &  2.240  & 0.516   \\
    0.30 &  2.351 &  2.429 & 0.488 &  2.432 & 0.485 &  2.231  & 0.509   \\
    0.35 &  2.317 &  2.359 & 0.462 &  2.382 & 0.471 &  2.203  & 0.493   \\
    0.40 &  2.267 &  2.271 & 0.435 &  2.307 & 0.448 &  2.138  & 0.470  \\
    0.45 &  2.196 &  2.183 & 0.419 &  2.237 & 0.439 &  2.062  & 0.452   \\
    0.50 &  2.146 &  2.123 & 0.401 &  2.199 & 0.433 &  2.026  & 0.446   \\
    0.55 &  2.194 &  2.165 & 0.405 &  2.240 & 0.440 &  2.075  & 0.456   \\
    0.60 &  2.264 &  2.249 & 0.419 &  2.309 & 0.459 &  2.148  & 0.473   \\
    0.65 &  2.314 &  2.342 & 0.441 &  2.370 & 0.468 &  2.204  & 0.490   \\
    0.70 &  2.343 &  2.408 & 0.464 &  2.418 & 0.487 &  2.241  & 0.505   \\
    0.75 &  2.351 &  2.429 & 0.475 &  2.430 & 0.494 & 2.244  & 0.513  \\
    0.80 &  2.342 &  2.405 & 0.474 &  2.413 & 0.481 &  2.235  & 0.510   \\
    0.85 &  2.309 &  2.346 & 0.459 &  2.364 & 0.473 &  2.198  & 0.497   \\
    0.90 &  2.258 &  2.262 & 0.431 &  2.299 & 0.456 &  2.133  & 0.482   \\
    0.95 &  2.188 &  2.167 & 0.399 &  2.227 & 0.444 &  2.061  & 0.463   \\
    1.00 &  2.149 &  2.123 & 0.392 &  2.188 & 0.440 &  2.025  & 0.462   \\
\hline
\end{tabular}
\end{table*}

\subsection{Size and power behavior of the bootstrap based test for the two-sample problem}
Note that Theorem 4.2  justifies  the  use of percentage points of the  distribution of $ U_{n_1,n_2}^\ast$ in order to obtain bootstrap critical values of  the test $U_{n_1,n_2}$. Furthermore, if $ H_1$ is true, that is if $ \|\mu_X-\mu_Y\| >0$  and $ U_{n_1,n_2} \stackrel{p}{\rightarrow} \infty$ as $ n_1,n_2 \rightarrow \infty$, see for instance  Theorem 4 of Horv\'ath et al. (2013),  then  the  consistency 
 of the  test $U_{n_1,n_2}$ based  on sieve bootstrap estimated critical values,  follows.    
 
To investigate the size and power behavior of the bootstrap based, fully functional test $ U_{n_1,n_2}$, we conducted a small numerical  experiment 
by adopting  the simulation design of  Horv\'ath et al. (2013) and considering  the functional moving average model 
\[ X_{t} =\Theta_1(\varepsilon_{t-1}) + \varepsilon_t,\]
with $\Theta_1$ the  integral operator  with kernel 
$$ \theta_1(t,s)=\frac{\displaystyle \mbox{exp}\{-(t^2+s^2)/2\}}{\displaystyle 4\int \mbox{exp}(-x^2)dx }$$  
and $\{\varepsilon_t\}$ i.i.d.  Brownian bridges. Pairs of functional time series of length $n_1$  and   $n_2$  have  been generated using the above FMA(1)  model with mean functions given by $ \mu_1=0$ for the first and $ \mu_2(\tau)=\gamma \tau(1-\tau)$, $\tau \in [0,1]$, for the second time series.   Notice that the value $\gamma=0$ corresponds to the null hypothesis while the degree of deviation from the null under the alternative is controlled by the parameter $\gamma$.   The rejection frequencies  obtained for different  sample sizes based on $R=200$ repetitions and $B=1000$ bootstrap replications are reported in Table 3 for different choices of the parameters $m$ and $p$ and for three different nominal levels. Notice that 
the data driven values of $m$ and $p$ chosen using the $\widehat{m}_n$ and $AICC$ rule are denoted in this table by $(\widehat{m},\widehat{p})$,  while 
$ (m,p)=(3,1)$ and $(m,p)=(3,2)$  are the most frequently chosen values of the corresponding parameters using the same rule  for $n_1=n_2=100$  and for $ n_1=n_2=200$, respectively. As this table shows, using the critical values 
 obtained by means of the functional sieve bootstrap procedure, the fully functional test $U_{n_1,n_2}$  retains the nominal  size  and shows at the same time a  nice  power behavior; the  power of the test increases as the deviation from  the null   and/or  the sample size increases.

%\begin{table*}
%\caption{Size and power behavior of the FSB-based test  for  the two-sample mean problem 
%($ R=1000$ replications).
%} 
%\begin{tabular}{llcclcc}
%\hline
%                   &  \multicolumn{3}{l}{$n_1=100$}  &  \multicolumn{3}{l}{$n_1=200$}    \\
%                    &  \multicolumn{3}{l}{$n_2=100$}  &  \multicolumn{3}{l}{$n_2=200$}    \\
%$\gamma$ & (p,m)   & $\alpha=0.05 $  & $\alpha=0.10$  &   (p,m) &$\alpha=0.05 $  & $\alpha=0.10$       \\
%\hline
%      0 &  (1,3) & 0 & 0 & (2,3)& 0 &  0     \\
%        &   $(\widehat{p},\widehat{m})$ & 0 & 0 &  $(\widehat{p},\widehat{m})$ & 0 &  0     \\ 
%        0.2 &  (1,3) & 0 & 0 & (2,3)& 0 &  0     \\
%        &   $(\widehat{p},\widehat{m})$ & 0 & 0 &  $(\widehat{p},\widehat{m})$ & 0 &  0     \\      
%  0.5 &  (1,3) & 0 & 0 & (2,3)& 0 &  0     \\
%        &   $(\widehat{p},\widehat{m})$ & 0 & 0 &  $(\widehat{p},\widehat{m})$ & 0 &  0     \\
%  0.8 &  (1,3) & 0 & 0 & (2,3)& 0 &  0     \\
%        &   $(\widehat{p},\widehat{m})$ & 0 & 0 &  $(\widehat{p},\widehat{m})$ & 0 &  0     \\ 
%   1.0&  (1,3) & 0 & 0 & (2,3)& 0 &  0     \\
%        &   $(\widehat{p},\widehat{m})$ & 0 & 0 &  $(\widehat{p},\widehat{m})$ & 0 &  0     \\                         
% \hline
%\end{tabular}
%\end{table*}

\begin{table*}
\caption{Size and power behavior of the FSB-based test  for  the two-sample mean problem (R=200 replications,  B=1000 bootstrap 
samples). 
%($ R=1000$ replications).
} 
\begin{tabular}{llccclccc}
\hline
                   &  \multicolumn{4}{l}{$n_1=100$}  &  \multicolumn{4}{l}{$n_1=200$}    \\
                    &  \multicolumn{4}{l}{$n_2=100$}  &  \multicolumn{4}{l}{$n_2=200$}    \\
                         &   &\multicolumn{3}{l}{$\alpha=$}  &  & \multicolumn{3}{l}{$\alpha=$}    \\
$\gamma$ & (m,p)   &   0.01   & 0.05  & 0.10 &  (m,p) &0.01   & 0.05 & 0.10     \\
\hline
      0 &  (3,1) & 0.008 & 0.055 & 0.125 & (3,2) & 0.010 &  0.050 &   0.112 \\
        &   $(\widehat{m},\widehat{p})$ & 0.010 & 0.050& 0.095&  $(\widehat{m},\widehat{p})$ & 0.015 &  0.045  &0.087   \\ 
        0.2 & (3,1) & 0.018 & 0.085 & 0.170 & (3,2)& 0.055 &  0.135 &0.210      \\
        &   $(\widehat{m},\widehat{p})$ & 0.035 & 0.080 & 0.180 & $(\widehat{m},\widehat{p})$ & 0.045 &  0.150  & 0.245  \\      
  0.5 &  (3,1) & 0.180 & 0.325 &  0.455 & (3,2)& 0.435 &  0.635 &0.770     \\
        &   $(\widehat{m},\widehat{p})$& 0.215 & 0.455 & 0.575  &  $(\widehat{m},\widehat{p})$ & 0.355 &  0.575    & 0.715\\
  0.8 &  (3,1) & 0.535 & 0.790 & 0.870& (3,2)& 0.915 &  0.955 &0.980      \\
        &   $(\widehat{m},\widehat{p})$ & 0.495 & 0.690 &0.815  & $(\widehat{m},\widehat{p})$ & 0.865 &  0.960   &  0.985\\ 
   1.0& (3,1) & 0.715 & 0.880 &0.940&  (3,2)& 0.980 &  1.000 &1.000     \\
        &   $(\widehat{m},\widehat{p})$ & 0.735 & 0.835 &  0.930 & $(\widehat{m},\widehat{p})$& 0.985 &  1.000   &  1.000 \\                         
 \hline
\end{tabular}
\end{table*}

\end{document}